\documentclass[11pt]{article}
\usepackage[margin = 1in]{geometry}

\usepackage{lipsum}
\usepackage{amsfonts}
\usepackage{graphicx}
\usepackage{epstopdf}

\usepackage{mathtools,mathrsfs,mathabx}
\usepackage{empheq}
\usepackage{amsfonts}
\usepackage{amsgen,amsthm,amsmath,amstext,amsbsy,amsopn,amssymb,subfigure,stmaryrd,bbm,dsfont}
\usepackage{comment}
\usepackage{enumitem}
\usepackage{natbib}
\usepackage{booktabs}
\usepackage{blkarray}
\usepackage{array}
\usepackage{algorithm}
\usepackage{algpseudocode}
\usepackage{url}
\usepackage{longtable}
\usepackage{mathtools}
\usepackage{multirow}
\usepackage{color}
\usepackage{hyperref}

\ifpdf
  \DeclareGraphicsExtensions{.eps,.pdf,.png,.jpg}
\else
  \DeclareGraphicsExtensions{.eps}
\fi

\newtheorem{Definition}{Definition}
\newtheorem{Theorem}{Theorem}

\newtheorem{Lemma}{Lemma}
\newtheorem{Remark}{Remark}
\newtheorem{Corollary}{Corollary}

\newtheorem{Proposition}{Proposition}

\DeclareMathAlphabet\mathbfcal{OMS}{cmsy}{b}{n}

\newcommand{\be}{\begin{equation}}
\newcommand{\ee}{\end{equation}}
\newcommand{\bea}{\begin{eqnarray}}
\newcommand{\eea}{\end{eqnarray}}
\newcommand{\beas}{\begin{eqnarray*}}
	\newcommand{\eeas}{\end{eqnarray*}}

\newcommand{\bbR}{\mathbb{R}}
\newcommand{\bbS}{\mathbb{S}}

\newcommand{\cU}{\mathcal{U}}
\newcommand{\cV}{\mathcal{V}}
\newcommand{\cL}{\mathcal{L}}
\newcommand{\cP}{\mathcal{P}}
\newcommand{\cJ}{\mathcal{J}}

\newcommand{\0}{{\mathbf{0}}}

\renewcommand{\a}{{\mathbf{a}}}

\renewcommand{\u}{{\mathbf{u}}}
\renewcommand{\v}{{\mathbf{v}}}

\newcommand{\x}{{\mathbf{x}}}
\newcommand{\y}{{\mathbf{y}}}

\newcommand{\X}{{\mathbf{X}}}
\newcommand{\B}{{\mathbf{B}}}

\newcommand{\D}{{\mathbf{D}}}

\newcommand{\I}{{\mathbf{I}}}

\renewcommand{\P}{{\mathbf{P}}}

\newcommand{\Q}{{\mathbf{Q}}}
\newcommand{\R}{{\mathbf{R}}}
\renewcommand{\L}{{\mathbf{L}}}
\renewcommand{\S}{{\mathbf{S}}}
\newcommand{\U}{{\mathbf{U}}}
\newcommand{\V}{{\mathbf{V}}}
\newcommand{\W}{{\mathbf{W}}}

\newcommand{\A}{{\mathbf{A}}}
\newcommand{\Y}{{\mathbf{Y}}}
\newcommand{\Z}{{\mathbf{Z}}}

\newcommand{\G}{{\mathbf{G}}}

\newcommand{\cM}{{\cal M}}

\newcommand{\scA}{{\mathscr{A}}}
\newcommand{\scC}{{\mathscr{C}}}
\newcommand{\scS}{{\mathscr{S}}}

\newcommand{\bSigma}{\boldsymbol{\Sigma}}
\newcommand{\bDelta}{\boldsymbol{\Delta}}

\newcommand{\rank}{{\rm rank}}

\newcommand{\Hess}{{\rm Hess}}
\newcommand{\grad}{{\rm grad}}
\newcommand{\reg}{{\rm reg}}
\newcommand{\rnull}{{\rm null}}
\newcommand{\F}{{\rm F}}
\newcommand{\rspan}{{\rm span}}

\newcommand{\bbO}{\mathbb{O}}

\makeatletter
\newcommand*{\rom}[1]{\expandafter\@slowromancap\romannumeral #1@}
\makeatother

\allowdisplaybreaks

\graphicspath{ {images/} }

\begin{document}
	\title{Nonconvex Factorization and Manifold Formulations are Almost Equivalent in Low-rank Matrix Optimization}
	
	\author{Yuetian Luo$^1$, ~ Xudong Li$^2$, ~ and ~ Anru R. Zhang$^{3}$}
	
	\date{}
	\maketitle

	\footnotetext[1]{Data Science Institute, University of Chicago.  (\texttt{yuetian@uchicago.edu})}
	\footnotetext[2]{School of Data Science, Fudan University. (\texttt{lixudong@fudan.edu.cn})}
	\footnotetext[3]{Departments of Biostatistics \& Bioinformatics, Computer Science, Electrical \& Computer Engineering, Mathematics, and Statistical Science, Duke University. (\texttt{anru.zhang@duke.edu})}
	
	\bigskip

\begin{abstract}
In this paper, we consider the geometric landscape connection of the widely studied manifold and factorization formulations in low-rank positive semidefinite (PSD) and general matrix optimization. We establish a sandwich relation on the spectrum of Riemannian and Euclidean Hessians at first-order stationary points (FOSPs). As a result of that, we obtain an equivalence on the set of FOSPs, second-order stationary points (SOSPs), and strict saddles between the manifold and the factorization formulations. In addition, we show the sandwich relation can be used to transfer more quantitative geometric properties from one formulation to another. Similarities and differences in the landscape connection under the PSD case and the general case are discussed. To the best of our knowledge, this is the first geometric landscape connection between the manifold and the factorization formulations for handling rank constraints, and it provides a geometric explanation for the similar empirical performance of factorization and manifold approaches in low-rank matrix optimization observed in the literature. In the general low-rank matrix optimization, the landscape connection of two factorization formulations (unregularized and regularized ones) is also provided. By applying these geometric landscape connections, in particular, the sandwich relation, we are able to solve unanswered questions in literature and establish stronger results in the applications on geometric analysis of phase retrieval, well-conditioned low-rank matrix optimization, and the role of regularization in factorization arising from machine learning and signal processing.
\end{abstract}

\section{Introduction}

Low-rank optimization problems are ubiquitous in a variety of fields, such as optimization, machine learning, signal processing, scientific computation, and statistics. One popular formulation is the following rank constrained optimization: \begin{equation} \label{eq: PSD-manifold-formulation}
	\text{PSD case}: \quad \quad  \min_{\X \in \bbS^{p \times p} \succcurlyeq 0, \rank(\X) = r} f(\X), \quad 0 < r \leq p,
	\end{equation}
\begin{equation} \label{eq: general prob}
	\text{general case}: \quad \min_{\X \in \bbR^{p_1 \times p_2}, \rank(\X) = r} f(\X), \quad 0 < r \leq \min\{p_1,p_2\}.
\end{equation} 
In the positive semi-definite (PSD) case, without loss of generality, we assume $f$ is symmetric in $\X$, i.e., $f(\X) = f(\X^\top)$; otherwise, we can set $\tilde{f}(\X) = \frac{1}{2}(f(\X) + f(\X^\top)) $ and have $\tilde{f}(\X) = f(\X)$ for all $\X \succcurlyeq 0 $ \citep{bhojanapalli2016dropping}. In both cases, we assume $f$ is twice continuously differentiable with respect to $\X$ and the Euclidean metric. 
Viewed as optimization problems over low-rank matrix manifolds under the embedded geometry \citep{absil2009optimization,boumal2020introduction},  \eqref{eq: PSD-manifold-formulation} and \eqref{eq: general prob} can be solved via various manifold optimization methods. On the other hand, to accelerate the computation and to better cope with the rank constraint, a line of research studied the following nonconvex factorization formulation \citep{burer2005local}:
\begin{equation} \label{eq: PSD factorization}
	\text{PSD case}: \quad \quad	\min_{\Y \in \bbR^{p \times r} } g(\Y) :=  f(\Y \Y^\top),
\end{equation}
\begin{equation} \label{eq: general factor formu}
	\text{general case}: 	\min_{\L \in \bbR^{p_1 \times r}, \R\in \bbR^{p_2 \times r} } g(\L, \R):= f(\L \R^\top).
\end{equation}
In the general asymmetric case, to promote balance between two factors $\L$ and $\R$ in \eqref{eq: general factor formu}, the following regularized optimization problem has also been widely studied \citep{tu2016low}:
\begin{equation} \label{eq: general factor with reg}
		\min_{\L \in \bbR^{p_1 \times r}, \R\in \bbR^{p_2 \times r} } g_{\reg}(\L, \R):=f(\L \R^\top) + \frac{\mu}{2} \|\L^\top \L - \R^\top \R\|_\F^2,
\end{equation}
where $\mu >0$ is some properly chosen regularization parameter.
Note that \eqref{eq: PSD factorization}, \eqref{eq: general factor formu}, and \eqref{eq: general factor with reg} are unconstrained, and thus can be tackled by running unconstrained optimization algorithms. Indeed, under proper assumptions, a number of algorithms with theoretical guarantees have been proposed for both the manifold and the factorization formulations   \citep{chi2019nonconvex,cai2018exploiting}. See Section \ref{sec: related-literature} for a review of existing results. 

On the other hand, the manifold and the factorization formulations are more or less treated as two different approaches for low-rank matrix optimization in the literature and they are not obviously related. Similar algorithmic guarantees under these two formulations, including convergence rate and sample complexity for successful recovery, were observed in a number of matrix inverse problems \citep{wei2016guarantees,luo2020recursive,cai2018solving,zhang2018robust,keshavan2009matrix,ma2019implicit,cai2015rop,chen2015fast,hardt2014understanding,zhao2015nonconvex,zheng2015convergent,wang2017unified,tong2020accelerating}, while there are little studies on the reason behind. Moreover, most of the existing geometric analyses in low-rank matrix optimization are performed under the factorization formulation \citep{bhojanapalli2016global,ge2017no,zhang2019sharp,zhu2018global,zhu2017global,li2019non,park2017non}. It has been asked by \cite{cai2018solving,li2019toward} whether it is possible to investigate the geometric landscape directly on the low-rank matrix manifolds as the manifold formulation avoids unidentifiable parameterizations of low-rank matrices and explicit regularizations to cope with the unbalanced factorization in \eqref{eq: general factor with reg}. In this work, we make the first attempt to answer these questions by investigating the geometric landscape connections between the manifold and the factorization formulations in low-rank matrix optimization.

\subsection{Our Contributions} \label{sec: contribution}
{First,  we establish a sandwich relation on the spectrum between Riemannian and Euclidean Hessians at  FOSPs under the manifold and the factorization formulations. In particular, sandwich inequalities between Riemannian and Euclidean Hessians are established for \eqref{eq: PSD-manifold-formulation}, \eqref{eq: PSD factorization} and \eqref{eq: general prob}, \eqref{eq: general factor formu}; a partial sandwich inequality is built between \eqref{eq: general prob} and \eqref{eq: general factor with reg}. As an immediate corollary, we obtain an equivalence on the set of first-order stationary points (FOSPs), second-order stationary points (SOSPs) and strict saddles between the manifold formulation under the embedded geometry and the factorization formulation in both the PSD and the general low-rank matrix optimization. In addition, we demonstrate the sandwich relation is useful in transferring more geometric landscape properties, such as the strict saddle property, from one formulation to another. To the best of our knowledge, this is the first equivalence geometric landscape connection between the manifold and the factorization formulations for low-rank matrix optimization. Key technical ingredients to establish these results include a characterization of the zero eigenspace of the Hessian of the factorization objective and a bijection between its orthogonal complement and the tangent space of the fixed-rank $r$ manifold at the reference point. In addition, a few similarities and key differences in the landscape connection under the PSD case and the general case are identified.}

We also provide a geometric landscape connection between the unregularized and the regularized factorization formulations (\eqref{eq: general factor formu} and \eqref{eq: general factor with reg}) and give a sandwich inequality on the spectrum of Euclidean Hessians of two factorization formulations at rank $r$ FOSPs.

Furthermore, we apply our main results in three applications from machine learning and signal processing. By our geometric landscape connections between the manifold and the factorization formulations, we provide the first global optimality result for phase retrieval with a rate-optimal sample complexity under the manifold formulation and specifically show there is a unique Riemannian SOSP that is the global optima and all other Riemannian FOSPs are strict saddles with an explicit upper bound on the negative eigenvalue. We also prove the global optimality result for generic well-conditioned low-rank matrix optimization under the manifold formulation in both exact-parameterization and over-parameterization settings. Finally, we provide a geometric analysis on the role of regularization in the factorization formulation for a general $f$; when $f$ is further well-conditioned, we give a global optimality result under the formulation \eqref{eq: general factor formu}. {All of these results rely critically on the sandwich inequalities we establish between the Riemannian and Euclidean Hessians under the manifold and the factorization formulations.}

In a broad sense, manifold and factorization can be treated as two different approaches in handling the rank constraint in optimization problems. This paper bridges them from a geometric point of view and demonstrates that the manifold and the factorization approaches are indeed strongly connected in solving low-rank matrix optimization problems.

\subsection{Related Literature} \label{sec: related-literature}
This work is related to a range of literature on low-rank matrix optimization, manifold/nonconvex optimization, and geometric landscape analysis arising from a number of communities, such as optimization, machine learning and signal processing.

First, from an algorithmic perspective, a number of algorithms, including the penalty approaches, gradient descent, alternating minimization, and Gauss-Newton, have been developed either for solving the manifold formulation \citep{bi2020multistage,gao2010majorized,boumal2011rtrmc,mishra2014fixed,meyer2011linear,mishra2014fixed,vandereycken2013low,huang2018blind,luo2020recursive} or the factorization formulation \citep{candes2015phase,jain2013low,sun2015guaranteed,tran2016extended,tu2016low,wen2012solving,bauch2020rank}. We refer readers to \cite{chi2019nonconvex,cai2018exploiting} for the recent algorithmic development under two formulations. Many algorithms developed under the manifold formulation involve Riemannian optimization techniques and can be more complex than the ones developed under the factorization formulation. On the other hand, similar guarantees were observed for both lines of algorithms under two formulations in various matrix inverse problems \citep{miao2016rank,wei2016guarantees,luo2020recursive,hou2020fast,cai2018solving,zhang2018robust,keshavan2009matrix,bhojanapalli2016dropping,park2018finding,li2019rapid,ma2019implicit,sanghavi2017local,chen2015fast,hardt2014understanding,zhao2015nonconvex,zheng2015convergent,wang2017unified,tong2020accelerating}. Our results on the geometric connection between two formulations shed light on this phenomenon by showing that these two approaches are in essence closely related. 

Second, from a geometric landscape perspective, a body of work showed that factorization will not introduce spurious local minima compared to the original rank constrained optimization problem when the objective $f$ is well-conditioned \citep{bhojanapalli2016global,ge2017no,zhang2019sharp,zhu2018global,chen2019model,park2017non,zhang2018primal,zhang2021general}. Similar benign landscape results were proved for factorization in solving semidefinite programs and convex programs on PSD matrices or with a nuclear norm regularization \citep{boumal2020deterministic,journee2010low,yamakawa2021equivalent,li2019non}. On the other hand, it is much less explored for the geometric analysis under the manifold formulation. \cite{maunu2019well,ahn2021riemannian} provided landscape analyses for robust subspace recovery and matrix factorization over the Grassmannian manifold. Under the embedded manifold, \cite{uschmajew2018critical} showed the benign landscape of \eqref{eq: general prob} when $f$ is quadratic and satisfies certain restricted spectral bounds properties. Different from both lines of work focusing on the landscape under either the factorization or the manifold formulation when $f$ is well-conditioned, i.e., $f$ satisfies the restricted strong convexity and smoothness or restricted spectral bounds properties, here we study the geometric landscape connection between the factorization and the manifold formulations in low-rank matrix optimization for a general $f$.

The closest work in the literature related to ours is \cite{ha2020equivalence}, where they study the relationship between Euclidean FOSPs and SOSPs under the factorization formulation and fixed points of the projected gradient descent (PGD) in the general low-rank matrix optimization. They show while the sets of FOSPs of \eqref{eq: general factor formu} and \eqref{eq: general factor with reg} can be larger, the sets of SOSPs of these two factorization formulations are contained in the set of fixed points of the PGD with a small stepsize. Complementary to their results, here we consider the geometric landscape connection between the manifold and the factorization formulations under both the PSD and the general low-rank matrix optimization and establish a stronger equivalence on sets of FOSPs as well as SOSPs of manifold and factorization formulations.

\subsection{Organization of the Paper}\label{sec: organization}
The rest of this article is organized as follows. After a brief introduction of notation, we introduce Riemannian optimization and some preliminary results on the Riemannian geometry of low-rank matrices in Section \ref{sec: Riemannian-opt-background}. Our main results on the geometric landscape connection between the manifold and the factorization formulations in low-rank PSD and general matrix optimization are presented in Sections \ref{sec: connection-PSD} and \ref{sec: connection-general}, respectively. In Section \ref{sec: application}, we present three applications of our main results in machine learning and signal processing. Conclusion and future work are given in Section \ref{sec: conclusion}. We present the proofs of the main results in the main text and additional proofs and lemmas are presented in Appendices \ref{sec: additional-proofs} and \ref{sec: additional-lemmas}, respectively.

\section{Notation and Preliminaries} \label{sec: notation}

The following notation will be used throughout this article. $\bbR^{p_1 \times p_2}$ and $\bbS^{p \times p}$ denote the spaces of $p_1$-by-$p_2$ real matrices and $p$-by-$p$ real symmetric matrices, respectively. Uppercase and lowercase letters (e.g., $A, B, a, b$), lowercase boldface letters (e.g. $\u, \v$), uppercase boldface letters (e.g., $\U, \V$) are used to denote scalars, vectors, matrices, respectively. We denote $[p_k]$ as the set $\{1,\ldots, p_k\}$. For any $a, b \in \bbR$, let $a \wedge b := \min\{a,b\}, a \vee b := \max\{a,b\}$. For any vector $\v$, denote its $\ell_1$ and $\ell_2$ norms as $\|\v\|_1$ and $\|\v\|_2$, respectively. For any matrix $\X \in \mathbb{R}^{p_1\times p_2}$ with singular value decomposition (SVD) $\sum_{i=1}^{p_1 \land p_2} \sigma_i(\X)\u_i \v_i^\top$, where $\sigma_1(\X) \geq \sigma_2(\X) \geq \cdots \geq \sigma_{p_1 \wedge p_2} (\X)$, denote $\|\X\|_\F = \sqrt{\sum_{i} \sigma^2_i(\X)}$ and $\|\X\| = \sigma_1(\X)$ as its Frobenius norm and spectral norm, respectively. Also, we use $\X^{-1}$, $\X^{-\top}$ and $\X^\dagger$ to denote the inverse, transpose inverse, and Moore-Penrose inverse of $\X$, respectively. For any real symmetric matrix $\X \in \bbS^{p \times p}$ having eigendecomposition $\U \bSigma \U^\top$ with non-increasing eigenvalues on the diagonal of $\bSigma$, let $\lambda_i(\X)$ be the $i$th largest eigenvalue of $\X$, $\lambda_{\min}(\X)$ be the least eigenvalue of $\X$ and $\X^{1/2} = \U \bSigma^{1/2} \U^\top$. We say a symmetric matrix $\X$ is positive semidefinite (PSD) and denote $\X \succcurlyeq 0$ if and only if (iff) for any vector $\y \in \bbR^{p}$, $\y^\top \X \y \geq 0$. For two symmetric matrices $\X, \Y$, we say $\X \succcurlyeq \Y$ iff $\X - \Y \succcurlyeq 0$. Throughout the paper, the SVD or eigendecomposition of a rank $r$ matrix $\X$ refers to its economic version. We use bracket subscripts to denote sub-matrices. For example, $\X_{[i_1,i_2]}$ is the entry of $\X$ on the $i_1$-th row and $i_2$-th column; $\X_{[(r+1):p_1, :]}$ contains the $(r+1)$-th to the $p_1$-th rows of $\X$. In addition, $\I_r$ is the $r$-by-$r$ identity matrix and $\mathcal{I}$ denotes an identity operator. Let $\mathbb{O}_{p, r} = \{\U \in \bbR^{p \times r}: \U^\top \U=\I_r\}$ be the set of all $p$-by-$r$ matrices with orthonormal columns and $\mathbb{O}_r := \mathbb{O}_{r,r}$. For any $\U\in \mathbb{O}_{p, r}$, $P_{\U} = \U\U^\top$ represents the orthogonal projector onto the column space of $\U$; we also note $\U_\perp\in \mathbb{O}_{p, p-r}$ as the orthonormal complement of $\U$. For any linear operator ${\cal L}$, we denote $\cL^*$ as its adjoint operator. Finally, for a linear space $\cV$, we denote its dimension as $\dim(\cV)$. For two linear spaces $\cV_1, \cV_2$, the sum of $\cV_1$ and $\cV_2$ is denoted by $\cV_1+\cV_2 := \{\v_1 + \v_2| \v_1 \in \cV_1,\v_2 \in \cV_2\}$. If every vector in $\cV_1 + \cV_2$ can be uniquely decomposed into $\v_1 + \v_2$, where $\v_1 \in \cV_1,\v_2 \in \cV_2$, then we call the sum of $\cV_1$ and $\cV_2$ as the direct sum, denoted by $\cV_1 \oplus \cV_2 $, and have $\dim(\cV_1 \oplus \cV_2) = \dim(\cV_1) + \dim(\cV_2)$. For two Euclidean spaces $\cV_1$ and $\cV_2$, we say $\cV_1$ is orthogonal to $\cV_2$ and denote it by $\cV_1 \perp \cV_2 $ iff $\langle \v_1, \v_2 \rangle =0 $ for any $\v_1 \in \cV_1,\v_2 \in \cV_2$.

Given differentiable scalar and matrix-valued functions $f: \bbR^{p_1 \times p_2} \to \bbR$ and $\phi: \bbR^{p_1 \times p_2} \to \bbR^{q_1 \times q_2}$. The Euclidean gradient of $f$ at $\X$ is denoted as $\nabla f(\X)$ and $(\nabla f(\X))_{[i,j]} = \frac{\partial f(\X) }{\partial \X_{[i,j]}} $ for $i \in [p_1], j \in [p_2]$. The Euclidean gradient of $\phi$ is a linear operator from $\bbR^{p_1 \times p_2}$ to $\bbR^{q_1 \times q_2}$ defined as $(\nabla \phi (\X) [\Z] )_{[i,j]} = \sum_{k \in [p_1],l \in [p_2]} \frac{\partial ( \phi (\X) )_{[i,j]} }{\partial \X_{[k,l]}} \Z_{[k,l]}$ for any $\Z \in \bbR^{p_1 \times p_2}, i \in [q_1], j \in [q_2]$. Using this notation, given a twice continuously differentiable scalar function $f: \bbR^{p_1 \times p_2} \to \bbR$, we denote its Euclidean Hessian by $\nabla^2 f(\X)[\cdot] $, which is the gradient of $\nabla f(\X)$ and can be viewed as a linear operator from $\bbR^{p_1 \times p_2}$ to $\bbR^{p_1 \times p_2}$ satisfying
\begin{equation*}
	( \nabla^2 f(\X)[\Z] )_{[i,j]} =  \sum_{k \in [p_1],l \in [p_2]} \frac{\partial (\nabla f (\X) )_{[i,j]} }{\partial \X_{[k,l]}} \Z_{[k,l]} = \sum_{k \in [p_1],l \in [p_2]} \frac{\partial^2 f (\X) }{\partial \X_{[k,l]} \partial \X_{[i,j]} } \Z_{[k,l]}.
\end{equation*} We also define the bilinear form for the Hessian of $f$ as $\nabla^2 f(\X)[\Z_1, \Z_2]:= \langle \nabla^2 f(\X)[\Z_1], \Z_2 \rangle $ for any $\Z_1, \Z_2 \in \bbR^{p_1 \times p_2}$. Apart from the matrix representation of $\nabla^2 f(\X)$ above, we can also view $\nabla^2 f(\X)$ as a $(p_1 p_2)$-by-$(p_1 p_2)$ symmetric matrix and define its spectrum in the classic way. We say $\X$ is a Euclidean first-order stationary point (FOSP) of $f$ iff $\nabla f(\X) = \0$ and a Euclidean second-order stationary point (SOSP) of $f$ iff $\nabla f(\X) = \0$ and $\nabla^2 f(\X) \succcurlyeq 0$. Finally, we say a pair of matrices $(\L,\R)\in \bbR^{p_1 \times r} \times \bbR^{p_2 \times r}$ is of rank $r$ if $\L\R^\top$ has rank $r$. 

\subsection{Riemannian Optimization and Riemannian Geometry of Low-rank PSD and General Matrices} \label{sec: Riemannian-opt-background}
In this section, we first give a brief introduction to Riemannian optimization and then present the necessary preliminaries to perform Riemannian optimization on \eqref{eq: PSD-manifold-formulation} and \eqref{eq: general prob}. Finally, we provide the Euclidean/Riemannian gradient and Hessian expressions for the optimization problems \eqref{eq: PSD-manifold-formulation}-\eqref{eq: general factor with reg} considered in this paper.

Riemannian optimization concerns optimizing a real-valued function $f$ defined on a Riemannian manifold $\cM$, for which the readers are referred to \cite{absil2009optimization,boumal2020introduction,hu2020brief} for more details. Algorithms for continuous optimization over the Riemannian manifold often require calculations of Riemannian gradients and Riemannian Hessians. Suppose $\X \in \cM$ and the Riemannian metric and tangent space of $\cM$ at $\X$ are $\langle \cdot, \cdot \rangle_\X$ and $T_\X \cM$, respectively. Then the Riemannian gradient of a smooth function $f:\cM \to \bbR$ at $\X$ is defined as the unique tangent vector ${\rm grad}\, f(\X) \in T_\X \cM$ such that $\langle \grad \, f(\X),\Z \rangle_\X = {\rm D} \, f(\X)[\Z], \forall\, \Z \in T_\X \cM$, where ${\rm D} f(\X)[\Z]$ denotes the directional derivative of $f$ at point $\X$ along the direction $\Z$. The Riemannian Hessian of $f$ at $\X\in\cM$ is the linear map ${\rm Hess} \,f(\X)$ of $T_\X\cM$ onto itself defined as
\begin{equation} \label{def: Riemannain-Hessian}
	{\rm Hess}\, f(\X)[\Z] = \widebar{\nabla}_{\Z} {\rm grad}\, f, \,\quad \forall \Z \in T_\X \cM,
\end{equation}where $\widebar{\nabla}$ is the Riemannian connection on $\cM$ and can be viewed as the generalization of the directional derivative on the manifold \cite[Section 5.3]{absil2009optimization}. The bilinear form of Riemannian Hessian is defined as $\Hess f(\X)[\Z_1, \Z_2]:= \langle \Hess f(\X)[\Z_1], \Z_2 \rangle_\X$ for any $\Z_1, \Z_2 \in T_{\X}\cM$. We say $\X \in \cM$ is a Riemannian FOSP of $f$ iff $\grad f(\X) = \0$ and a Riemannian SOSP of $f$ iff $\grad f(\X) = \0$ and $\Hess f(\X) \succcurlyeq 0$. Moreover, we call a Riemannian or Euclidean FOSP a strict saddle if and only if the Riemannian or Euclidean Hessian evaluated at this point has a strict negative eigenvalue. 

For optimization problems in \eqref{eq: PSD-manifold-formulation} and \eqref{eq: general prob}, two manifolds of particular interests are the set of rank $r$ PSD matrices $\cM_{r+}:=\left\{ 
\X\in \bbS^{p\times p} \mid {\rm rank}(\X) = r, \X \succcurlyeq 0
\right\}$ and the set of rank $r$ matrices $\cM_r:=\left\{ 
\X\in \bbR^{p_1\times p_2}\mid {\rm rank}(\X) = r
\right\}$. \cite{lee2013smooth,helmke2012optimization,vandereycken2010riemannian} showed that $\cM_{r+}$ and $\cM_r$ are smooth embedded submanifolds of $\bbR^{p \times p}$ and $\bbR^{p_1 \times p_2}$, respectively and established their Riemannian geometry as follows. 
\begin{Lemma} {\rm(\cite[Chapter 5]{helmke2012optimization}, \cite[Proposition 5.2]{vandereycken2010riemannian}, \cite[Example 8.14]{lee2013smooth})}\label{lm: Mr+ and Mr manifold}
	$\cM_{r+},\cM_r$ are smooth embedded submanifolds of $\bbR^{p \times p}$ and $\bbR^{p_1 \times p_2}$ with dimensions $(pr - r(r-1)/2)$ and $(p_1 + p_2 -r)r$, respectively. The tangent space $T_{\X}\cM_{r+}$ at $\X \in \cM_{r+}$ with the eigendecomposition $\X = \U\bSigma \U^\top$ is given by
	\begin{equation}
	\label{eq:tangent Mr+}
	T_\X \cM_{r+} = \left\{
	[\U\quad \U_{\perp}] \begin{bmatrix}
	\S & \D^\top \\[2pt]
	\D & \0
	\end{bmatrix}
	[\U\quad \U_{\perp}]^\top: \S \in \bbS^{r \times r}, \D \in \bbR^{(p-r) \times r}
	\right\}.
	\end{equation}
	 The tangent space $T_{\X}\cM_{r}$ at $\X \in \cM_{r}$ with SVD $\X = \U\bSigma \V^\top$ is given by
	\begin{equation}
	\label{eq:tangent Mr}
	T_\X \cM_r = \left\{
	[\U\quad \U_{\perp}] \begin{bmatrix}
	\S & \D_2^\top \\[2pt]
	\D_1 & \0
	\end{bmatrix}
	[\V\quad \V_{\perp}]^\top: \S \in \bbR^{r \times r}, \D_1 \in \bbR^{(p_1 - r)\times r},  \D_2 \in \bbR^{(p_2 - r)\times r}
	\right\}.
	\end{equation}
\end{Lemma}

Throughout the paper, we equip the Riemannian manifolds $\mathcal{M}_{r+}$ and $\mathcal{M}_r$ with the metric induced by the Euclidean inner product, i.e., $\langle \U, \V \rangle = \mathrm{trace}(\U^\top \V)$. Given $\X \in \cM_{r+}$ with eigendecomposition $\U\bSigma \U^\top$ or $\X \in \cM_{r}$ with SVD $\U\bSigma \V^\top$, the orthogonal projectors $P_{T_\X \cM_{r+}}(\cdot)$ and $P_{T_\X \cM_r}(\cdot)$, which project any matrix onto $T_\X \cM_{r+}$ and $T_\X \cM_r$, are given as follows
\begin{equation} \label{eq: tangent space projection}
	\begin{split}
		 P_{T_\X \cM_{r+}}(\Z) &= P_\U \Z P_\U + P_{\U_\perp} \Z P_\U + P_\U \Z P_{\U_\perp},\quad \forall\, \Z\in \bbS^{p \times p },\\
		P_{T_\X \cM_r}(\Z) &= P_\U \Z P_\V + P_{\U_\perp} \Z P_\V + P_\U \Z P_{\V_\perp},\quad \forall\, \Z\in\bbR^{p_1\times p_2}.
 	\end{split}
\end{equation}

Next, we collectively give the expressions for gradients and Hessians under both manifold formulations \eqref{eq: PSD-manifold-formulation} and \eqref{eq: general prob} and factorization formulations \eqref{eq: PSD factorization}, \eqref{eq: general factor formu} and \eqref{eq: general factor with reg} in Propositions \ref{prop: gradient-exp} and \ref{prop: Hessian-exp}, respectively. The proof is postponed to Appendix.  

\begin{Proposition}[Riemannian and Euclidean Gradients] \label{prop: gradient-exp}
	The Riemannian and Euclidean gradients under the manifold and the factorization formulations are:
	\begin{itemize}[leftmargin=*]
		\item PSD case: \begin{equation*}
	\begin{split}
		\grad f(\X) &= P_{\U} \nabla f(\X)P_{\U} + P_{\U_\perp} \nabla f(\X)P_{\U} + P_{\U} \nabla f(\X)P_{\U_\perp},\\
		\nabla g(\Y) &= 2 \nabla f(\Y \Y^\top) \Y.
	\end{split}
	\end{equation*} Here $\U\in \bbO_{p,r}$ is formed by the top $r$ eigenvectors of $\X$.
	\item General case: \begin{equation*}
	\begin{split}
		\grad f(\X) &= P_{\U} \nabla f(\X)P_{\V} + P_{\U_\perp} \nabla f(\X)P_{\V} + P_{\U} \nabla f(\X)P_{\V_\perp},\\
		\nabla g(\L,\R) &= \begin{bmatrix}
				\nabla_{\L} g(\L, \R)\\
				\nabla_{\R} g(\L, \R)
			\end{bmatrix} = \begin{bmatrix}
				\nabla f(\L\R^\top )\R\\
				(\nabla f(\L\R^\top ))^\top \L
			\end{bmatrix}, \\
		\nabla g_\reg(\L,\R) &= \begin{bmatrix}
			\nabla_{\L} g_\reg(\L, \R)\\
			\nabla_{\R} g_\reg(\L, \R)
		\end{bmatrix} = \begin{bmatrix}
			\nabla_\L g(\L, \R) + 2\mu \L(\L^\top \L- \R^\top \R) \\
			\nabla_\R g(\L, \R) - 2\mu \R(\L^\top \L- \R^\top \R)
		\end{bmatrix}.
	\end{split}
	\end{equation*} Here $\U$ and $\V$ are the top $r$ left and right singular vectors of $\X$.
	\end{itemize}
\end{Proposition}

\begin{Proposition}[Riemannian and Euclidean Hessians] \label{prop: Hessian-exp}
The Riemannian and Euclidean Hessians under the manifold and the factorization formulations are:
\begin{itemize}[leftmargin=*]
		\item PSD case: Suppose $\X \in \cM_{r+}$ has eigendecomposition $\U \bSigma \U^\top$, $\xi = [\U \quad \U_\perp] \begin{bmatrix}
			\S & \D^\top\\
			\D & \0
		\end{bmatrix} [\U \quad \U_\perp]^\top \in T_{\X}\cM_{r+}$ and $\A \in \bbR^{p \times r}$. Then
	\begin{equation} \label{eq: Hessian-PSD}
	\begin{split}
		\Hess f(\X)[\xi, \xi] &= \nabla^2 f(\X)[\xi, \xi] + 2\langle \nabla f(\X), \U_\perp \D \bSigma^{-1} \D^\top \U_\perp^\top \rangle,\\
		\nabla^2 g(\Y)[\A, \A] &= \nabla^2 f(\Y \Y^\top)[\Y \A^\top + \A \Y^\top,\Y \A^\top + \A \Y^\top] + 2 \langle \nabla f(\Y \Y^\top), \A\A^\top \rangle.
	\end{split}
	\end{equation} 
	\item General case: Suppose $\X \in \cM_{r}$ has SVD $\U \bSigma \V^\top$, $\xi = [\U \quad \U_\perp] \begin{bmatrix}
			\S & \D_2^\top\\
			\D_1 & \0
		\end{bmatrix} [\V \quad \V_\perp]^\top \in T_{\X}\cM_{r}$ and $\A = [\A_L^\top \quad \A_R^\top]^\top$ with $\A_L \in \bbR^{p_1 \times r},\A_R \in \bbR^{p_2 \times r}$. Then
 \begin{equation}\label{eq: Hessian-general}
	\begin{split}
		\Hess f(\X)[\xi, \xi] &= \nabla^2 f(\X)[\xi, \xi] + 2\langle \nabla f(\X), \U_\perp \D_1 \bSigma^{-1} \D_2^\top \V_\perp^\top \rangle,\\
		\nabla^2 g(\L,\R)[\A,\A] &= \nabla^2 f(\L\R^\top)[\L\A_R^\top + \A_L \R^\top,\L\A_R^\top + \A_L \R^\top] + 2 \langle \nabla f(\L\R^\top), \A_L\A_R^\top \rangle, \\
		\nabla^2 g_\reg(\L,\R)[\A,\A] &= \nabla^2 g(\L,\R)[\A,\A] + \mu \| \L^\top \A_L + \A_L^\top \L- \R^\top \A_R - \A_R^\top \R \|_\F^2\\
		& \quad \quad   + 2\mu \langle \L^\top \L- \R^\top \R, \A_L^\top \A_L - \A_R^\top \A_R \rangle .
	\end{split}
	\end{equation}
	\end{itemize}
\end{Proposition}

In Proposition \ref{prop: Hessian-exp}, we give the quadratic expressions of the Hessians as we use them exclusively throughout the paper. It is relatively easy to obtain the general bilinear expressions by noting that $\nabla^2 g(\Y)[\A,\B] =(\nabla^2 g(\Y)[\A+\B,\A+\B] - \nabla^2 g(\Y)[\A-\B,\A-\B]  )/4 $ and similarly for the Riemannian Hessian. 

\section{Geometric Connection of Manifold and Factorization Formulations: PSD Case} \label{sec: connection-PSD}

In this section, we present the geometric landscape connections between the manifold formulation \eqref{eq: PSD-manifold-formulation} and the factorization formulation \eqref{eq: PSD factorization} in low-rank PSD matrix optimization. Here we single out the PSD case from the general case and give it a detailed discussion for a few reasons: first, the low-rank PSD matrix manifold is different from the low-rank matrix manifold; second, low-rank PSD optimization appears in many real applications such as phase retrieval \citep{fienup1982phase}, rank-$1$ covariance sensing \citep{chen2015exact,cai2015rop} and is interesting on its own; third, the results in PSD case are easier to follow and also convey most of the essential messages in the general setting. 

{We begin with a few more definitions. Suppose $\Y \in \bbR^{p \times r}$ is of rank $r$, $\X = \Y \Y^\top$ has eigendecomposition $\U \bSigma \U^\top$, and $\P = \U^\top \Y$. For any $\A \in \bbR^{p \times r}$, define
\begin{equation} \label{def: A-xi-PSD}
\xi^{\A}_{\Y}:= \Y \A^\top + \A \Y^\top  = [\U \quad \U_\perp] \begin{bmatrix}
		\P \A^\top \U + \U^\top \A \P^\top & \, \P \A^\top \U_\perp \\
		\U_\perp^\top \A \P^\top & \0
	\end{bmatrix}[\U \quad \U_\perp]^\top  \in T_{\X}\cM_{r+}.
\end{equation}
For any $\xi = [\U \quad \U_\perp] \begin{bmatrix}
			\S & \D^\top\\
			\D & \0
		\end{bmatrix} [\U \quad \U_\perp]^\top \in T_{\X}\cM_{r+}$ with $\S\in  \bbS^{r \times r}$, define
		\begin{equation}\label{def: xi-setA-PSD}
			\scA^\xi_{\Y} := \left\{ \A : \Y \A^\top + \A \Y^\top = \xi\right\}.
		\end{equation}
In the following Lemma \ref{lm: tangent-vector-equiva-PSD}, we quickly check $\scA^\xi_{\Y}$ is nonempty and forms an $(r^2 - r)/2 $ dimensional subspace.
\begin{Lemma}\label{lm: tangent-vector-equiva-PSD}
	Suppose $\Y \in \bbR^{p \times r}$ is of rank $r$, $\X = \Y \Y^\top$ has eigendecomposition $\U \bSigma \U^\top$, and $\P = \U^\top \Y$. 
	Given any $\xi \in T_{\X}\cM_{r+}$, it holds $\scA^\xi_{\Y} := \left\{ \A: \A = (\U\S_1 + \U_\perp \D )\P^{-\top}  \in \bbR^{p \times r}, \S_1 + \S_1^\top = \S \right \}$. 
\end{Lemma} 
The motivation behind the constructions of $\xi^{\A}_{\Y}$ and $\scA_{\Y}^{\xi}$ is to find a correspondence between $\bbR^{p \times r}$ and $T_\X \cM_{r+}$. Later, this correspondence will be used to establish the connections of Riemannian and Euclidean Hessians \eqref{eq: R-E-Hessian-PSD} in Theorem \ref{th: RHessian-EHessian PSD}.}

From Lemma \ref{lm: tangent-vector-equiva-PSD}, we can see there is no one-to-one correspondence between $\bbR^{p \times r}$ and $T_\X \cM_{r+}$ due to their mismatched dimensions, and in particular given $\Y$, $\A$ and $\xi$, $\xi^{\A}_{\Y}$ is a single matrix while there exists an $(r^2 - r)/2 $ dimensional subspace $\scA_{\Y}^{\xi}$ such that $\Y \A^{'\top} + \A' \Y^\top = \xi$ for any $\A' \in \scA_{ \Y}^{\xi}$. To handle the mismatch, we introduce the following decomposition of $\bbR^{p \times r}$:
\begin{Lemma}\label{lm: decomposition-Rp-r-PSD}
	Suppose the conditions in Lemma \ref{lm: tangent-vector-equiva-PSD} hold. Then $\bbR^{p \times r} = \scA_{\rnull}^{\Y} \oplus \scA_{\overline{\rnull}}^{\Y}$ for
	\begin{equation*}
	\begin{split}
	\scA_{\rnull}^{\Y} &= \left\{ \A: \A = \U \S\P^{-\top} , \S +\S^\top = \0 \in \bbR^{r \times r}    \right \},\\
	\scA_{\overline{\rnull}}^{\Y} &= \left\{ \A: \A = (\U\S + \U_\perp \D )\P^{-\top}  , \D \in \bbR^{(p-r) \times r}, \S \bSigma^{-1} \in \bbS^{r \times r}  \right \}.
	\end{split}
	\end{equation*}
	Moreover, $\dim(\scA_{\rnull}^{\Y}) = (r^2-r)/2$, $\dim(\scA_{\overline{\rnull}}^{\Y}) = pr- (r^2-r)/2$, and $\scA_{\rnull}^{\Y}$ is orthogonal to $\scA_{\overline{\rnull}}^{\Y}$, i.e., $\scA_{\rnull}^{\Y} \perp \scA_{\overline{\rnull}}^{\Y}$.
\end{Lemma}

\begin{figure}
	\centering
	\includegraphics[width=0.6\textwidth]{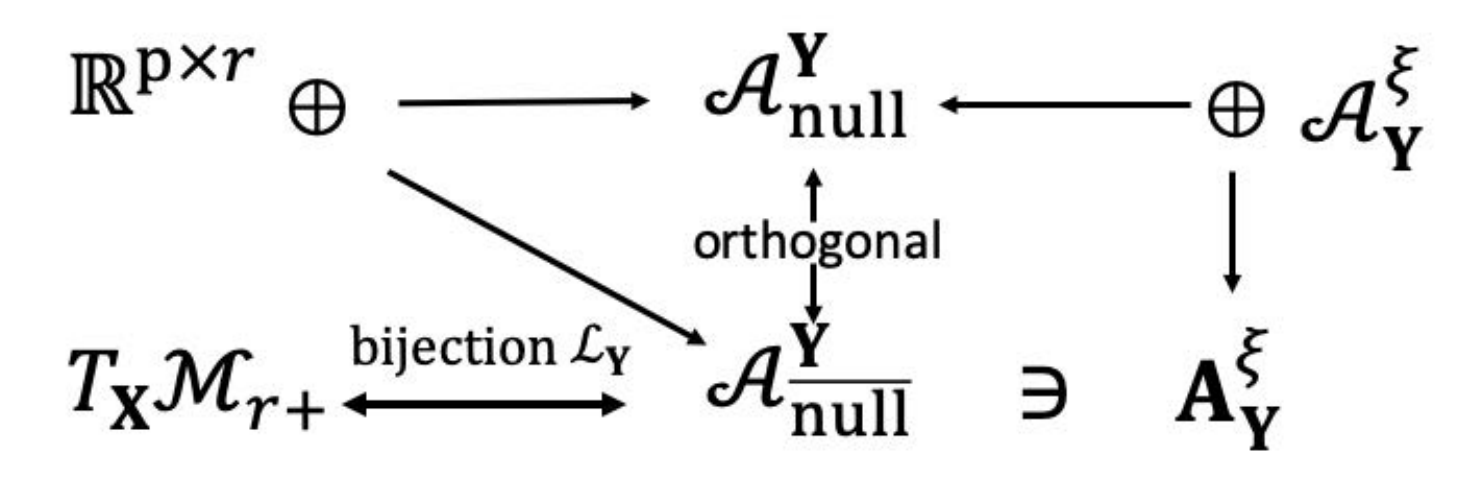}
	\caption{Relationship of $\bbR^{p \times r}$, $T_\X \cM_{r+}$, $\scA_{\rnull}^{\Y}$, $\scA_{\overline{\rnull}}^{\Y}$, $\A_\Y^\xi$, and $\scA_\Y^\xi$.}
	\label{fig: PSD_decom_ill}
\end{figure}

By decomposing $\bbR^{p \times r}$ into $\scA_{\rnull}^{\Y}$ and $\scA_{\overline{\rnull}}^{\Y}$, it is easy to check $\xi_{\Y}^{\A} = \0$ if and only if $\A \in \scA_{\rnull}^{\Y}$. In the following Proposition \ref{prop: bijection-PSD}, we show $\scA_\Y^\xi$ can be decomposed as the direct sum of $\scA_{\rnull}^{\Y}$ and a singleton from $\scA_{\overline{\rnull}}^{\Y}$. In addition, there is a bijective linear map, $\cL_\Y$, between $\scA_{\overline{\rnull}}^{\Y}$ and $T_\X \cM_{r+}$ and $\scA_{\rnull}^{\Y}$ is the null space of $\cL_\Y$. A pictorial illustration of the relationship of these subspaces is given in Figure \ref{fig: PSD_decom_ill}.

\begin{Proposition}[Decomposition of $\scA_\Y^\xi$ and Bijection Between $\scA_{\overline{\rnull}}^{\Y}$ and $T_\X \cM_{r+}$] \label{prop: bijection-PSD} Suppose the conditions in Lemma \ref{lm: tangent-vector-equiva-PSD} hold and $\xi = [\U \quad \U_\perp] \begin{bmatrix}
			\S & \D^\top\\
			\D & \0
		\end{bmatrix} [\U \quad \U_\perp]^\top \in T_{\X}\cM_{r+}$. Then 
	$\scA_\Y^\xi = \A_\Y^\xi \oplus \scA_{\rnull}^{\Y}$, where $\A_\Y^\xi =  (\U\widebar{\S} + \U_\perp \D )\P^{-\top} \in \scA_{\overline{\rnull}}^{\Y}$ and $\widebar{\S}$ is the unique solution of the linear equation system $\widebar{\S} \bSigma^{-1} = \bSigma^{-1} \widebar{\S}^\top $ and $\widebar{\S} + \widebar{\S}^\top = \S$. 
	
	Moreover, there is a bijective linear map $\mathcal{L}_\Y$ between $\scA_{\overline{\rnull}}^{\Y}$ and $T_\X \cM_{r+}$ given as follows
	\begin{equation*}
		\cL_\Y : \A \in \scA_{\overline{\rnull}}^{\Y} \longrightarrow \xi_\Y^\A \in T_\X\cM_{r+} \quad \text{ and } \quad \cL^{-1}_\Y : \xi \in T_\X\cM_{r+} \longrightarrow \A_\Y^\xi \in \scA_{\overline{\rnull}}^{\Y}.
	\end{equation*} Finally, recall $\X = \Y \Y^\top$, we have the following spectrum bounds for $\cL_\Y$:
	\begin{equation} \label{ineq: bijection-spectrum}
		2 \sigma_r(\X) \|\A\|_\F^2 \leq \|\cL_\Y(\A)\|_\F^2 \leq 4 \sigma_1(\X) \|\A\|_\F^2, \quad \forall \A \in \scA_{\overline{\rnull}}^{\Y}.
	\end{equation}
\end{Proposition}
{\noindent \bf Proof of Proposition \ref{prop: bijection-PSD}.} We divide the proof into two steps: in Step 1, we prove the decomposition for $\scA_\Y^\xi $; in Step 2, we show $\cL_\Y$ is a bijection and prove the spectrum bounds.

{\bf Step 1}. Noting that $\bSigma=\{\{\bSigma_{[i,i]}\}_{i=1}^r\}$ is a diagonal matrix, the linear equation system $\widebar{\S} \bSigma^{-1} = \bSigma^{-1} \widebar{\S}^\top$, $\widebar{\S} + \widebar{\S}^\top = \S$ is equivalent to 
$$\widebar{\S}_{[i,j]}\bSigma_{[j,j]}^{-1} = \bSigma_{[i,i]}^{-1}\widebar{\S}_{[j,i]}, \quad \widebar{\S}_{[i,j]} + \widebar{\S}_{[j,i]} = \S_{[i,j]},\quad 1\leq i,j\leq r.$$
Here, we use the fact that $\S$ is symmetric, i.e., $\S_{[i,j]} = \S_{[j,i]}$. By calculations, we know the equation system above is further equivalent to
\begin{equation}\label{eq:equ_solution}
\bar{\S}_{[i,j]} = \S_{[i,j]}\bSigma_{[j,j]}/(\bSigma_{[i,i]}+\bSigma_{[j,j]}),\quad \bar{\S}_{[j,i]} = \S_{[i,j]}\bSigma_{[i,i]}/(\bSigma_{[i,i]}+\bSigma_{[j,j]}), \quad 1\leq i \leq j \leq r.
\end{equation}
Namely, $\widebar{\S} \bSigma^{-1} = \bSigma^{-1} \widebar{\S}^\top $, $\widebar{\S} + \widebar{\S}^\top = \S$ has the unique solution as presented in \eqref{eq:equ_solution}.
Therefore, $\A_\Y^\xi$ is well-defined for any $\xi \in T_{\X}\cM_{r_+}$.

	At the same time, given $\A = (\U\S_1 + \U_\perp \D )\P^{-\top} \in \scA_\Y^\xi$, we can check $\A - \A_\Y^\xi \in \scA_{\rnull}^{\Y}$. In addition, $\A_\Y^\xi + \A \in \scA_\Y^\xi$ for any $\A \in \scA_{\rnull}^{\Y}$. This shows $\scA_\Y^\xi = \A_\Y^\xi \oplus \scA_{\rnull}^{\Y}$.
	
{\bf Step 2}. Notice both $ \scA_{\overline{\rnull}}^{\Y}$ and $T_\X\cM_{r+}$ are of dimension $(pr - (r^2-r)/2)$. Suppose $\cL_\Y' : \xi \in T_\X\cM_{r+} \longrightarrow \A_\Y^\xi \in \scA_{\overline{\rnull}}^{\Y}$. For any $\xi = [\U \quad \U_\perp] \begin{bmatrix}
			\S & \D^\top\\
			\D & \0
		\end{bmatrix} [\U \quad \U_\perp]^\top \in T_\X\cM_{r+}$, we have
\begin{equation} \label{eq: bijection-PSD}
	\cL_\Y ( \cL_\Y'(\xi) ) = \cL_\Y (\A_\Y^\xi) =  [\U \quad \U_\perp] \begin{bmatrix}
		\P \A_\Y^{\xi\top} \U + \U^\top \A_\Y^\xi \P^\top & \P \A_\Y^{\xi\top} \U_\perp \\
		\U_\perp^\top \A_\Y^\xi \P^\top & \0
	\end{bmatrix}[\U \quad \U_\perp]^\top = \xi.
\end{equation}
Since $\cL_\Y$ and $\cL_\Y'$ are linear maps, \eqref{eq: bijection-PSD} implies $\cL_\Y$ is a bijection and $\cL_\Y' = \cL_\Y^{-1}$. 

Next, we provide the spectrum bounds for $\cL_\Y$. Suppose $\A = (\U\S + \U_\perp \D )\P^{-\top} \in \scA_{\overline{\rnull}}^{\Y}$. Then,
	\begin{equation} \label{ineq: upper-bound-A-norm}
	\begin{split}
		\|\A\|_\F^2 \leq ( \|\S\|_\F^2 + \|\D\|_\F^2 ) \sigma^2_1(\P^{-\top}) \overset{\P \P^\top = \bSigma }= ( \|\S\|_\F^2 + \|\D\|_\F^2 )/\sigma_r(\X);\\
		\langle \S^\top, \S \rangle \overset{(a)}= \langle \bSigma \S \bSigma^{-1}, \S \rangle = \langle \bSigma^{1/2} \S \bSigma^{-1/2}, \bSigma^{1/2} \S \bSigma^{-1/2} \rangle \geq 0.
	\end{split}
	\end{equation} Here (a) is because $ \S \bSigma^{-1} = \bSigma^{-1} \S^\top $ by the definition of $\scA_{\overline{\rnull}}^{\Y}$. So
\begin{equation*}
\begin{split}
	\|\cL_\Y(\A)\|_\F^2 = \|\xi^{\A}_{\Y}\|_\F^2 & \overset{ \eqref{def: A-xi-PSD}  }= \|\P \A^\top \U + \U^\top \A \P^\top \|_\F^2 + 2 \|\U_\perp^\top \A \P^\top \|_\F^2\\
	& = \|\S^\top + \S\|_\F^2 + 2\|\D\|_\F^2\\
	& \overset{\eqref{ineq: upper-bound-A-norm} } \geq 2 (\|\S\|_\F^2 + \|\D\|_\F^2) \overset{\eqref{ineq: upper-bound-A-norm} } \geq 2\sigma_r(\X) \|\A\|_\F^2,
\end{split}
\end{equation*} and
\begin{equation*}
		\begin{split}
		\|\cL_\Y(\A)\|_\F^2  = \|\xi^{\A}_{\Y}\|_\F^2 & \overset{ \eqref{def: A-xi-PSD}  }= \|\P \A^\top \U + \U^\top \A \P^\top \|_\F^2 + 2 \|\U_\perp^\top \A \P^\top \|_\F^2 \\
		& \leq 4 \|\U^\top \A \P^\top\|_\F^2 + 2 \|\U_\perp^\top \A \P^\top \|_\F^2\\
		& = 2 \|\U^\top \A \P^\top\|_\F^2 + 2 \|\A \P^\top\|_\F^2 \\
		& \leq 4 \sigma^2_1(\P) \|\A\|_\F^2 = 4 \sigma_1(\X) \|\A\|_\F^2.
	\end{split}
\end{equation*} This finishes the proof of this proposition. \quad $\blacksquare$

Next, we present our first main result on the geometric landscape connection between formulations \eqref{eq: PSD-manifold-formulation} and \eqref{eq: PSD factorization}. 

\begin{Theorem}({\bf Geometric Landscape Connection Between Manifold and Factorization Formulations (PSD Case)}) \label{th: RHessian-EHessian PSD}
	Suppose $\Y \in \bbR^{p \times r}$ is of rank $r$ and $\X = \Y \Y^\top$. Then
	\begin{equation} \label{eq: gradient-connect-PSD}
	\begin{split}
		\grad f(\X) = \left( \nabla g(\Y) \Y^\dagger + (\nabla g(\Y) \Y^\dagger)^\top(\I_p - \Y \Y^\dagger) \right)/2\quad \text{ and } \quad
		\nabla g(\Y) = 2 \grad f(\X) \Y.
	\end{split}
	\end{equation}
	
Furthermore, if $\Y$ is a Euclidean FOSP of \eqref{eq: PSD factorization}, we have: 
\begin{equation} \label{eq: R-E-Hessian-PSD}
	\begin{split}
		\nabla^2 g(\Y)[\A,\A] &= \Hess f(\X)[\xi_\Y^\A,\xi_\Y^\A ], \quad \forall \A \in \bbR^{p \times r} ;\\
		\Hess f(\X)[\xi,\xi] &= \nabla^2 g(\Y)[\cL_\Y^{-1}(\xi) ,\cL_\Y^{-1}(\xi)], \quad \forall \xi \in T_\X\cM_{r+}.
	\end{split}
	\end{equation}
	More precisely,
	\begin{equation} \label{eq: R-E-Hessian-PSD-2}
		\begin{split}
			\nabla^2 g(\Y)[\A] &= \0, \quad \forall \A \in \scA_{\rnull}^\Y;\\
		\nabla^2 g(\Y)[\A,\A] &= \Hess f(\X)[\cL_\Y (\A),\cL_\Y (\A) ],\quad \forall \A \in \scA_{\overline{\rnull}}^{\Y}.
		\end{split}
	\end{equation}
	
	Finally, $\Hess f(\X)$ has $(pr - (r^2 -r)/2) $ eigenvalues and $ \nabla^2 g(\Y)$ has $pr$ eigenvalues. $\nabla^2 g(\Y)$ has at least $(r^2 -r)/2$ zero eigenvalues which correspond to the zero eigenspace $\scA_{\rnull}^\Y$. Denote the rest of the $(pr - (r^2 -r)/2) $ possibly non-zero eigenvalues of $\nabla^2 g(\Y)$ from the largest to the smallest as $\widebar{\lambda}_1(\nabla^2 g(\Y)),\ldots, \widebar{\lambda}_{pr - (r^2 -r)/2}(\nabla^2 g(\Y))$. Then for $i=1,\ldots,pr - (r^2 -r)/2$:
	
	{\it $\widebar{\lambda}_i(\nabla^2 g(\Y))$ is sandwiched between $2\sigma_r(\X) \lambda_i(\Hess f(\X)) $ and $4\sigma_1(\X) \lambda_i(\Hess f(\X)) $.}
\end{Theorem}

{\noindent \bf Proof of Theorem \ref{th: RHessian-EHessian PSD}.} First, suppose $\Y \Y^\top$ has eigendecomposition $\U \bSigma \U^\top$. Then $\Y$ lies in the column space spanned by $\U$ and $\Y \Y^\dagger = P_\U$. So  \eqref{eq: gradient-connect-PSD} is by direct calculation from the gradient expressions in Proposition \ref{prop: gradient-exp}. The rest of the proof is divided into two steps: in Step 1, we prove \eqref{eq: R-E-Hessian-PSD} and \eqref{eq: R-E-Hessian-PSD-2}; in Step 2, we prove the individual eigenvalue connection between $\Hess f(\X)$ and $\nabla^2 g(\Y)$.

{\bf Step 1.} We begin by proving the first equality in \eqref{eq: R-E-Hessian-PSD}. Since $\Y$ is a Euclidean FOSP of \eqref{eq: PSD factorization}, by \eqref{eq: gradient-connect-PSD} we have $\X = \Y \Y^\top$ is a Riemannian FOSP of \eqref{eq: PSD-manifold-formulation} and $ \nabla f(\X) = P_{\U_\perp} \nabla f(\X) P_{\U_\perp}$. Let $\P = \U^\top \Y$. Given $\A \in \bbR^{p \times r}$, we have
\begin{equation} \label{eq: Hessian-con-gradient-1}
	\langle \nabla f(\X), P_{\U_\perp} \A \P^{\top} \bSigma^{-1} \P \A^\top  P_{\U_\perp}    \rangle = \langle \nabla f(\X), \A \A^\top \rangle,
\end{equation} where the equality is because $\P$ is nonsingular, $\P \P^\top = \bSigma$ and $ \nabla f(\X) = P_{\U_\perp} \nabla f(\X) P_{\U_\perp} $.

Then by Proposition \ref{prop: Hessian-exp},
\begin{equation}\label{eq: Hessian-con-Hessian-1}
\begin{split}
	\nabla^2 g(\Y)[\A,\A] &=  \nabla^2 f(\Y \Y^\top)[\Y \A^\top + \A \Y^\top,\Y \A^\top + \A \Y^\top] + 2 \langle \nabla f(\Y \Y^\top), \A\A^\top \rangle\\
	& \overset{ \text{Lemma } \ref{lm: tangent-vector-equiva-PSD}, \eqref{eq: Hessian-con-gradient-1} }= \nabla^2 f(\X)[\xi_\Y^\A,\xi_\Y^\A] + 2\langle \nabla f(\X), P_{\U_\perp} \A \P^{\top} \bSigma^{-1} \P \A^\top  P_{\U_\perp}    \rangle\\
	& = \Hess f(\X)[\xi_\Y^\A,\xi_\Y^\A],
\end{split}
\end{equation} 
where the last equality follows from the expressions of  $\Hess f(\X)$  and $\xi_\Y^\A$  in \eqref{eq: Hessian-PSD} and \eqref{def: A-xi-PSD}, respectively.
This finishes the proof of the first equality in \eqref{eq: R-E-Hessian-PSD}.

Moreover, for any $\A \in \scA_{\rnull}^\Y, \B \in \bbR^{p \times r}$:
\begin{equation} \label{eq: zero-eigenspace-PSD-argument}
\begin{split}
	\nabla^2 g(\Y)[\A,\B] &= \left( \nabla^2 g(\Y)[\A+ \B, \A + \B] - \nabla^2 g(\Y)[\A- \B, \A - \B]  \right) /4 \\
	& \overset{ \eqref{eq: Hessian-con-Hessian-1} }=  \left( \Hess f(\X)[\xi_\Y^{\A+\B}, \xi_\Y^{\A+\B}] - \Hess f(\X)[\xi_\Y^{\A-\B}, \xi_\Y^{\A-\B}]  \right) /4\\
	& = \left( \Hess f(\X)[\xi_\Y^{\A} + \xi_\Y^{\B}, \xi_\Y^{\A} + \xi_\Y^{\B}] - \Hess f(\X)[\xi_\Y^{\A} - \xi_\Y^{\B}, \xi_\Y^{\A} - \xi_\Y^{\B}]  \right) /4\\
	& = \Hess f(\X)[\xi_\Y^{\A},\xi_\Y^{\B}  ] \overset{(a)}= 0,
\end{split}
\end{equation} where (a) is because $\xi_\Y^\A = \0$ for $\A \in  \scA_{\rnull}^\Y$. This implies the first equality in \eqref{eq: R-E-Hessian-PSD-2}. The second equality in \eqref{eq: R-E-Hessian-PSD-2} follows directly from the first equality in \eqref{eq: R-E-Hessian-PSD} and the definition of $\cL_\Y$. Finally, since $\cL_\Y$ is a bijection, the second equality in \eqref{eq: R-E-Hessian-PSD} follows from the second equality in \eqref{eq: R-E-Hessian-PSD-2}. 

{\bf Step 2.} $\Hess f(\X)$ and $\nabla^2 g(\Y)$ are by definition linear maps from $T_{\X} \cM_{r+}$ and $\bbR^{p \times r}$ to $T_{\X} \cM_{r+}$ and $\bbR^{p \times r}$, respectively. Because $T_{\X} \cM_{r+}$ is of dimension $(pr - (r^2 -r)/2)$, the number of eigenvalues of $\Hess f(\X)$ and $\nabla^2 g(\Y)$ are $(pr - (r^2 -r)/2)$ and $pr$, respectively. By the first equality in \eqref{eq: R-E-Hessian-PSD-2}, we have $\scA_{\rnull}^\Y$ is the eigenspace of $(r^2 -r)/2$ zero eigenvalues of $\nabla^2 g(\Y)$ and the rest of the $(pr-(r^2-r)/2)$ possibly non-zero eigenvalues of $\nabla^2 g(\Y)$ span the eigenspace $\scA_{\overline{\rnull}}^{\Y}$. 
Restricting to $\scA_{\overline{\rnull}}^{\Y}$ and $T_\X \cM_{r+}$ and using \eqref{eq: R-E-Hessian-PSD-2}, \eqref{ineq: bijection-spectrum} and Lemma \ref{lm: two-symmetric-matrix-spectrum-connection} in the Appendix, we have $\widebar{\lambda}_i(\nabla^2 g(\Y))$ is sandwiched between $2\sigma_r(\X) \lambda_i(\Hess f(\X)) $ and $4\sigma_1(\X) \lambda_i(\Hess f(\X)) $. \quad $\blacksquare$.

\begin{Remark}({\bf Necessity of First-order Property in Connecting Riemannian and Euclidean Hessians}) \label{rem: necessity-of-FOSP-in-connection-Hessian}
	The following example shows that the assumption on the first-order stationary property is necessary for establishing the connection of the Riemannian and the Euclidean Hessians in Theorem \ref{th: RHessian-EHessian PSD}. Consider a special case that $p=r=1$, the objective functions are $f(x)$ ($x > 0$) and $g(y) = f(y^2)$, where both parameters $x, y$ are scalars. In such a scenario, when $x = y^2$ we have
	\begin{equation*}
	  g'(y) = 2y f'(y^2); \quad  g''(y) = \frac{\partial^2}{\partial y^2}f(y^2) = 2f'(y^2) + 4y^2 f''(y^2); \quad \Hess f(x) = f''(x) = f''(y^2),
	\end{equation*} where $f'$ and $f''$ denote the first and second derivatives of $f$.
	If $y$ is not a first-order stationary point of $g$ (namely $f'(y^2)$ can be non-zero without any constraint), $2f'(y^2) + 4y^2 f''(y^2)$ and $f''(y^2)$ do not necessarily share the same sign or hold any sandwich inequality.
\end{Remark}

\begin{Remark}[Connection of $\scA_{\rnull}^{\Y}$ and Rotational Invariance of $g(\Y)$] {We note $\scA_{\rnull}^{\Y}$ is also called the vertical space in studying the Riemannian quotient geometry of $\bbR^{p\times r}_*/\bbO_r$, where $\bbR^{p\times r}_*$ is the set of $p$-by-$r$ full column rank matrices \citep{massart2020quotient}. It has also appeared in \cite[Theorem 2, Example 4]{li2019symmetry} in analyzing the landscape of low-rank PSD matrix factorization.} By assuming $f$ is convex, \cite{li2019symmetry} showed via invariance theory that $\scA_{\rnull}^{\Y}$ has the property $\nabla^2 g(\Y)[\A] = 0, \forall \A \in \scA_{\rnull}^\Y$ at FOSP $\Y$. Here, by establishing the connection between the Riemannian and the Euclidean Hessians, we can establish the same result without assuming $f$ is convex. Moreover, in the later Theorems \ref{th: general case spec-connect 1} and \ref{th: general case spec-connect 2} we extend our result to the general case and provide explicit expressions for the eigenspace corresponding to the zero eigenvalues there. Interested readers are also referred to the recent survey \cite{zhang2020symmetry} for the discussion on the effect of invariance in geometry of nonconvex problems. 
\end{Remark}

Theorem \ref{th: RHessian-EHessian PSD} immediately shows the following equivalence of FOSPs and SOSPs between the manifold and the factorization formulations in low-rank PSD matrix optimization.

\begin{Corollary}({\bf Equivalence on FOSPs, SOSPs and Strict Saddles of Manifold and Factorization Formulations (PSD Case)}) \label{coro: landscape connection PSD} (a) If $\Y$ is a rank $r$ Euclidean FOSP or SOSP or strict saddle of \eqref{eq: PSD factorization}, then $\X = \Y \Y^\top$ is a Riemannian FOSP or SOSP or strict saddle of \eqref{eq: PSD-manifold-formulation}; (b) if $\X$ is a Riemannian FOSP or SOSP or strict saddle of \eqref{eq: PSD-manifold-formulation}, then any $\Y$ such that $\Y \Y^\top = \X$ is a Euclidean FOSP or SOSP or strict saddle of \eqref{eq: PSD factorization}.
\end{Corollary}

\begin{Remark}({\bf Geometric Landscape Connection Between Manifold and Factorization Formulations on FOSPs, SOSPs and Strict Saddles}) We show in Corollary \ref{coro: landscape connection PSD} that when constraining the Euclidean FOSPs/SOSPs/strict saddles of $g(\Y)$ to be rank $r$, the sets of matrices $\Y \Y^\top$ are exactly the same as the sets of Riemannian FOSPs/SOSPs/strict saddles under the manifold formulation. On the other hand, we note the factorization formulation \eqref{eq: PSD factorization} can have many rank degenerate FOSPs: one canonical example is $\Y = \0$. 
	
	{In addition, we would like to mention that Corollary \ref{coro: landscape connection PSD} can also be obtained via \cite[Proposition 9.6]{boumal2020introduction}. But to our knowledge, that result is included recently after the first preprint of our manuscript. Having said that, our sandwich inequalities established in Theorem \ref{th: RHessian-EHessian PSD} are novel and not covered by theirs. Our sandwich inequalities reveal a finer connection on the spectrum of the Riemannian and the Euclidean Hessians at FOSPs: (1) $\nabla g^2 (\Y)$ has $(r^2-r)/2$ zero eigenvalues with zero eigenspace $\scA_{\rnull}^\Y$; (2) each of the other eigenvalues of $\nabla g^2 (\Y)$ is sandwiched by the corresponding eigenvalues of $\Hess f(\Y \Y^\top)$ with explicit sandwich constants. In Section \ref{sec: application}, we will illustrate the power of these sandwich inequalities in transferring the strict saddle property \citep{ge2015escaping,lee2019first} quantitatively from the factorization formulation to the manifold formulation.}  
\end{Remark}

\begin{Remark}({\bf Implication on Connection of Different Approaches for Rank Constrained Optimization}) \label{rem: implication-on-connection-diff-approaches}
	Broadly speaking, manifold and factorization are two different ways to handle the rank constraint in matrix optimization problems (see also discussion in the Introduction of paper \cite{absil2007trust} on the relationship between manifold optimization and constrained optimization in the Euclidean space). Manifold formulation deals with the rank constraint explicitly via running Riemannian optimization algorithms on the manifold, while the factorization formulation treats the constraint implicitly via factorizing $\X$ into $\Y \Y^\top$ and running the unconstrained optimization algorithms in the Euclidean space. Theorem \ref{th: RHessian-EHessian PSD} establishes a strong geometric landscape connection between two formulations and this provides an example under which the two different approaches are indeed connected in treating the rank constraint.
\end{Remark}

\begin{Remark}
	Currently, the problem \eqref{eq: PSD-manifold-formulation} we considered only has the PSD and rank constraints. Standard SDP problems may have additional linear constraints such as $ \langle \A_i, \X \rangle = b_i $ for $i = 1, \ldots, m$. An interesting research direction is to extend current results to such settings. One strategy to handle these linear constraints is adding a quadratic penalty to the objective and considering solving
	\begin{equation} \label{eq: penalized-PSD-objective}
	 \min_{\X \in \bbS^{p \times p} \succcurlyeq 0, \rank(\X) = r} f(\X) + \mu \sum_{i = 1}^m (\langle \A_i, \X \rangle - b_i  )^2,
	\end{equation} where $\mu > 0$ is a penalty parameter. The landscape of \eqref{eq: penalized-PSD-objective} under the factorization formulation has been considered in \cite{bhojanapalli2018smoothed}. Using our results, one can also transfer the landscape characterization to the corresponding manifold formulation.
\end{Remark}

\section{Geometric Connection of Manifold and Factorization Formulations: General Case} \label{sec: connection-general}
In this section, we present the geometric landscape connection of the manifold formulation \eqref{eq: general prob} and factorization formulations without regularization \eqref{eq: general factor formu} or with regularization \eqref{eq: general factor with reg}. Given $\L\in \bbR^{p_1 \times r}$, $\R \in \bbR^{p_2 \times r}$, suppose that $\X = \L\R^\top$ is of rank $r$ and has SVD $\X = \U \bSigma \V^\top$. Let $\P_1 = \U^\top \L$ and $\P_2 = \V^\top \R$. {Given any $\A = [\A_L^\top \quad \A_R^\top]^\top$ with $\A_L \in \bbR^{p_1\times r}, \A_R \in \bbR^{p_2 \times r}$, define
\begin{equation} \label{def: A-xi-general}
	\xi^{\A}_{\L,\R}:= \L\A_R^\top + \A_L \R^\top = [\U \quad \U_\perp] \begin{bmatrix}
		\P_1 \A_R^\top \V + \U^\top \A_L \P_2^\top & \P_1 \A_R^\top \V_\perp \\
		\U_\perp^\top \A_L \P_2^\top & \0
	\end{bmatrix}[\V \quad \V_\perp]^\top \in T_\X \cM_r;
\end{equation} at the same time, given any $\xi = [\U \quad \U_\perp] \begin{bmatrix}
			\S & \D_2^\top\\
			\D_1 & \0
		\end{bmatrix} [\V \quad \V_\perp]^\top \in T_{\X}\cM_{r}$, define
		\begin{equation}\label{def: xi-setA-general}
			\scA^\xi_{\L,\R} = \{\A = [\A_L^\top \quad\A_R^\top]^\top: \L \A_R^\top + \A_L \R^\top = \xi \}
		\end{equation}
Note that \eqref{def: A-xi-general} and \eqref{def: xi-setA-general} are generalizations of \eqref{def: A-xi-PSD} and \eqref{def: xi-setA-PSD}, respectively and are used to connect the landscape geometry of the manifold formulation \eqref{eq: general prob} and the unregularized factorization formulation \eqref{eq: general factor formu}. To further incorporate the geometry of the regularized formulation \eqref{eq: general factor with reg}, we introduce:
\begin{equation*} 
	\widetilde{\scA}^{\,\,\xi}_{\L,\R} = \{\A = [\A_L^\top \quad\A_R^\top]^\top:  \L \A_R^\top + \A_L \R^\top = \xi\, \text{ and } \, \L^\top \A_L + \A_L^\top \L- \R^\top \A_R -\A_R^\top \R = \0 \}.
\end{equation*}
Compared to $\scA^\xi_{\L,\R}$, there is one additional constraint in the definition of $\widetilde{\scA}^{\,\,\xi}_{\L,\R}$ corresponding to $\nabla^2 g_\reg(\L,\R)$ and it is useful in connecting $\nabla^2 g_\reg(\L,\R)$ with $\nabla^2 g(\L,\R)$ and $\Hess f(\X)$ as we will see in Theorem \ref{th: general case spec-connect 2}. The following lemma shows the affine space $\widetilde{\scA}^{\,\,\xi}_{\L,\R}$ is nonempty and establishes the dimension and some properties of $\scA^{\,\,\xi}_{\L,\R}, \widetilde{\scA}^{\,\,\xi}_{\L,\R}$. 

\begin{Lemma}\label{lm: tangent-vector-equiva-general}
	Given $\L\in \bbR^{p_1 \times r}$, $\R \in \bbR^{p_2 \times r}$, suppose that $\X = \L\R^\top$ is of rank $r$ and has SVD $\U \bSigma \V^\top$. Let $\P_1 = \U^\top \L$ and $\P_2 = \V^\top \R$. Given any $\xi \in T_{\X}\cM_r$, we have $\dim(\scA^\xi_{\L,\R}) = r^2$, $\dim(\widetilde{\scA}^{\,\,\xi}_{\L,\R}) = (r^2-r)/2$,
	\begin{equation}
	\scA^\xi_{\L,\R} := \left\{ \A =\begin{bmatrix}
	\A_L \\
	\A_R
\end{bmatrix}: \begin{array}{c}
	 \A_L = (\U \S_1 + \U_\perp \D_1 ) \P_2^{-\top} \in \bbR^{p_1 \times r} \\
	\A_R = (\V \S_2^\top + \V_\perp \D_2 ) \P_1^{-\top}\in \bbR^{p_2 \times r}
\end{array} \text{ and  } \S_1 + \S_2 = \S   \right \},
\end{equation} and 
\begin{equation} \label{def: xi-setA2-general}
	\widetilde{\scA}^{\,\,\xi}_{\L,\R} := \left\{ \A =\begin{bmatrix}
	\A_L \\
	\A_R
\end{bmatrix}: \begin{array}{l}
	 \A_L = (\U \S_1 + \U_\perp \D_1 ) \P_2^{-\top} \in \bbR^{p_1 \times r}, \\
	 \A_R = (\V \S_2^\top + \V_\perp \D_2 ) \P_1^{-\top} \in \bbR^{p_2 \times r}, \quad \S_1 + \S_2 = \S,\\
	 \P_1^\top \S_1 \P_2^{-\top} + \P_2^{-1} \S_1^\top \P_1 - \P_2^\top \S_2^\top \P_1^{-\top} - \P_1^{-1} \S_2 \P_2 = \0
\end{array}  \right \}.
\end{equation} 
\end{Lemma}
}
Similar to the PSD case discussed in Section \ref{sec: connection-PSD}, we construct $\xi_{\L,\R}^\A, \scA^\xi_{\L,\R}$ and $\widetilde{\scA}^{\,\,\xi}_{\L,\R}$ to find a correspondence between $\bbR^{(p_1 + p_2) \times r}$ and $T_\X\cM_{r}$. On the other hand, we note given $(\L, \R), \A \in \bbR^{(p_1+p_2) \times r}$ and $\xi \in T_{\X} \cM_r$, $\xi^{\A}_{\L,\R}$ is a single matrix while $\scA^\xi_{\L, \R}$ forms a subspace of $\bbR^{(p_1 + p_2) \times r}$ with dimension $r^2$ and $\widetilde{\scA}^{\,\,\xi}_{\L,\R} \subset \scA^\xi_{\L,\R}$ forms a subspace of $\bbR^{(p_1 + p_2) \times r}$ with dimension $ (r^2 - r)/2 $. To deal with this ambiguity, we introduce the following two decompositions for $\bbR^{(p_1 + p_2) \times r}$ tailored to $\scA^\xi_{\L,\R}$ and $\widetilde{\scA}^{\,\,\xi}_{\L,\R}$, respectively. 
\begin{Lemma}\label{lm: decomposition-Rp1p2-general}
	 Under the conditions in Lemma \ref{lm: tangent-vector-equiva-general}, it holds that: 
	 \begin{itemize}[leftmargin=*]
	 	\item $\bbR^{(p_1 + p_2) \times r} = \scA_{\rnull}^{\L,\R} \oplus \scA_{\overline{\rnull}}^{\L,\R}$ with $\dim( \scA_{\rnull}^{\L,\R}) = r^2$, $\dim(\scA_{\overline{\rnull}}^{\L,\R}) = (p_1 + p_2-r)r$ and $\scA_{\rnull}^{\L,\R} \perp \scA_{\overline{\rnull}}^{\L,\R}$, where
	 \begin{equation*}
	 		\begin{split}
	 		\scA_{\rnull}^{\L,\R}&=	\left\{ \A: \A =\begin{bmatrix}
	\U \S \P_2^{-\top} \\
	-\V \S^\top \P_1^{-\top}
\end{bmatrix}, \S \in \bbR^{r \times r}   \right \};\\
\scA_{\overline{\rnull}}^{\L,\R}&=	\left\{\A: \A =\begin{bmatrix}
	(\U \S \P_2 \P_2^\top + \U_\perp \D_1 ) \P_2^{-\top} \\
	(\V \S^\top \P_1 \P_1^\top  + \V_\perp \D_2 ) \P_1^{-\top}
\end{bmatrix}, \begin{array}{l}
	 \D_1 \in \bbR^{(p_1-r) \times r}, \D_2 \in \bbR^{(p_2-r) \times r},
	\S \in \bbR^{r \times r}
\end{array}    \right \}.
	 		\end{split}
	 	\end{equation*}
	 \item Define 
	 \begin{equation} \label{def: space-S-general}
	 	\scS_{\L,\R} = \left\{\S \in \bbR^{r \times r}: \P_1^\top \S \P_2^{-\top} + (\P_1^\top \S \P_2^{-\top})^\top + \P_1^{-1} \S \P_2 + (\P_1^{-1} \S \P_2)^\top = \0 \right\}. 
	 \end{equation} Then $\dim (\scS_{\L,\R}) = (r^2-r)/2$. Moverover, $\bbR^{(p_1 + p_2) \times r} = \widetilde{\scA}_{\rnull}^{\,\,\L,\R} \oplus \widetilde{\scA}_{\overline{\rnull}}^{\,\,\L,\R}$ with $\dim(\widetilde{\scA}_{\rnull}^{\,\,\L,\R}) = (r^2 - r)/2$, $\dim(\widetilde{\scA}_{\overline{\rnull}}^{\,\,\L,\R}) = (p_1+p_2)r - (r^2 - r)/2$ and $\widetilde{\scA}_{\rnull}^{\,\,\L,\R} \perp \widetilde{\scA}_{\overline{\rnull}}^{\,\,\L,\R}$, where
		 \begin{equation*}
	 		\begin{split}
	 			\widetilde{\scA}_{\rnull}^{\,\,\L,\R}&=	\left\{ \A: \A =\begin{bmatrix}
	\U \S \P_2^{-\top} \\
	-\V \S^\top \P_1^{-\top}
\end{bmatrix}, \S \in \scS_{\L,\R}  \right \} \subseteq \scA_{\rnull}^{\L,\R};\\
\widetilde{\scA}_{\overline{\rnull}}^{\,\,\L,\R}&=	\left\{\A: \A =\begin{bmatrix}
	(\U \S_1 \P_2 \P_2^\top + \U_\perp \D_1 ) \P_2^{-\top} \\
	(\V \S_2^\top \P_1 \P_1^\top + \V_\perp \D_2 ) \P_1^{-\top}
\end{bmatrix},\begin{array}{l}
    \D_1 \in \bbR^{(p_1-r) \times r}, \\
     \D_2 \in \bbR^{(p_2-r) \times r}, (\S_1 - \S_2) \perp \scS_{\L,\R}
\end{array}  \right \}.
	 		\end{split}
	 	\end{equation*} 
	 \end{itemize}
\end{Lemma}

\begin{figure}
	\centering
	\includegraphics[width=0.90\textwidth]{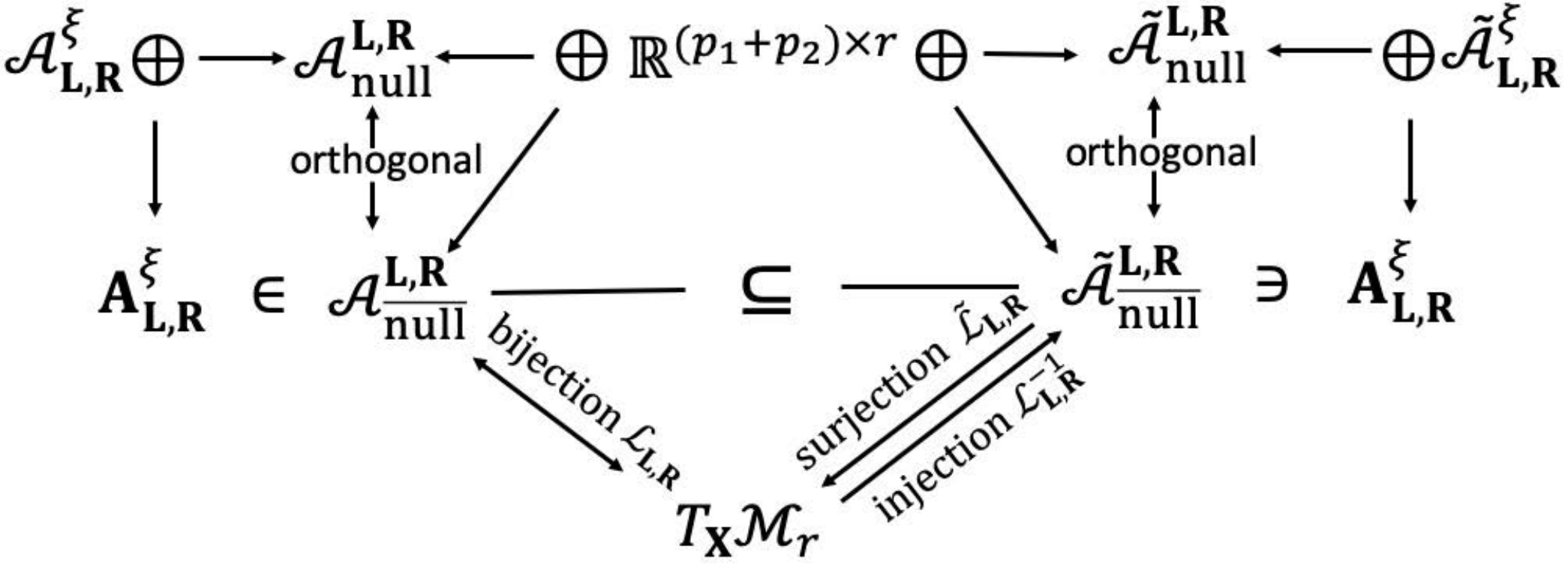}
	\caption{Relationship of subspaces involved in two decompositions in Lemma \ref{lm: decomposition-Rp1p2-general}. Left hand side: first decomposition in Lemma \ref{lm: decomposition-Rp1p2-general} on relationship between $\bbR^{(p_1+p_2) \times r}$, $T_\X \cM_{r}$, $\scA_{\rnull}^{\L,\R}$, $\scA_{\overline{\rnull}}^{\L,\R}$, $\A_{\L,\R}^\xi$, and $\scA_{\L,\R}^\xi$; Right hand side: second decomposition in Lemma \ref{lm: decomposition-Rp1p2-general} on relationship between $\bbR^{(p_1+p_2) \times r}$, $T_\X \cM_{r}$, $\widetilde{\scA}_{\rnull}^{\,\,\L,\R}$, $\widetilde{\scA}_{\overline{\rnull}}^{\,\,\L,\R}$, $\A_{\L,\R}^{\xi}$, and $\widetilde{\scA}_{\L,\R}^{\,\,\xi}$. }
	\label{fig: general_decom_ill}
\end{figure}

As we will see in Theorems \ref{th: general case spec-connect 1} and \ref{th: general case spec-connect 2}, $ \scA_{\rnull}^{\L,\R}$ and $\widetilde{\scA}_{\rnull}^{\,\,\L,\R}$ correspond to the eigenspace of zero eigenvalues of $\nabla^2 g(\L,\R)$ and $\nabla^2 g_\reg(\L,\R)$, respectively. By the decompositions in Lemma \ref{lm: decomposition-Rp1p2-general}, we derive the following three results in Proposition \ref{prop: bijection-general}: first, we show $\scA_{\L,\R}^\xi$ and $\widetilde{\scA}_{\L,\R}^{\,\,\xi}$ can be further decomposed as the direct sum of $\scA_{\rnull}^{\L,\R}$ and $\widetilde{\scA}_{\rnull}^{\,\,\L,\R}$ with the same single matrix from $\scA_{\overline{\rnull}}^{\L,\R}$; second, we find there exists a bijection between $\scA_{\overline{\rnull}}^{\L,\R}$ and $T_\X \cM_r$; third, we construct a ``pseudobijection'' between $\widetilde{\scA}_{\overline{\rnull}}^{\,\,\L,\R}$ and $T_\X \cM_r$ since a bijection between them is impossible due to the mismatch of their dimensions. A pictorial illustration of the relationship of subspaces in Lemma \ref{lm: decomposition-Rp1p2-general} is given in Figure \ref{fig: general_decom_ill}.

\begin{Proposition}[Decompositions of $\scA_{\L,\R}^\xi$ and $\widetilde{\scA}_{\L,\R}^{\,\,\xi}$ and Maps Between $\scA_{\overline{\rnull}}^{\L,\R}$, $\widetilde{\scA}_{\overline{\rnull}}^{\,\,\L,\R}$ and $T_\X \cM_r$] \label{prop: bijection-general} Under the conditions in Lemma \ref{lm: tangent-vector-equiva-general}, let $\xi = [\U \quad \U_\perp] \begin{bmatrix}
			\S & \D_2^\top\\
			\D_1 & \0
		\end{bmatrix} [\V \quad \V_\perp]^\top \in T_{\X}\cM_{r}$. Then
\begin{equation} \label{eq: A-LR-xi-general-no-reg}
\begin{array}{c}
	\scA_{\L,\R}^\xi= \A_{\L,\R}^{\xi} \oplus \scA_{\rnull}^{\L,\R} \\
	\widetilde{\scA}_{\L,\R}^{\,\,\xi}= \A_{\L,\R}^\xi \oplus \widetilde{\scA}_{\rnull}^{\,\,\L,\R}
\end{array}
    ~~ \text{for}~~\A_{\L,\R}^\xi = \begin{bmatrix}
		(\U\widebar{\S}\P_2\P_2^\top + \U_\perp\D_1)\P_2^{-\top}\\
		(\V\widebar{\S}^\top\P_1\P_1^\top + \V_\perp\D_2)\P_1^{-\top}
		\end{bmatrix} \in \scA_{\overline{\rnull}}^{\L,\R} \subseteq \widetilde{\scA}_{\overline{\rnull}}^{\,\,\L,\R},
\end{equation}	where $\widebar{\S}$ is the unique solution of the following Sylvester equation $\P_2 \P_2^\top \widebar{\S}^\top   +  \widebar{\S}^\top\P_1 \P_1^\top  = \S^\top$. 
	
	Moreover, there is a bijective linear map $\cL_{\L,\R}$ between $\scA_{\overline{\rnull}}^{\L,\R}$ and $T_\X \cM_r$ given as follows
	\begin{equation} \label{eq: bijection-general-case}
		\cL_{\L,\R}: \A \in \scA_{\overline{\rnull}}^{\L,\R} \longrightarrow \xi_{\L,\R}^\A \in T_\X \cM_r \quad  \text{ and } \quad \cL^{-1}_{\L,\R}: \xi \in T_\X \cM_r\longrightarrow \A_{\L,\R}^\xi \in \scA_{\overline{\rnull}}^{\L,\R}.
	\end{equation} In addition, there is a surjective linear map $\widetilde{\cL}_{\L,\R}$ between $\widetilde{\scA}_{\overline{\rnull}}^{\,\,\L,\R}$ and $T_\X \cM_r$ given as follows
	\begin{equation} \label{def: linear-maps-general-reg}
		\widetilde{\cL}_{\L,\R}: \A \in \widetilde{\scA}_{\overline{\rnull}}^{\,\,\L,\R} \longrightarrow \xi_{\L,\R}^\A \in T_\X \cM_r,
	\end{equation} and it satisfies $\widetilde{\cL}_{\L,\R}(\cL^{-1}_{\L,\R}(\xi) ) = \xi$ for any $\xi \in T_\X \cM_r$.
	
	Finally, we have the following spectrum bounds for $\cL_{\L,\R}$, $\cL^{-1}_{\L,\R}$ and $\widetilde{\cL}_{\L,\R}$, respectively:
	\begin{equation} \label{ineq: bijection-spectrum-bound-general1}
	\begin{split}
		( \sigma_r(\L) \wedge \sigma_r(\R) )^2 \|\A\|_\F^2 \leq &\|\cL_{\L,\R}(\A)\|_\F^2 \leq 2 (\sigma_1(\L) \vee \sigma_1(\R) )^2 \|\A\|_\F^2, \quad \forall \A \in \scA_{\overline{\rnull}}^{\L,\R}, \\
		(\sigma_1(\L) \vee \sigma_1(\R)  )^{-2} \|\xi\|_\F^2/2 \leq &\|\cL^{-1}_{\L,\R}(\xi)\|_\F^2 \leq (\sigma_r(\L) \wedge \sigma_r(\R) )^{-2} \|\xi\|_\F^2, \quad  \forall \xi \in T_\X \cM_r,\\
		&\|\widetilde{\cL}_{\L,\R}(\A)\|_\F^2 \leq 2( \sigma_1(\L) \vee \sigma_1(\R) )^2 \|\A\|_\F^2, \quad  \forall \A \in \widetilde{\scA}_{\overline{\rnull}}^{\,\,\L,\R}.
	\end{split}
	\end{equation}
\end{Proposition}
{\noindent \bf Proof of Proposition \ref{prop: bijection-general}.} We divide the proof into two steps: in Step 1, we prove the decomposition results for $\scA_{\L,\R}^\xi$ and $\widetilde{\scA}_{\L,\R}^{\,\,\xi}$; in Step 2, we show $\cL_{\L,\R}$ is a bijection, $ \widetilde{\cL}_{\L,\R}(\cL^{-1}_{\L,\R}(\xi) ) = \xi $, and prove their spectrum bounds.

{\bf Step 1}. First, the uniqueness of $\widebar{\S}$ is guaranteed by the fact $\P_1 \P_1^\top$ and $- \P_2 \P_2^\top$ have disjoint spectra and \cite[Theorem VII.2.1]{bhatia2013matrix}. Next, we prove $\scA_{\L,\R}^\xi= \A_{\L,\R}^{\xi} \oplus \scA_{\rnull}^{\L,\R}$. Recall $\A_{\L,\R}^{\xi}  = \begin{bmatrix}
		(\U\widebar{\S}\P_2\P_2^\top + \U_\perp\D_1)\P_2^{-\top}\\
		(\V\widebar{\S}^\top\P_1\P_1^\top + \V_\perp\D_2)\P_1^{-\top}
		\end{bmatrix} \in \scA_{\overline{\rnull}}^{\L,\R}$, given $\A = \begin{bmatrix}
		             (\U \S_1 + \U_\perp \D_1 ) \P_2^{-\top}\\
		              (\V \S_2^\top + \V_\perp \D_2 ) \P_1^{-\top}
		\end{bmatrix}  \in \scA_{\L,\R}^\xi$, we have
		\begin{equation} \label{eq: general-bijection-eq1}
		    \A - \A_{\L,\R}^\xi = \begin{bmatrix}
		       \U(\S_1-\widebar{\S}\P_2\P_2^\top )\P_2^{-\top} \\
		       \V(\S_2^\top-\widebar{\S}^\top\P_1\P_1^\top)\P_1^{-\top}
		\end{bmatrix} \overset{(a)}=  \begin{bmatrix}
		       \U(- \S_2+\P_1 \P_1^\top \widebar{\S})\P_2^{-\top} \\
		       \V(\S_2^\top-\widebar{\S}^\top\P_1\P_1^\top)\P_1^{-\top}
		\end{bmatrix}\in \scA_{\rnull}^{\L,\R},
		\end{equation} where (a) is because $\widebar{\S} \P_2 \P_2^\top + \P_1 \P_1^\top \widebar{\S} = \S$ and $\S_1 + \S_2 = \S$. Moreover, $\A_{\L,\R}^\xi + \A \in \scA_{\L,\R}^\xi$ for any $\A \in \scA_{\rnull}^{\L,\R}$. This proves $ \scA_{\L,\R}^\xi= \A_{\L,\R}^{\xi} \oplus \scA_{\rnull}^{\L,\R}$. 
		
		Next, we prove the second decomposition result $\widetilde{\scA}_{\L,\R}^{\,\,\xi}= \A_{\L,\R}^\xi \oplus \widetilde{\scA}_{\rnull}^{\,\,\L,\R}$. Given $\A_{\L,\R}^{\xi}  = \begin{bmatrix}
		(\U\widebar{\S}\P_2\P_2^\top + \U_\perp\D_1)\P_2^{-\top}\\
		(\V\widebar{\S}^\top\P_1\P_1^\top + \V_\perp\D_2)\P_1^{-\top}
		\end{bmatrix}$ and $\A = \begin{bmatrix}
		             (\U \S_1 + \U_\perp \D_1 ) \P_2^{-\top}\\
		              (\V \S_2^\top + \V_\perp \D_2 ) \P_1^{-\top}
		\end{bmatrix}  \in \widetilde{\scA}_{\L,\R}^{\,\,\xi}$, we have
		\begin{equation} \label{eq: check-scS-LR}
		    \begin{split}
		       &\P_1^\top(\P_1 \P_1^\top \widebar{\S} - \S_2) \P_2^{-\top} + (\P_1^\top(\P_1 \P_1^\top \widebar{\S} - \S_2) \P_2^{-\top})^\top \\
		       & + \P_1^{-1}(\P_1 \P_1^\top \widebar{\S} - \S_2) \P_2 + (\P_1^{-1}(\P_1 \P_1^\top \widebar{\S} - \S_2) \P_2)^\top\\
		       \overset{(a)}= & \P_1^\top (\S - \widebar{\S} \P_2 \P_2^\top ) \P_2^{-\top} + \P_2^{-1}(\S^\top - \P_2 \P_2^\top \widebar{\S}^\top ) \P_1 + \P_1^\top \widebar{\S} \P_2 + \P_2^\top \widebar{\S}^\top \P_1 \\
		       &- \P_1^\top \S_2 \P_2^{-\top} - (\P_1^\top \S_2 \P_2^{-\top})^\top - \P_1^{-1} \S_2 \P_2 - (\P_1^{-1} \S_2 \P_2)^\top \\
		       = &  \P_1^\top \S\P_2^{-\top} + \P_2^{-1} \S^\top \P_1 - \P_1^\top \S_2 \P_2^{-\top} - (\P_1^\top \S_2 \P_2^{-\top})^\top - \P_1^{-1} \S_2 \P_2 - (\P_1^{-1} \S_2 \P_2)^\top\\
		      \overset{(b)} =  & \P_1^\top \S\P_2^{-\top} + \P_2^{-1} \S^\top \P_1  - \P_1^\top (\S-\S_1) \P_2^{-\top} - (\P_1^\top (\S-\S_1) \P_2^{-\top})^\top - \P_1^{-1} \S_2 \P_2 - (\P_1^{-1} \S_2 \P_2)^\top \\
		      = & \P_1^\top \S_1  \P_2^{-\top} + (\P_1^\top \S_1 \P_2^{-\top})^\top - \P_1^{-1} \S_2 \P_2 - (\P_1^{-1} \S_2 \P_2)^\top \overset{(c)}= \0.
		    \end{split}
		\end{equation} Here (a) is because $\widebar{\S} \P_2 \P_2^\top + \P_1 \P_1^\top \widebar{\S} = \S$, (b) is because $\S_1 + \S_2 = \S$ and (c) is by the constraints of $\S_1, \S_2$ given in $\widetilde{\scA}_{\L,\R}^{\,\,\xi}$. So \eqref{eq: check-scS-LR} shows $\P_1 \P_1^\top \widebar{\S} - \S_2 \in \scS_{\L,\R}$ and we have
		\begin{equation*}
		   \A - \A_{\L,\R}^\xi \overset{ \eqref{eq: general-bijection-eq1} }=  \begin{bmatrix}
		       \U(-\S_2 + \P_1 \P_1^\top \widebar{\S})\P_2^{-\top} \\
		       \V(\S_2^\top-\widebar{\S}^\top\P_1\P_1^\top )\P_1^{-\top}
		\end{bmatrix} \overset{\eqref{eq: check-scS-LR}}\in \widetilde{\scA}_{\rnull}^{\,\,\L,\R}.
		\end{equation*} 
        Moreover, for $\A_{\L,\R}^\xi$ and any $\A = \begin{bmatrix}
	\U \S \P_2^{-\top} \\
	-\V \S^\top \P_1^{-\top}
\end{bmatrix} \in \widetilde{\scA}_{\rnull}^{\,\,\L,\R}$ with $\S \in \scS_{\L,\R}$, we have
\begin{equation} \label{eq: check-scA-LR}
    \begin{split}
        &(\widebar{\S}\P_2\P_2^\top + \S) + (\widebar{\S}^\top\P_1\P_1^\top - \S^\top)^\top \overset{(a)}= \S; \\
        &\P_1^\top (\widebar{\S}\P_2\P_2^\top + \S) \P_2^{-\top} + \P_2^{-1} (\P_2\P_2^\top \widebar{\S}^\top + \S^\top) \P_1 - \P_2^\top (\widebar{\S}^\top \P_1\P_1^\top - \S^\top ) \P_1^{-\top} - \P_1^{-1} (\P_1\P_1^\top \widebar{\S} - \S)\P_2 \\
        = & \P_1^\top \S \P_2^{-\top} + (\P_1^\top \S \P_2^{-\top})^\top + \P_1^{-1} \S \P_2 + (\P_1^{-1} \S \P_2)^\top \overset{(b)}= \0.
    \end{split}
\end{equation} Here (a) is because $\widebar{\S} \P_2 \P_2^\top + \P_1 \P_1^\top \widebar{\S} = \S$ and (b) is because $\S \in \scS_{\L,\R}$. Thus,
\begin{equation*}
    \A_{\L,\R}^\xi + \A =  \begin{bmatrix}
		(\U(\widebar{\S}\P_2\P_2^\top + \S) + \U_\perp\D_1)\P_2^{-\top}  \\
		(\V(\widebar{\S}^\top\P_1\P_1^\top - \S^\top) + \V_\perp\D_2)\P_1^{-\top}
		\end{bmatrix} \overset{ \eqref{eq: check-scA-LR} }\in \widetilde{\scA}_{\L,\R}^{\,\,\xi}.
\end{equation*} This finishes the proof for $\widetilde{\scA}_{\L,\R}^{\,\,\xi}= \A_{\L,\R}^\xi \oplus \widetilde{\scA}_{\rnull}^{\,\,\L,\R}$.

{\bf Step 2}. We begin by proving $\cL_{\L,\R}$ is a bijection. Note that both $ \scA_{\overline{\rnull}}^{\L,\R}$ and $T_\X\cM_{r}$ have dimension $(p_1 + p_2 -r)r$. Suppose $\cL_{\L,\R}' : \xi \in T_\X\cM_{r} \longrightarrow \A_{\L,\R}^\xi \in \scA_{\overline{\rnull}}^{\L,\R}$. Then for any $\xi = [\U \quad \U_\perp] \begin{bmatrix}
			\S & \D_2^\top\\
			\D_1 & \0
		\end{bmatrix} [\V \quad \V_\perp]^\top \in T_\X\cM_{r}$, we have
\begin{equation} \label{eq: bijection-general}
\begin{split}
	& \cL_{\L,\R} ( \cL_{\L,\R}'(\xi) ) = \cL_{\L,\R} (\A_{\L,\R}^\xi) 
	=  [\U \quad \U_\perp] \begin{bmatrix}
	\widebar{\S} \P_2 \P_2^\top + \P_1 \P_1^\top \widebar{\S} & \D_2^\top \\
		\D_1 & \0
	\end{bmatrix}[\V \quad \V_\perp]^\top = \xi.
\end{split}
\end{equation} Since $\cL_{\L,\R}$ and $\cL_{\L,\R}'$ are linear maps, \eqref{eq: bijection-general} implies $\cL_{\L,\R}$ is bijection and $\cL_{\L,\R}' = \cL_{\L,\R}^{-1}$. Following a similar proof of \eqref{eq: bijection-general}, we can also show $\widetilde{\cL}_{\L,\R}(\cL^{-1}_{\L,\R}(\xi) ) = \xi $ holds for any $\xi \in T_\X\cM_{r}$ and this implies $\widetilde{\cL}_{\L,\R}$ is surjective.

Next, we provide the spectrum bounds for $\cL_{\L,\R}$. Suppose $\A = [\A_L^\top \quad \A_R^\top]^\top \in  \scA_{\overline{\rnull}}^{\L,\R}$, where $\A_L =  (\U \S \P_2 \P_2^\top + \U_\perp \D_1 ) \P_2^{-\top}$, $\A_R = (\V \S^\top \P_1 \P_1^\top + \V_\perp \D_2 ) \P_1^{-\top}$. Then
\begin{equation}\label{ineq: bijection-general-ineq1}
	\begin{split}
		\|\A\|_\F^2 &\leq (\|\S \P_2 \P_2^\top\|_\F^2 + \|\D_1\|_\F^2) \sigma^2_1(\P_2^{-\top}) + (\| \P_1 \P_1^\top \S\|_\F^2 + \|\D_2\|_\F^2) \sigma^2_1(\P_1^{-\top}) \\
		&\overset{(a)}\leq (\|\S \P_2 \P_2^\top\|_\F^2 + \|\P_1 \P_1^\top \S\|_\F^2 + \|\D_1\|_\F^2 + \|\D_2\|_\F^2 )/(\sigma_r(\L) \wedge \sigma_r(\R) )^2;
	\end{split}
\end{equation} where in (a), we use the fact $\L= \U \P_1, \R = \V \P_2$ and $\L,\R$ share the same spectrum as $\P_1,\P_2$. In addition,
\begin{equation}\label{ineq: bijection-general-ineq2}
	\langle \P_1 \P_1^\top \S, \S \P_2 \P_2^\top \rangle  = \langle (\P_1 \P_1^\top)^{1/2} \S (\P_2 \P_2^\top)^{1/2},(\P_1 \P_1^\top)^{1/2} \S (\P_2 \P_2^\top)^{1/2}\rangle \geq 0.
\end{equation}  So
\begin{equation*}
	\begin{split}
		\|\cL_{\L,\R}(\A)\|_\F^2 = \|\xi^\A_{\L,\R}\|_\F^2 &\overset{ \eqref{def: A-xi-general} } = \|\P_1 \A_R^\top \V + \U^\top \A_L \P_2^\top\|_\F^2 + \|\U_\perp^\top \A_L \P_2^\top\|_\F^2 + \| \P_1 \A_R^\top \V_\perp\|_\F^2\\
		& = \|\P_1 \P_1^\top \S + \S \P_2 \P_2^\top\|_\F^2  + \|\D_1\|_\F^2 + \|\D_2\|_\F^2\\
		& \overset{ \eqref{ineq: bijection-general-ineq2} } \geq  \|\S \P_2 \P_2^\top\|_\F^2 + \|\P_1 \P_1^\top \S\|_\F^2 + \|\D_1\|_\F^2 + \|\D_2\|_\F^2 \overset{ \eqref{ineq: bijection-general-ineq1} }\geq (\sigma_r(\L) \wedge \sigma_r(\R) )^2\|\A\|_\F^2,
	\end{split}
\end{equation*} and
\begin{equation*}
\begin{split}
	\|\cL_{\L,\R}(\A)\|_\F^2 =\|\xi^\A_{\L,\R}\|_\F^2 &\overset{ \eqref{def: A-xi-general} } = \|\P_1 \A_R^\top \V + \U^\top \A_L \P_2^\top\|_\F^2 + \|\U_\perp^\top \A_L \P_2^\top\|_\F^2 + \| \P_1 \A_R^\top \V_\perp\|_\F^2\\
	 & \leq 2(\|\P_1 \A_R^\top \V \|_\F^2 + \|\U^\top \A_L \P_2^\top \|_\F^2 ) + \|\U_\perp^\top \A_L \P_2^\top\|_\F^2 + \| \P_1 \A_R^\top \V_\perp\|_\F^2\\
	 & = \|\P_1 \A_R^\top \V \|_\F^2 + \|\U^\top \A_L \P_2^\top \|_\F^2 + \|\A_L\P_2^\top\|_\F^2 + \|\P_1 \A_R^\top\|_\F^2\\
	 & \leq 2( \sigma_1(\L) \vee \sigma_1(\R) )^2 \|\A\|_\F^2.
\end{split}
\end{equation*} By the relationship of the spectrum of an operator and its inverse, the spectrum bounds for $\cL_{\L,\R}^{-1}$ follow from the ones of $\cL_{\L,\R}$. Finally, since $ \widetilde{\cL}_{\L,\R}$ is surjective and the ``pseudoinverse'' of $\cL^{-1}_{\L,\R}$, its spectrum upper bound follows from the spectrum lower bound of $\cL^{-1}_{\L,\R}$. This finishes the proof of this proposition. \quad $\blacksquare$

Now, we are ready to present our main results on the geometric landscape connection of the manifold and the factorization formulations in the general low-rank matrix optimization. 
\begin{Theorem}[Geometric Landscape Connection of Formulations \eqref{eq: general prob} and \eqref{eq: general factor formu}] \label{th: general case spec-connect 1}
	 Suppose $\L\in \bbR^{p_1 \times r}$, $\R \in \bbR^{p_2 \times r}$ and $\X = \L\R^\top$ are of rank $r$. Then
	\begin{equation} \label{eq: gradient-connection-general}
	\begin{split}
		\grad f(\X) = \nabla_\L g(\L,\R) \R^\dagger + (\nabla_\R g(\L,\R) \L^\dagger)^\top(\I_{p_2} - \R \R^\dagger) \quad \text{and} \quad
		\nabla g(\L,\R) = \left[\begin{array}{c}
			\grad f(\X) \R\\
			(\grad f(\X))^\top \L
		\end{array} \right].
	\end{split}
	\end{equation}
	
Furthermore, if $(\L,\R)$ is a Euclidean FOSP \eqref{eq: general factor formu}, then we have
\begin{equation} \label{eq: R-E-Hessian-general-1}
	\begin{split}
	\nabla^2 g(\L,\R)[\A,\A] &=  \Hess f(\X)[\xi_{\L,\R}^\A, \xi_{\L,\R}^\A], \quad \forall \A \in \bbR^{(p_1 + p_2) \times r}; \\
		\Hess f(\X)[\xi,\xi] &= \nabla^2 g({\L,\R})[\cL_{\L,\R}^{-1}(\xi) ,\cL_{\L,\R}^{-1}(\xi)], \quad \forall \xi \in T_\X\cM_{r}.
	\end{split}
	\end{equation}
	More precisely,
\begin{equation} \label{eq: R-E-Hessian-general-2}
	\begin{split}
		\nabla^2 g(\L,\R)[\A] &= \0, \quad \forall \A \in \scA_{\rnull}^{\L,\R};\\
		\nabla^2 g(\L,\R)[\A,\A] &= \Hess f(\X)[\cL_{\L,\R} (\A),\cL_{\L,\R} (\A) ],\quad \forall \A \in \scA_{\overline{\rnull}}^{\L,\R}.
	\end{split}
	\end{equation}	
	
	 Finally, $\Hess f(\X)$ has $(p_1 + p_2 -r)r $ eigenvalues and $ \nabla^2 g({\L,\R})$ has $(p_1 + p_2)r$ eigenvalues. $\nabla^2 g({\L,\R})$ has at least $r^2$ zero eigenvalues 
	 with the corresponding zero eigenspace $\scA_{\rnull}^{\L,\R}$. Denote the rest of the $(p_1+p_2-r)r $ possibly non-zero eigenvalues of $\nabla^2 g({\L,\R})$ from the largest to the smallest as $\widebar{\lambda}_1(\nabla^2 g({\L,\R})),\ldots, \widebar{\lambda}_{(p_1+p_2-r)r}(\nabla^2 g({\L,\R}))$. Then for $i=1,\ldots,(p_1+p_2-r)r$: $\widebar{\lambda}_i(\nabla^2 g({\L,\R}))$ is sandwiched between $( \sigma_r(\L) \wedge \sigma_r(\R) )^2 \lambda_i(\Hess f(\X)) $ and $2 (\sigma_1(\L) \vee \sigma_1(\R) )^2 \lambda_i(\Hess f(\X)) $.
\end{Theorem}
{\noindent \bf Proof of Theorem \ref{th: general case spec-connect 1}.}  First, suppose $\X = \L \R^\top$ has SVD $\U \bSigma \V^\top$. Then $\L \L^\dagger = P_\U, \R\R^\dagger = P_\V$ and $\L$ and $\R$ lie in the column spaces of $\U$ and $\V$, respectively. So \eqref{eq: gradient-connection-general} is by direct calculation from the expressions of Riemannian and Euclidean gradients given in Proposition \ref{prop: gradient-exp}. The rest of the proof is divided into two steps: in Step 1, we prove \eqref{eq: R-E-Hessian-general-1} and \eqref{eq: R-E-Hessian-general-2}; in Step 2, we prove the individual eigenvalue connection between $\Hess f(\X)$ and $\nabla^2 g(\L,\R)$.

{\bf Step 1.} We begin by showing the first equality in \eqref{eq: R-E-Hessian-general-1}. Let $\P_1 = \U^\top \L,\P_2 = \V^\top \R$, it is easy to verify $\P_1 \P_2^\top = \U^\top \L\R^\top \V = \bSigma$.
		Since $(\L,\R)$ is a Euclidean FOSP of \eqref{eq: general factor formu}, by \eqref{eq: gradient-connection-general} we have $\X$ is a Riemannian FOSP of \eqref{eq: general prob}. So $\nabla f(\X) = P_{\U_\perp} \nabla f(\X) P_{\V_\perp}$. Given $\A = [\A_L^\top \quad \A_R^\top]^\top \in \bbR^{(p_1 + p_2) \times r}$, we have
\begin{equation} \label{eq: Hessian-connnect-gradient-1}
\begin{split}
	\langle \nabla f(\X), \U_\perp \U_\perp^\top \A_L \P_2^\top \bSigma^{-1} \P_1 \A_R^\top \V_\perp  \V_\perp^\top \rangle \overset{(a)}=  \langle \nabla f(\X), \A_L \A_R^\top \rangle.
\end{split}
\end{equation} Here $(a)$ is because $\nabla f(\X) = P_{\U_\perp} \nabla f(\X) P_{\V_\perp}$ and $\P_1 \P_2^\top = \bSigma$.

Then by Proposition \ref{prop: Hessian-exp},
\begin{equation} \label{eq: Hessian-connect-hessian}
	\begin{split}
		\nabla^2 g(\L,\R)[\A,\A] &= 2 \langle \nabla f(\L\R^\top), \A_L \A_R^\top \rangle + \nabla^2 f(\L\R^\top)[\L\A_R^\top + \A_L \R^\top, \L\A_R^\top + \A_L \R^\top]\\
		& \overset{ \text{Lemma } \ref{lm: tangent-vector-equiva-general}, \eqref{eq: Hessian-connnect-gradient-1} }=  2\langle \nabla f(\X), \U_\perp \U_\perp^\top \A_L \P_2^\top \bSigma^{-1} \P_1 \A_R^\top \V_\perp  \V_\perp^\top \rangle + \nabla^2 f(\X)[\xi_{\L,\R}^\A, \xi_{\L,\R}^\A]\\
		& = \Hess f(\X)[\xi_{\L,\R}^\A, \xi_{\L,\R}^\A],
	\end{split}
\end{equation} 
where the last equality follows from the expressions of $\Hess f(\X)$  and $\xi_{\L,\R}^\A$ in \eqref{eq: Hessian-general} and \eqref{def: A-xi-general}, respectively.
This finishes the proof for the first equality in \eqref{eq: R-E-Hessian-general-1}. Meanwhile, by a similar argument as \eqref{eq: zero-eigenspace-PSD-argument}, we have
\begin{equation*}
	\nabla^2 g(\L,\R)[\A,\B] = \Hess f(\X)[\xi_{\L,\R}^{\A},\xi_{\L,\R}^{\B} ] = 0, \quad  \forall \A \in \scA_{\rnull}^{\L,\R}, \forall \B \in \bbR^{(p_1 + p_2) \times r}.
\end{equation*} This implies the first equality in \eqref{eq: R-E-Hessian-general-2}. 

The second equality in \eqref{eq: R-E-Hessian-general-2} follows directly from the first equality in \eqref{eq: R-E-Hessian-general-1} and the definition of $\cL_{\L,\R}$. Finally, by the bijectivity of $\cL_{\L,\R}$, the second equality in \eqref{eq: R-E-Hessian-general-1} follows from the second equality in \eqref{eq: R-E-Hessian-general-2}. 

{\bf Step 2.} $\Hess f(\X)$ and $\nabla^2 g({\L,\R})$ are by definition linear maps from $T_{\X} \cM_{r}$ and $\bbR^{(p_1 + p_2) \times r}$ to $T_{\X} \cM_{r}$ and $\bbR^{(p_1 + p_2) \times r}$, respectively. Because $T_{\X} \cM_{r}$ is of dimension $(p_1 +p_2-r)r$, the number of eigenvalues of $\Hess f(\X)$ and $\nabla^2 g({\L,\R})$ are $(p_1 +p_2-r)r$ and $(p_1+p_2)r$, respectively. By the first equality in \eqref{eq: R-E-Hessian-general-2}, we have $\scA_{\rnull}^{\L,\R}$ is the eigenspace of $r^2$ zero eigenvalues and the rest of the $(p_1+p_2-r)r$ possibly non-zero eigenvalues of $\nabla^2 g({\L,\R})$ span the eigenspace $\scA_{\overline{\rnull}}^{\L,\R}$. Restricting to $\scA_{\overline{\rnull}}^{\L,\R}$ and $T_\X \cM_{r}$ and using \eqref{eq: R-E-Hessian-general-2}, \eqref{ineq: bijection-spectrum-bound-general1} and Lemma \ref{lm: two-symmetric-matrix-spectrum-connection} in the Appendix, we have $\widebar{\lambda}_i(\nabla^2 g(\L,\R))$ is sandwiched between $( \sigma_r(\L) \wedge \sigma_r(\R) )^2 \lambda_i(\Hess f(\X)) $ and $2 (\sigma_1(\L) \vee \sigma_1(\R) )^2 \lambda_i(\Hess f(\X)) $. This finishes the proof. \quad $\blacksquare$ 

\begin{Theorem}[Geometric Landscape Connection of Formulations \eqref{eq: general prob} and \eqref{eq: general factor with reg}]\label{th: general case spec-connect 2}
	Suppose $\L\in \bbR^{p_1 \times r}, \R \in \bbR^{p_2 \times r}$ and $ \X = \L\R^\top$ are of rank $r$ and $(\L,\R)$ is a Euclidean FOSP of \eqref{eq: general factor with reg}. First, we have
	\begin{equation} \label{eq: reg-FOSP-property}
	\L^\top \L= \R^\top \R \quad \text{ and } \quad \nabla_\reg g(\L,\R) = \nabla g(\L,\R),
\end{equation} and for any  $\A = [\A_L^\top \quad \A_R^\top]^\top \in \bbR^{(p_1 + p_2) \times r}$, 
\begin{equation} \label{eq: reg-Hessian-on-FOSP}
	\nabla^2 g_\reg (\L, \R)[\A, \A] = \nabla^2 g (\L, \R)[\A, \A] + \mu \|\L^\top \A_L + \A_L ^\top \L- \R^\top \A_R - \A_R^\top \R \|_\F^2.
\end{equation}

	Second,
	\begin{equation} \label{eq: R-E-Hessian-general-reg}
		\begin{split}
			\nabla^2 g_\reg(\L,\R)[\A,\A] &= \Hess f(\X)[\xi_{\L,\R}^\A, \xi_{\L,\R}^\A] \\
			& \quad \quad + \mu \|\L^\top \A_L + \A_L^\top \L- \R^\top \A_R - \A_R^\top \R \|_\F^2, \quad \forall \A \in \bbR^{(p_1 + p_2) \times r};\\
			\Hess f(\X)[\xi, \xi] & 	=  \nabla^2 g_\reg(\L,\R)[\cL_{\L,\R}^{-1}(\xi),\cL_{\L,\R}^{-1}(\xi)], \quad \forall \xi \in T_\X \cM_r,
		\end{split}
	\end{equation} where $\cL_{\L,\R}^{-1}$ is the bijective map given in \eqref{eq: bijection-general-case}. More precisely,
	\begin{equation} \label{eq: R-E-Hessian-general-reg-2}
		\begin{split}
			\nabla^2g_\reg(\L,\R)[\A]= & \0, \quad \forall \A \in 	\widetilde{\scA}_{\rnull}^{\,\,\L,\R};\\
			\nabla^2g_\reg(\L,\R)[\A,\A]= & \Hess f(\X)[\widetilde{\cL}_{\L,\R}(\A),\widetilde{\cL}_{\L,\R}(\A)] \\
			& + \mu \|\L^\top \A_L + \A_L ^\top \L- \R^\top \A_R - \A_R^\top \R \|_\F^2, \quad \forall \A \in \widetilde{\scA}_{\overline{\rnull}}^{\,\, \L,\R}.
		\end{split}
	\end{equation}

Finally, $\nabla^2 g_{\reg}(\L,\R)$ has $(p_1 + p_2)r$ eigenvalues and at least $(r^2 -r)/2$ of them are zero spanning the eigenspace $\widetilde{\scA}_{\rnull}^{\,\,\L,\R}$. Denote the rest of the $((p_1+p_2)r-(r^2-r)/2) $ possibly non-zero eigenvalues of $\nabla^2 g_\reg({\L,\R})$ from the largest to the smallest as $\widebar{\lambda}_1(\nabla^2 g_\reg({\L,\R})),\ldots, \widebar{\lambda}_{(p_1+p_2)r-(r^2-r)/2}(\nabla^2 g_\reg({\L,\R}))$. Then
\begin{itemize}[leftmargin=*]
	\item the following lower bounds for $\widebar{\lambda}_i(\nabla^2 g_\reg({\L,\R}))$ hold: \begin{equation} \label{ineq: spectrum-bound-regularized-1}
	\begin{split}
		\widebar{\lambda}_i(\nabla^2 g_\reg({\L,\R})) &\geq \left( \sigma_r(\X) \lambda_i(\Hess f(\X))\right) \wedge \left(2\sigma_1(\X) \lambda_i(\Hess f(\X))\right) , \text{ for } 1 \le i \le (p_1 +p_2-r)r,  \\
		\widebar{\lambda}_{i}(\nabla^2 g_\reg({\L,\R})) &\geq 2\sigma_1(\X) \lambda_{\min}(\Hess f(\X)) \wedge 0, \text{ for } (p_1 +p_2-r)r + 1 \le i \le (p_1+p_2)r-(r^2-r)/2;
	\end{split}
	\end{equation}
	 \item the following upper bounds for $\widebar{\lambda}_i(\nabla^2 g_\reg({\L,\R}))$ hold: 
	 \begin{equation}\label{ineq: spectrum-bound-regularized-2}
	 	\begin{split}
	 		\widebar{\lambda}_i(\nabla^2 g_\reg({\L,\R})) &\leq 2 \sigma_1(\X)( (\lambda_{1}(\Hess f(\X)) \vee 0 ) + 4\mu ), \text{ for } 1\leq  i \leq (r^2+r)/2,\\
	 		\widebar{\lambda}_{i}(\nabla^2 g_\reg({\L,\R})) 
	 	&\leq  \left( 2\sigma_1(\X)  \lambda_{(i - (r^2+2)/2)}(\Hess f(\X)) \right)  \vee \left( \sigma_r(\X)\lambda_{(i - (r^2+2)/2)}(\Hess f(\X)) \right), \\
	 		&\quad \quad \text{ for } (r^2+r)/2 +1 \leq i \leq (p_1+p_2)r-(r^2-r)/2.
	 	\end{split}
	 \end{equation}
\end{itemize} 
\end{Theorem}

{\noindent \bf Proof of Theorem \ref{th: general case spec-connect 2}}. Since $(\L, \R)$ is a FOSP of \eqref{eq: general factor with reg}, the first result in \eqref{eq: reg-FOSP-property} is by Theorem 3 of \cite{zhu2018global}. The second result in \eqref{eq: reg-FOSP-property} is by $\L^\top \L= \R^\top \R$ and Proposition \ref{prop: gradient-exp}. In addition, for any $\A =  [\A_L^\top \quad \A_R^\top]^\top\in \bbR^{(p_1 + p_2) \times r}$, \eqref{eq: reg-Hessian-on-FOSP} follows from \eqref{eq: reg-FOSP-property} and Proposition \ref{prop: Hessian-exp}.

 The rest of the proof is divided into two steps. In Step 1, we prove the second part of Theorem \ref{th: general case spec-connect 2}, i.e., \eqref{eq: R-E-Hessian-general-reg} and \eqref{eq: R-E-Hessian-general-reg-2}; in Step 2, we prove the final part of the theorem, i.e., the spectrum bounds in \eqref{ineq: spectrum-bound-regularized-1} and \eqref{ineq: spectrum-bound-regularized-2}.

{\bf Step 1}. First, by the first equality in \eqref{eq: R-E-Hessian-general-1} and \eqref{eq: reg-Hessian-on-FOSP}, we obtain the first equality in  \eqref{eq: R-E-Hessian-general-reg}. Since $(\L, \R)$ is a FOSP of \eqref{eq: general factor with reg}, from \eqref{eq: reg-FOSP-property}, we see that $(\L,\R)$ is also a Euclidean FOSP of \eqref{eq: general factor formu}. Recalling the definition of $\cL_{\L,\R}^{-1}$ in \eqref{eq: bijection-general-case}, given $\xi \in T_{\X} \cM_r$, we see that $\cL_{\L,\R}^{-1}(\xi) = \A_{\L,\R}^\xi = [\A_L^{\xi\top} \quad \A_R^{\xi\top}]^\top$ satisfies $$\L^\top \A_L^\xi + \A_L^{\xi\top} \L- \R^\top \A^\xi_R - \A_R^{\xi\top} \R = \P_1^\top \widebar{\S} \P_2 + (\P_1^\top \widebar{\S} \P_2)^\top - \P_1^\top \widebar{\S} \P_2 - (\P_1^\top \widebar{\S} \P_2)^\top = \0.$$ So the second equality in \eqref{eq: R-E-Hessian-general-reg} follows from \eqref{eq: reg-Hessian-on-FOSP} and the second equality in \eqref{eq: R-E-Hessian-general-1}. 

 The second equality in \eqref{eq: R-E-Hessian-general-reg-2} directly follows from the first equality in \eqref{eq: R-E-Hessian-general-reg} and the definition of $\widetilde{\cL}_{\L,\R}$. Next, we prove the first equality in \eqref{eq: R-E-Hessian-general-reg-2}. For any $\A = [\A_L^\top \quad \A_R^\top]^\top \in \widetilde{\scA}_{\rnull}^{\,\,\L,\R}, \B = [\B_L^\top \quad \B_R^\top ]^\top \in \bbR^{(p_1 + p_2) \times r}$,
 \begin{equation*}
 	\begin{split}
 		\nabla^2 g_\reg(\L,\R)[\A,\B] &= \left( \nabla^2 g_\reg(\L,\R)[\A + \B,\A+\B] - \nabla^2 g_\reg(\L,\R)[\A-\B,\A-\B]  \right)/4\\
 		& \overset{(a) }= \big( \Hess f(\X)[\xi_{\L,\R}^{\A+ \B},\xi_{\L,\R}^{\A+ \B}]  + \mu \|\L^\top \B_L + \B_L ^\top \L- \R^\top \B_R - \B_R^\top \R\|_\F^2 \\
 		& \quad -\Hess f(\X)[\xi_{\L,\R}^{\A- \B},\xi_{\L,\R}^{\A- \B}]  - \mu \|\L^\top \B_L + \B_L ^\top \L- \R^\top \B_R - \B_R^\top \R\|_\F^2 \big)/4\\
 		& = \left(\Hess f(\X)[\xi_{\L,\R}^{\A} + \xi_{\L,\R}^{\B}, \xi_{\L,\R}^{\A} + \xi_{\L,\R}^{\B}] - \Hess f(\X)[\xi_{\L,\R}^{\A} - \xi_{\L,\R}^{\B}, \xi_{\L,\R}^{\A} -\xi_{\L,\R}^{\B}]  \right)/4\\
 		& = \Hess f(\X)[\xi_{\L,\R}^{\A}, \xi_{\L,\R}^{\B}] \overset{(b)}= 0,
 	\end{split}
 \end{equation*} where (a) is because of \eqref{eq: R-E-Hessian-general-reg} and $\L^\top \A_L + \A_L ^\top \L- \R^\top \A_R - \A_R^\top \R = \0$ for any $ \A \in \widetilde{\scA}_{\rnull}^{\,\,\L,\R}$ and (b) is because $\xi_{\L,\R}^\A = \0$ for any $ \A \in \widetilde{\scA}_{\rnull}^{\,\,\L,\R}$. This implies the first equality in \eqref{eq: R-E-Hessian-general-reg-2} and finishes the proof of this part.  

{\bf Step 2.} It is easy to check the number of eigenvalues of $\Hess f(\X)$ and $\nabla^2 g_\reg({\L,\R})$ are $(p_1 +p_2-r)r$ and $(p_1+p_2)r$, respectively. By the first equality in \eqref{eq: R-E-Hessian-general-reg-2}, we have $\widetilde{\scA}_{\rnull}^{\,\,\L,\R}$ is the eigenspace of $(r^2-r)/2$ zero eigenvalues and the rest of the $((p_1+p_2)r - (r^2-r)/2)$ possibly non-zero eigenvalues of $\nabla^2 g_\reg({\L,\R})$ span the eigenspace $\widetilde{\scA}_{\overline{\rnull}}^{\,\,\L,\R}$. Restricting $\Hess f(\X)$ and $\nabla^2 g_\reg(\L,\R)$ to $T_{\X} \cM_r$ and $\widetilde{\scA}_{\overline{\rnull}}^{\,\,\L,\R}$, respectively, next we prove the inequalities in \eqref{ineq: spectrum-bound-regularized-1} and \eqref{ineq: spectrum-bound-regularized-2} sequentially. 

{\bf Proof of the first inequality in \eqref{ineq: spectrum-bound-regularized-1}}. Define the linear map $ \cP:T_\X \cM_r \to \widetilde{\scA}_{\overline{\rnull}}^{\,\,\L,\R}$ as $\cP(\xi)=\A_{\L,\R}^\xi$. By the definition of $\cL_{\L,\R}^{-1}$ and second equality in \eqref{eq: R-E-Hessian-general-reg}, we have
\begin{equation*} \label{eq: proof-first-ineq-partial-sand}
    \Hess f(\X)[\xi, \xi]=  \nabla^2 g_\reg(\L,\R)[\cL_{\L,\R}^{-1}(\xi),\cL_{\L,\R}^{-1}(\xi)] = \nabla^2 g_\reg(\L,\R)[\cP(\xi), \cP(\xi)], \quad \forall \xi \in T_\X \cM_r,
\end{equation*}
i.e., 
\begin{equation}
\label{eq:hessfeqpgp}
\Hess f(\X) = \cP^* \nabla^2 g_\reg(\L,\R) \cP.
\end{equation}
Moreover, by the construction of $\cP$ and \eqref{ineq: bijection-spectrum-bound-general1}, we have for any $\xi \in T_\X \cM_r$,
\begin{equation} \label{ineq: reg-spectrum-connect-ineq}
	(2\sigma_1(\X))^{-1} \|\xi\|_\F^2 \overset{(a)}= (\sigma_1(\L) \vee \sigma_1(\R) )^{-2} \|\xi \|_\F^2/2 \leq \|\cP (\xi) \|_\F^2  \leq ( \sigma_r(\L) \wedge \sigma_r(\R))^{-2} \|\xi\|_\F^2 \overset{(b)}= \sigma_r(\X)^{-1} \|\xi\|_\F^2.
\end{equation}In (a) and (b), we use the fact that when \eqref{eq: reg-FOSP-property} holds, we have
\begin{equation} \label{eq: scale-R-L-balanced}
	\sigma_1(\L)= \sigma_1(\R) = \sigma_1^{1/2}(\X), \quad \sigma_r(\L)= \sigma_r(\R) = \sigma_r^{1/2}(\X).
\end{equation} Finally, by \eqref{eq:hessfeqpgp}, \eqref{ineq: reg-spectrum-connect-ineq} and Lemma \ref{lm: two-diff-symmetric-matrix-spectrum-connection1}(i) in the Appendix, we have obtained the first inequality in \eqref{ineq: spectrum-bound-regularized-1}. 

{\bf Proof of the second inequality in \eqref{ineq: spectrum-bound-regularized-1}}. By the second equality in \eqref{eq: R-E-Hessian-general-reg-2}, we have $\nabla^2g_\reg(\L,\R) \succcurlyeq \widetilde{\cL}_{\L,\R}^* \Hess f(\X) \widetilde{\cL}_{\L,\R} $. Then by \eqref{ineq: bijection-spectrum-bound-general1}, \eqref{eq: scale-R-L-balanced} and Lemma \ref{lm: two-diff-symmetric-matrix-spectrum-connection1}(iii) in the Appendix, we have for $(p_1 + p_2 -r)r +1 \leq i \leq (p_1+p_2)r-(r^2-r)/2$, $$\widebar{\lambda}_{i}(\nabla^2 g_\reg({\L,\R})) \geq \widebar{\lambda}_{(p_1+p_2)r-(r^2-r)/2}(\nabla^2 g_\reg({\L,\R})) \geq  2\sigma_1(\X) \lambda_{\min}(\Hess f(\X)) \wedge 0.$$

{\bf Proof of the first inequality in \eqref{ineq: spectrum-bound-regularized-2}}. By the second equality in \eqref{eq: R-E-Hessian-general-reg-2}, we have for any $\A \in \widetilde{\scA}_{\overline{\rnull}}^{\,\, \L,\R}$:
\begin{equation*}
\begin{split}
\nabla^2g_\reg(\L,\R)[\A,\A]& = \Hess f(\X)[\widetilde{\cL}_{\L,\R}(\A),\widetilde{\cL}_{\L,\R}(\A)] + \mu \|\L^\top \A_L + \A_L ^\top \L- \R^\top \A_R - \A_R^\top \R \|_\F^2 \\
	& \overset{\text{Lemma } \ref{lm: reg-term-Hessian-bound}, \eqref{eq: scale-R-L-balanced}  } \leq  \Hess f(\X)[\widetilde{\cL}_{\L,\R}(\A),\widetilde{\cL}_{\L,\R}(\A)] + 8\mu \sigma_1(\X) \|\A\|_\F^2.
\end{split}
\end{equation*}

So we have $(\nabla^2 g_\reg(\L,\R) - 8\mu \sigma_1(\X) \mathcal{I}) \preccurlyeq\widetilde{\cL}_{\L,\R}^* \Hess f(\X) \widetilde{\cL}_{\L,\R}$ where $\mathcal{I}$ denotes an identity operator. Then by \eqref{ineq: bijection-spectrum-bound-general1}, \eqref{eq: scale-R-L-balanced} and Lemma \ref{lm: two-diff-symmetric-matrix-spectrum-connection1}(iv) in the Appendix, we have for $1\leq i \leq (r^2 + r)/2$:
$$\widebar{\lambda}_i(\nabla g_\reg({\L,\R})) - 8\mu \sigma_1(\X) \leq \widebar{\lambda}_{1}(\nabla g_\reg({\L,\R})) - 8\mu \sigma_1(\X) \leq 2\sigma_1(\X) \lambda_1(\Hess f(\X)) \vee 0.$$

{\bf Proof of the second inequality in \eqref{ineq: spectrum-bound-regularized-2}}. 
The desired inequality can be obtained by \eqref{eq:hessfeqpgp}, 
\eqref{ineq: reg-spectrum-connect-ineq} and Lemma \ref{lm: two-diff-symmetric-matrix-spectrum-connection1}(ii) in the Appendix. This finishes the proof. \quad $\blacksquare$

\begin{Remark}[Comparison of Regularized and Unregularized Factorization Formulations] \label{rem: two-factorization-comparison}
	Compared to Theorem \ref{th: general case spec-connect 2}, the gap of the sandwich inequality in Theorem \ref{th: general case spec-connect 1} depends explicitly on the spectrum of $\L$ and $\R$ and can be arbitrarily large for ill-conditioned $(\L,\R)$ pairs. Such an issue makes the geometry analysis for the unregularized factorization \eqref{eq: general factor formu} hard \citep{zhang2020symmetry}. On the other hand, any Euclidean FOSP of the regularized formulation \eqref{eq: general factor with reg} is always balanced \citep{zhu2018global}, i.e., satisfying $\L^\top \L = \R^\top \R$, and the gap of the sandwich inequality in Theorem \ref{th: general case spec-connect 2} only depends on $\X = \L \R^\top$, not individual $\L$ or $\R$.
	
	In addition, comparing two factorization formulations \eqref{eq: general factor formu} and \eqref{eq: general factor with reg}, $\nabla^2 g_\reg(\L,\R)$ has $(r^2 + r)/2$ less zero eigenvalues than $\nabla^2 g(\L,\R)$ as the regularization reduces the ambiguity set from invertible transforms to rotational transforms. On the other hand, it is difficult to control these potentially non-zero $(r^2 + r)/2$ eigenvalues in $\nabla^2 g_\reg(\L,\R)$ due to the complex interaction between the regularization and the original objective function. So that is why in Theorem \ref{th: general case spec-connect 2}, we can only get a partial sandwich inequality between $\widebar{\lambda}_i(\nabla^2 g_\reg(\L,\R) )$ and $\lambda_i(\Hess f(\X))$ while in Theorems \ref{th: RHessian-EHessian PSD} and \ref{th: general case spec-connect 1} we have full sandwich inequalities. 
\end{Remark}

\begin{Remark}[Comparison of PSD and General Cases] \label{rem: comparison-PSD-general-case}
	There are a few similarities and key differences in the landscape connection under the PSD case and the general case. First, in both cases, we tackle the problem via finding a connection between Riemannian and Euclidean Hessians on some carefully constructed points. However, exact Riemannian and Euclidean Hessian connections between \eqref{eq: PSD-manifold-formulation} and \eqref{eq: PSD factorization} as well as \eqref{eq: general prob} and \eqref{eq: general factor formu} are available, while the Hessian connection between \eqref{eq: general prob} and \eqref{eq: general factor with reg} is weaker. Second, although sandwich inequalities between the spectrum of Riemannian and Euclidean Hessians can be established in both the PSD case (\eqref{eq: PSD-manifold-formulation} and \eqref{eq: PSD factorization}) and the general case (\eqref{eq: general prob} and \eqref{eq: general factor formu}), the gap of the sandwich inequality in the general case depends on the balancing of two factors $\L$, $\R$ as we mentioned in Remark \ref{rem: two-factorization-comparison} while there is no such an issue in the PSD case. Finally, compared to the PSD one, there are two factorization formulations in the general case (unregularized and regularized ones) and it is nontrivial to extend the results from the PSD case to the general case. In particular, the regularized factorization formulation can potentially have a distinct landscape geometry from the unregularized one and establishing the landscape connection between \eqref{eq: general prob} and \eqref{eq: general factor with reg} is much harder than  \eqref{eq: PSD-manifold-formulation} and \eqref{eq: PSD factorization} or \eqref{eq: general prob} and \eqref{eq: general factor formu}.
\end{Remark}

By Theorems \ref{th: general case spec-connect 1} and \ref{th: general case spec-connect 2}, we have the following Corollary \ref{coro: landscape connection general case} on the equivalence of FOSPs, SOSPs and strict saddles between the manifold and the factorization formulations in the general low-rank matrix optimization.
\begin{Corollary}({\bf Equivalence on FOSPs, SOSPs and Strict Saddles of Manifold and Factorization Formulations (General Case)}) \label{coro: landscape connection general case} (a) If $(\L,\R)$ is a rank $r$ Euclidean FOSP or SOSP or strict saddle of \eqref{eq: general factor formu} or \eqref{eq: general factor with reg}, then $\X = \L \R^\top$ is a Riemannian FOSP or SOSP or strict saddle of \eqref{eq: general prob}; (b) if $\X$ is a Riemannian FOSP or SOSP or strict saddle of \eqref{eq: general prob}, then any $(\L, \R)$ such that $\L \R^\top = \X$ is a Euclidean FOSP or SOSP or strict saddle of \eqref{eq: general factor formu} and any $(\L, \R)$ such that $\L \R^\top = \X, \L^\top \L = \R^\top \R$ is a Euclidean FOSP or SOSP or strict saddle of \eqref{eq: general factor with reg}.
\end{Corollary}

In Theorems \ref{th: general case spec-connect 1} and \ref{th: general case spec-connect 2}, we present the geometric landscape connection between the manifold and the two factorization formulations in the general low-rank matrix optimization. There is also a simple landscape connection between the two factorization formulations \eqref{eq: general factor formu} and \eqref{eq: general factor with reg}. This connection will be used to analyze the role of regularization in Section \ref{sec: role-regularization-low-rank-optimization}. 

\begin{Theorem}({\bf Geometric Landscape Connection of Unregularized Formulation \eqref{eq: general factor formu} and Regularized Formulation \eqref{eq: general factor with reg}}) \label{th: spectrum-connection-two-factorizations} Suppose $(\L,\R)$ and $(\L_\reg, \R_\reg)$ are rank $r$ Euclidean FOSPs of $g(\L,\R)$ and $g_\reg(\L,\R)$, respectively and $\L\R^\top = \L_\reg \R_\reg^\top$. Let $\bDelta = \L^{\dagger} \L_\reg $. Then $\bDelta$ is nonsingular and we can find a linear bijection $\cJ$ on $\bbR^{(p_1 + p_2) \times r}$:
\begin{equation*}
\begin{split}
	\cJ: &\A = \begin{bmatrix}
		\A_L\\
		\A_R
	\end{bmatrix} \in \bbR^{(p_1 + p_2) \times r} \longrightarrow \A'  = \begin{bmatrix}
		\A_L \bDelta^{-1} \\
		\A_R \bDelta^\top
	\end{bmatrix} \in \bbR^{(p_1 + p_2) \times  r},
\end{split}
\end{equation*} such that 
\begin{equation} \label{eq: two-Euclidean-Hessian-connection}
	\nabla^2 g_\reg (\L_\reg, \R_\reg) [\A,\A] -  \mu \|\L_\reg^\top \A_L + \A_L^\top \L_\reg - \R_\reg^\top \A_R - \A_R^\top \R_\reg \|_\F^2 = \nabla g^2 (\L,\R)[ \cJ(\A), \cJ(\A) ]
\end{equation} holds for any $\A = [\A_L^\top \quad \A_R^\top]^\top  \in \bbR^{(p_1 + p_2) \times r}$.

Moreover, we have the following spectrum bounds for $\cJ$:
\begin{equation} \label{ineq: bijection-spectrum-bound-two-Euclidean-Hessian}
	\theta_{\bDelta}^2 \|\A\|_\F^2 \leq \|\cJ (\A)\|_\F^2 \leq \Theta_{\bDelta}^2 \|\A\|_\F^2, \quad \forall \A \in \bbR^{(p_1 + p_2) \times r},
\end{equation} where $\Theta_{\bDelta} := \sigma_1(\bDelta) \vee (1/\sigma_{r}(\bDelta ))$ and $\theta_{\bDelta}:= 1/\Theta_{\bDelta} = (1/\sigma_1(\bDelta)) \wedge \sigma_{r}(\bDelta) $.

Finally, for $1 \leq i \leq (p_1 + p_2)r$, the following connections on individual eigenvalues between $ \nabla^2 g_\reg(\L_\reg, \R_\reg) $ and $\nabla^2 g(\L,\R)$ hold:
\begin{equation} \label{ineq: spectrum-bound-two-Euclidean-Hessian}
	\begin{split}
		\lambda_i ( \nabla^2 g(\L,\R)  ) &\leq \left(   \Theta^2_{\bDelta} \lambda_i( \nabla^2 g_\reg(\L_\reg,\R_\reg) )  \right)  \vee  \left( \theta^2_{\bDelta} \lambda_i( \nabla^2 g_\reg(\L_\reg,\R_\reg) )   \right),\\
		\lambda_i ( \nabla^2 g(\L,\R)  ) &\geq \left(   \Theta^2_{\bDelta}  \left( \lambda_i( \nabla^2 g_\reg(\L_\reg,\R_\reg) ) - 8 \mu \sigma_1(\L_\reg \R_\reg^\top) \right)  \right) \\
		& \quad \vee  \left( \theta^2_{\bDelta} \left( \lambda_i( \nabla^2 g_\reg(\L_\reg,\R_\reg) ) - 8 \mu \sigma_1(\L_\reg \R_\reg^\top) \right)  \right).
	\end{split}
\end{equation}
\end{Theorem}

\section{Applications} \label{sec: application}
In this section, we apply our main results to three specific problems from machine learning and signal processing. 

\subsection{Global Optimality for Phase Retrieval Under Manifold Formulation} \label{sec: phase-retrieval-global-optimality}
We first consider the following real-valued quadratic equation system
\begin{equation} \label{eq: phase retrieval model}
	\y_i = \langle\a_i, \x^*\rangle^2 \quad \text{for} \quad 1 \leq i \leq n,
\end{equation} 
where $\y \in \bbR^n$ and covariates $\{\a_i \}_{i=1}^n  \in \bbR^{p}$ are known whereas $\x^* \in \bbR^p$ is unknown. The goal is to recover $\x^*$ based on $\{\y_i, \a_i\}_{i=1}^n$. One important application is known as \emph{phase retrieval} arising from physical science due to the nature of optical sensors \citep{fienup1982phase}. A common formulation to solve \eqref{eq: phase retrieval model} is the following least squares formulation:
\begin{equation} \label{eq: least-square-phase-retrieval}
	\tilde{g}(\x) = \frac{1}{2n} \sum_{i=1}^n (\y_i - \langle\a_i, \x\rangle^2 )^2.
\end{equation}
  In the literature, both convex relaxation \citep{candes2013phaselift,waldspurger2015phase} and nonconvex approaches \citep{candes2015phase,chen2017solving,ma2019implicit,netrapalli2013phase,sanghavi2017local,wang2017solving,cai2018solving} have been proposed to solve \eqref{eq: least-square-phase-retrieval} with provable recovery guarantees. In terms of the geometric landscape analysis, \cite{sun2018geometric} showed that under the Gaussian design, i.e., $\a_i$s are drawn from i.i.d. Gaussian distribution, $\tilde{g}(\x)$ does not have any spurious local minima if $n \geq C p \log^3 p$ for some positive constant $C$. Later, the sample complexity requirement for the global optimality in phase retrieval under Gaussian design was improved to $n \geq C p$ for a slightly modified loss function \citep{li2019toward}:
\begin{equation} \label{eq: modi-formulation-phase-retrieval}
	g(\x) = \frac{1}{2n} \sum_{i=1}^n (\y_i - \langle\a_i, \x\rangle^2 )^2 h\left( \frac{ \langle\a_i, \x\rangle^2}{\|\x\|_2^2} \right) h\left( \frac{n \y_i}{\|\y\|_1} \right),
\end{equation} where for two predetermined universal parameters $1 < \beta < \gamma$, the twice continuously differential activation function $h(a)$ satisfies:
\begin{equation} \label{eq h-phase-retrieval}
	\left\{ \begin{array}{l c}
		h(a) =1 & \text{ if } 0 \leq a \leq \beta,\\
		h(a) \in [0,1] & \text{ if } a \in (\beta,\gamma),\\
		h(a) =0 & \text{ if } a \geq \gamma
	\end{array} \right.
\end{equation} and $|h'(a)|, |h''(a)|$ exist and are bounded. 

Comparing objectives in \eqref{eq: least-square-phase-retrieval} and \eqref{eq: modi-formulation-phase-retrieval}, $g(\x)$ incorporates a smooth activation function $h$ to handle the heavy-tailedness of the fourth moment of Gaussian random variables in $\tilde{g}(\x)$. On the other hand, for both \eqref{eq: least-square-phase-retrieval} and \eqref{eq: modi-formulation-phase-retrieval}, the geometric landscape analyses performed in \cite{sun2018geometric,li2019toward} are carried out in terms of $\x$ in the vector space. However, it is known that by lifting $\x$ to $\X = \x \x^\top$, both \eqref{eq: least-square-phase-retrieval} and \eqref{eq: modi-formulation-phase-retrieval} can be recast as a rank-1 PSD matrix recovery problem, e.g., $\min_{\x \in \bbR^{p}} g(\x)$ has the following equivalent PSD manifold formulation:
\begin{equation} \label{eq: phase-retrieval-manifold}
\min_{\X \in \bbS^{p \times p}: \rank{(\X)} = 1, \X \succcurlyeq 0 } f(\X) := \frac{1}{2n} \sum_{i=1}^n (\y_i - \langle\A_i, \X\rangle )^2 h\left( \frac{ \langle\A_i, \X\rangle}{\|\X\|_\F} \right) h\left( \frac{n \y_i}{\|\y\|_1} \right),
\end{equation} where $\A_i = \a_i \a_i^\top$ for $i = 1,\ldots, n$. Since $h$ is twice continuously differentiable, the objective function in \eqref{eq: phase-retrieval-manifold} is also twice continuously differentiable over rank-1 PSD matrices.\cite{cai2018solving,li2019toward} asked whether it is possible to investigate the geometric landscape of the phase retrieval problem directly on the rank-$1$ PSD matrix manifold. By our Theorem \ref{th: RHessian-EHessian PSD} and Corollary \ref{coro: landscape connection PSD}, we provide an affirmative answer to their question and give the first global optimality result for phase retrieval under the manifold formulation with a rate-optimal sample complexity.
\begin{Theorem}[Global Optimality for Phase Retrieval under Manifold Formulation]\label{th: manifold-landscape-phase-retrieval} In \eqref{eq: phase-retrieval-manifold}, suppose $\a_i \overset{i.i.d.} \sim N(0,\I_p)$, $\gamma > \beta > 1$ are sufficiently large in the smooth activation function $h$, and $n \geq C p$ for large enough positive constant $C$. Then with probability at least $1 - \exp(-C'n)$ for $C' > 0$, $\X^* = \x^* \x^{*\top}$ is the unique Riemannian SOSP of \eqref{eq: phase-retrieval-manifold} and any other Riemannian FOSP $\X$ is a strict saddle with $\lambda_{\min} (\Hess f(\X)) \leq  - \frac{3 \sigma_1(\X^*) }{4\sigma_1(\X)}$.
\end{Theorem}

\begin{Remark}[Transferring the Strict Saddle Property] \label{rem: transferring-strict-saddle} The key reason we can establish the global optimality and strict saddle results for phase retrieval under the manifold formulation is attributed to the spectrum connection of the Riemannian and the Euclidean Hessians given in Theorem \ref{th: RHessian-EHessian PSD}. As a result of that, we can transfer the strict saddle property \citep{ge2015escaping,lee2019first}, which states that the function has a strict negative curvature at all stationary points but local minima, from the factorization formulation of phase retrieval to the manifold one. This is fundamentally different from the results in \cite{ha2020equivalence} where only the connection between Euclidean SOSPs and fixed points of PGD was established without giving the estimation on the curvature of the Hessian. With this strict saddle property, various gradient descent and trust region methods are guaranteed to escape all strict saddles and converge to a SOSP \citep{o2020line,ge2015escaping,lee2019first,jin2017escape,paternain2019newton,sun2018geometric,sun2019escaping,criscitiello2019efficiently,boumal2019global,han2020escape}. {Finally, we note here an explicit upper bound on the negative eigenvalue of the strict saddle in Theorem \ref{th: manifold-landscape-phase-retrieval} can be helpful in determining the convergence rate of perturbed GD to the global minima \citep{jin2017escape,sun2019escaping,criscitiello2019efficiently}.}
\end{Remark}

\subsection{Global Optimality of General Well-Conditioned Low-rank Matrix Optimization Under Manifold Formulation} \label{sec: general-low-rank-opt-global-optimality}

In the existing literature on low-rank matrix optimization, most of the geometric landscape analyses focused on the factorization formulation. They showed that doing factorization for a rank constrained objective will not introduce spurious local minima when the objective $f$ satisfies the restricted strong convexity and smoothness property (see the upcoming Definition \ref{def: RSC-RSM}) \citep{bhojanapalli2016global,ge2017no,zhang2019sharp,zhu2018global,zhu2017global,park2017non}. On the other hand, the geometric analysis performed directly under the rank constrained manifold formulation is scarce. \cite{uschmajew2018critical} showed the benign landscape of \eqref{eq: general prob} under the embedded manifold geometry and exact-parameterization setting, i.e., $r = \rank(\X^*)$, where $\X^*$ is a low-rank parameter matrix of interest, when $f$ is quadratic and satisfies certain restricted spectral bounds property.

\begin{Definition} \label{def: RSC-RSM}
We say $f: \bbR^{p_1 \times p_2}\to \bbR$ satisfies the $(2r,4r)$-restricted strong convexity and smoothness property if for any $\X,\G \in \bbR^{p_1 \times p_2}$ with $\rank(\X) \leq 2r$ and $\rank(\G) \leq 4r$, the Euclidean Hessian of $f$ satisfies
\begin{equation} \label{eq: RSC-RSM}
	\alpha_1 \|\G\|_\F^2 \leq \nabla^2 f(\X)[\G,\G] \leq \alpha_2 \|\G\|_\F^2
\end{equation} 
for some $\alpha_2 \geq \alpha_1 > 0$.
\end{Definition}

By Theorem \ref{th: general case spec-connect 2} and Corollary \ref{coro: landscape connection general case}, we can leverage the existing results in \cite{zhu2018global} to provide the first global optimality result of \eqref{eq: general prob} under the manifold formulation for a generic objective $f$ satisfying the restricted strong convexity and smoothness property. Moreover, our results cover both over-parameterization ($r > \rank(\X^*)$) and exact-parameterization ($r = \rank(\X^*)$) settings compared with \cite{uschmajew2018critical}.

\begin{Theorem}({\bf Global Optimality of Well-Conditioned Low-rank Matrix Optimization under Manifold Formulation}) \label{th: global-optimality-generic-f-manifold-formulation}
	Consider the optimization problem  \eqref{eq: general prob}. Suppose there exists a rank $r^*$ ($r^* \leq r$) matrix $\X^*$ s.t. $\nabla f(\X^*) = \0$ and $f$ satisfies the $(2r,4r)$-restricted strong convexity and smoothness property \eqref{eq: RSC-RSM} with positive constants $\alpha_1$ and $\alpha_2$ satisfying $\alpha_2/\alpha_1 \leq 1.5$. Then,
	\begin{itemize}
		\item if $r = r^*$, $\X^*$ is the unique Riemannian SOSP of \eqref{eq: general prob} and any other Riemannian FOSP $\X$ is a strict saddle with $\lambda_{\min}(\Hess f(\X)) \leq -0.04 \alpha_1 \sigma_r(\X^*)/\sigma_1(\X) $;
		\item if $r > r^*$, there is no Riemannian SOSP of \eqref{eq: general prob} and any Riemannian FOSP $\X$ is a strict saddle with $\lambda_{\min}(\Hess f(\X)) \leq -0.05 \alpha_1 ( \sigma_r(\X) \wedge \sigma_{r^*}(\X^*) ) /\sigma_1(\X) $.
	\end{itemize}
\end{Theorem}

\begin{Remark}
	As guaranteed by Proposition 1 of \cite{zhu2018global}, the $(2r,4r)$-restricted strong convexity and smoothness property of $f$ ensures $\X^*$ in Theorem \ref{th: global-optimality-generic-f-manifold-formulation} is the unique global minimizer of $\min_{\X: \rank(\X) \leq r} f(\X)$. So Theorem \ref{th: global-optimality-generic-f-manifold-formulation} shows that if the input rank $r$ is equal to the true rank $r^*$, i.e., under exact-parameterization, then \eqref{eq: general prob} has no spurious local minimizer other than the global minima $\X^*$ and any other Riemannian FOSP is a strict saddle. These two facts together ensure the recovery of $\X^*$ by many iterative algorithms \citep{lee2019first,sun2018geometric,sun2019escaping,criscitiello2019efficiently}.
	
	On the other hand, when the input rank $r$ is greater than the true rank $r^*$, i.e., under over-parameterization, Theorem \ref{th: global-optimality-generic-f-manifold-formulation} shows that there is no Riemannian SOSP for \eqref{eq: general prob} and all Riemannian FOSPs are strictly saddles. In addition, the upper bound on the negative curvature of the strict saddle implies that when running algorithms with guaranteed strict-saddle escaping property, the least singular value, i.e., $\sigma_r(\X)$, of the iterates will converge to zero. This suggests that the iterates tend to enter a lower rank matrix manifold and we can adopt some rank-adaptive Riemannian optimization methods to accommodate this \citep{zhou2016riemannian,gao2021riemannian}. {We note this observation is only possible due to an explicit upper bound on the negative eigenvalue at strict saddles powered by the sandwich inequalities Theorem \ref{th: general case spec-connect 2}. }
\end{Remark}

\subsection{Role of Regularization in Nonconvex Factorization for Low-rank Matrix Optimization} \label{sec: role-regularization-low-rank-optimization}
As we have discussed in the introduction, for nonconvex factorization of the general low-rank matrix optimization, the regularized formulation \eqref{eq: general factor with reg} is often considered. The regularization is introduced to balance the scale of two factors $\L,\R$ and it facilitates both algorithmic and geometric analyses in the nonconvex factorization formulation \citep{tu2016low,zheng2015convergent,ma2019implicit,wang2017unified,park2018finding,zhu2018global,zhu2017global,park2017non,ge2017no}. On the other hand, it has been first observed empirically \citep{zhu2018global}, and then recently proved that the regularization is not necessary for iterative algorithms to converge in a number of smooth and non-smooth formulated matrix inverse problems \citep{du2018algorithmic,charisopoulos2021low,tong2020accelerating,tong2020low,ma2021beyond,ye2021global}. Moreover, \cite{li2020global} showed from a geometric point of view that without regularization, the landscape of the factorization formulation \eqref{eq: general factor formu} is still benign when $f$ satisfies the restricted strong convexity and smoothness property \eqref{eq: RSC-RSM}.

In this paper, we provide more geometric landscape connections between two factorization formulations \eqref{eq: general factor formu} and \eqref{eq: general factor with reg} under a general $f$. Specifically, by connecting them with the manifold formulation, we show in Corollary \ref{coro: landscape connection general case} that the sets of $\L\R^\top$s formed by rank $r$ Euclidean FOSPs and SOSPs of two factorization formulations are exactly the same. By Theorem \ref{th: spectrum-connection-two-factorizations}, we also have a connection on the spectrum of Hessians at Euclidean FOSPs under two factorization formulations. If we further assume $f$ is well-conditioned as in Section \ref{sec: general-low-rank-opt-global-optimality}, then we can have the following global optimality result under the unregularized formulation \eqref{eq: general factor formu}.
\begin{Theorem}\label{th: global-optimality-factorization-noregularization}
	Consider the optimization problem  \eqref{eq: general factor formu}. Suppose there exists a rank $r$ matrix $\X^*$ s.t. $\nabla f(\X^*) = \0$ and $f$ satisfies the $(2r,4r)$-restricted strong convexity and smoothness property \eqref{eq: RSC-RSM} with positive constants $\alpha_1$ and $\alpha_2$ satisfying $\alpha_2/\alpha_1 \leq 1.5$. Then for any rank $r$ Euclidean FOSP $(\L,\R)$ of $g(\L,\R)$, it is either a Euclidean SOSP and satisfies $\L\R^\top = \X^*$, or a strict saddle with $\lambda_{\min}( \nabla^2 g(\L,\R) ) \leq -0.08 \left( (\sigma^2_r(\L)/\sigma_1(\L\R^\top)) \wedge (\sigma_1(\L\R^\top)/\sigma^2_r(\L)) \right) \cdot \alpha_1\sigma_r(\X^*)$. 
\end{Theorem}

Part of the results in Theorem \ref{th: global-optimality-factorization-noregularization} have appeared in the recent work \cite{li2020global}, but here we provide a precise upper bound on the negative curvature of the strict saddle that is absent in \cite{li2020global}. {Again, the precise upper bound on the negative curvature of the strict saddle is helpful in determining the convergence rate of perturbed GD to the global minima as we mentioned in Remark \ref{rem: transferring-strict-saddle}}.

\section{Conclusion and Discussions} \label{sec: conclusion}
In this paper, we consider the geometric landscape connection of the manifold and the factorization formulations in low-rank matrix optimization. We establish sandwich inequalities on the corresponding eigenvalues of the Riemannian and Euclidean Hessians and show an equivalence on the sets of FOSPs, SOSPs and strict saddles between two formulations. These results provide partial reasons for the similar empirical performance of manifold and factorization approaches in low-rank matrix optimization. Finally, we apply our main results to three applications in machine learning and signal processing.

There are many interesting extensions to the results in this paper to be explored in the future. First, as we have mentioned in Remark \ref{rem: necessity-of-FOSP-in-connection-Hessian}, our results on the connection of Riemannian and Euclidean Hessians are established at FOSPs. It is interesting to explore whether it is possible to connect the geometry of the manifold and the factorization formulations of low-rank matrix optimization at non-stationary points. By achieving this we can (1) connect approximate SOSPs\footnote{An approximate SOSP means the gradient norm at the point is small and the least eigenvalue of the Hessian at the point is lower bounded by a small negative constant \citep{jin2017escape}. } between two formulations, which is useful in practice as standard optimization methods such as stochastic or perturbed gradient descent can only find approximate SOSPs \citep{ge2015escaping,jin2017escape,sun2019escaping,criscitiello2019efficiently}; (2) transfer the global geometry properties (the landscape property of the objective in the whole space rather than at stationary points) between two formulations \citep{zhu2017global,li2019symmetry}. Second, in this work, we consider the natural embedded geometry of low-rank matrices in the manifold formulation. Another choice for handling low-rank matrices is the quotient manifold \citep{mishra2014fixed}. The follow-up work \cite{luo2023geometric} investigates the landscape connection of an optimization problem under the embedded and quotient geometries. Third, it is interesting to explore how will the landscape connect under two formulations when the objective function is nonsmooth. The connection of FOSPs might still be possible based on the notion of Clarke subdifferential \citep{clarke1990optimization}, but some regularity condition on $f$ might be needed. Finally, the manifold approach is a general way to deal with geometric constraints in optimization problems and here we show a strong geometric connection of it to the factorization approach in dealing with the rank constraint in matrix optimization. From an algorithmic perspective, connections of manifold methods with the sequential quadratic programming (SQP) method for solving equality constrained optimization problems and common nonlinear programming methods for handling orthogonal constraints were revealed in \cite{edelman1998geometry,mishra2016riemannian} and \cite{edelman1998geometry}, respectively. It is interesting future work to find more instances under which the manifold approach is geometrically or algorithmically connected with other well-known approaches in general nonlinear optimization.

\section*{Acknowledgments}
The authors thank the editors and reviewers for their great suggestions and comments. 

\bibliographystyle{apalike}
\bibliography{reference.bib}

\appendix
\section{Additional Proofs} \label{sec: additional-proofs}

\subsection{Additional Proofs in Section \ref{sec: Riemannian-opt-background}  }

{\noindent \bf Proof of Proposition \ref{prop: gradient-exp}.} The expressions for Euclidean gradients are obtained via direct computation. For the Riemannian gradient, since $\cM_{r+}$ and $\cM_r$ are embedded submanifolds of $\bbR^{p \times p}$ and $\bbR^{p_1 \times p_2}$, respectively and the Euclidean metric is considered, from \cite[(3.37)]{absil2009optimization}, we know the Riemannian gradients are the projections of the Euclidean gradients onto the corresponding tangent spaces. The results follow by observing the projection operator onto  $T_\X \cM_{r+}$ and $T_\X \cM_r$ given in \eqref{eq: tangent space projection}. 
\quad $\blacksquare$

\vskip.3cm
{\noindent \bf Proof of Proposition \ref{prop: Hessian-exp}.} First the expressions for $ \nabla^2 g(\L,\R)[\A,\A]$ and $\nabla^2 g_\reg (\L,\R)[\A,\A]$ are given in \cite[Eq. (2.8)]{ha2020equivalence} and \cite[Section IV-A and Remark 8]{zhu2018global}, respectively. The expressions for $ \nabla^2 g(\Y)[\A',\A'] $ can be obtained by letting $\A = [\A^{'\top} \quad \A^{'\top} ]^\top$ and $\L= \R = \Y$ in $\nabla^2 g(\L,\R)[\A,\A]$.
\vskip.3cm

Next, we derive the Riemannian Hessian of $f$. The Riemannian Hessian of an objective function $f$ is usually defined in terms of the Riemannian connection as in \eqref{def: Riemannain-Hessian}. But in the case of embedded submanifolds, it can also be defined by means of the so-called second-order retractions.

Given a general smooth manifold $\cM$, a {\it retraction} $R$ is a smooth map from $T \mathcal{M}$ to $\mathcal{M}$ satisfying i) $R(\X, 0) = \X$ and ii) $\frac{d}{d t} R(\X, t \eta) \vert_{t = 0} = \eta$ for all $\X \in \mathcal{M}$ and $\eta \in T_\X \mathcal{M}$, where $T \mathcal{M} = \{(\X, T_\X \mathcal{M}) : \X \in \mathcal{M}\}$, is the tangent bundle of $\mathcal{M}$ \cite[Chapter 4]{absil2009optimization}. We also let $R_\X$ to be the restriction of $R$ to $T_\X \cM$ and it satisfies $R_\X: T_\X \cM \to \cM, \xi \longmapsto R(\X,\xi)$. Retraction is in general a first-order approximation of the exponential map \cite[Chapter 4]{absil2009optimization}. A {\it second-order retraction} is the retraction defined as a second-order approximation of the exponential map \citep{absil2012projection}. As far as convergence of Riemannian optimization methods goes, first-order retraction is sufficient \cite[Chapter 3]{absil2009optimization}, but second-order retraction enjoys the following nice property: the Riemannian Hessian of an objective function $f$ coincides with the Euclidean Hessian of the lifted objective $\hat{f}_\X := f \circ R_\X$. \begin{Lemma}[Proposition 5.5.5 of \cite{absil2009optimization}] \label{lm: computation-Riemannian-hessian} Let $R_\X$ be a second-order retraction on $\cM$. Then $\Hess f(\X) = \nabla^2 ( f \circ R_\X )(0)$ for all $ \X \in \cM.$
\end{Lemma}

We present the second-order retractions under both PSD and general low-rank matrix settings in the following Proposition \ref{prop: second-order-retraction}.
\begin{Proposition}[Second-order Retractions in PSD and General Low-rank Matrix Manifolds] \label{prop: second-order-retraction}
\quad
\begin{itemize}[leftmargin=*]
	\item PSD case: Suppose $\X \in \cM_{r+}$ has eigendecomposition $\U \bSigma \U^\top$. Then the mapping $R_\X^{(2)}: T_\X \cM_{r+} \to \cM_{r+}$ given by 
\begin{equation*}
	R_\X^{(2)}: \xi = [\U \quad \U_\perp] \begin{bmatrix}
			\S & \D^\top\\
			\D & \0
		\end{bmatrix} [\U \quad \U_\perp]^\top \to \W \X^\dagger \W^\top
\end{equation*} is a second-order retraction on $\cM_{r+}$, where $\W = \X + \frac{1}{2} \xi^s + \xi^p - \frac{1}{8} \xi^s \X^\dagger \xi^s - \frac{1}{2} \xi^p \X^\dagger \xi^s $, $\xi^s = P_\U \xi P_\U$ and $\xi^p = P_{\U_\perp} \xi P_{\U} + P_{\U} \xi P_{\U_\perp}$. Furthermore, we have
\begin{equation*}
	R_\X^{(2)}(\xi) = \X + \xi + \U_\perp \D \bSigma^{-1} \D^\top \U_\perp^\top + O(\|\xi\|_\F^3), \, \text{ as } \|\xi\|_\F \to 0.
\end{equation*}
\item  General case: Suppose $\X \in \cM_{r}$ has SVD $\U \bSigma \V^\top$. Then the mapping $R_\X^{(2)}: T_\X \cM_{r} \to \cM_{r}$ given by 
\begin{equation*}
	R_\X^{(2)}: \xi = [\U \quad \U_\perp] \begin{bmatrix}
			\S & \D_2^\top\\
			\D_1 & \0
		\end{bmatrix} [\V \quad \V_\perp]^\top \to \W \X^\dagger \W
\end{equation*} is a second-order retraction on $\cM_{r}$, where $\W = \X + \frac{1}{2} \xi^s + \xi^p - \frac{1}{8} \xi^s \X^\dagger \xi^s - \frac{1}{2} \xi^p \X^\dagger \xi^s - \frac{1}{2} \xi^s\X^\dagger \xi^p $, $\xi^s = P_\U \xi P_\V$ and $\xi^p = P_{\U_\perp} \xi P_{\V} + P_{\U} \xi P_{\V_\perp}$. Furthermore, we have
\begin{equation*}
	R_\X^{(2)}(\xi) = \X + \xi + \U_\perp \D_1 \bSigma^{-1} \D_2^\top \V_\perp^\top + O(\|\xi\|_\F^3), \, \text{ as } \|\xi\|_\F \to 0.
\end{equation*}
\end{itemize}	
\end{Proposition}
{\noindent \bf Proof of Proposition \ref{prop: second-order-retraction}}. The results for the PSD case can be found in \cite[Proposition 5.10]{vandereycken2010riemannian} and the results under the general case can be found in \cite[Proposition A.1]{vandereycken2013low} and \cite[Theorem 3]{shalit2012online}.\quad $\blacksquare$

By Lemma \ref{lm: computation-Riemannian-hessian} and the property of second-order retraction, the sum of the first three dominating terms in the Taylor expansion of $f \circ R_\X^{(2)} (\xi)$ w.r.t. $\xi$ are $f(\X) + \langle \grad f(\X),\xi \rangle + \frac{1}{2} \Hess f(\X)[\xi,\xi]$. By matching the corresponding terms and the expressions of $R_\X^{(2)}$ in Proposition \ref{prop: second-order-retraction}, we can get the quadratic expression for $\Hess f(\X)[\xi,\xi]$. 

Next, we discuss how to obtain $\Hess f(\X)[\xi,\xi]$ in PSD and general low-rank matrix manifolds, respectively.

{\bf PSD case}: Given small enough $ \xi = [\U \quad \U_\perp] \begin{bmatrix}
			\S & \D^\top\\
			\D & \0
		\end{bmatrix} [\U \quad \U_\perp]^\top$, define $\U_p = \U_\perp \D$. By Proposition \ref{prop: second-order-retraction} and Taylor expansion, we have
\begin{equation} \label{eq: Taylor-expan-PSD}
	\begin{split}
		f \circ R_\X^{(2)} (\xi) &= f(\X + \xi + \U_p \bSigma^{-1} \U_p^\top + O(\|\xi\|_\F^3) )\\
		& = f(\X + \xi + \U_p \bSigma^{-1} \U_p^\top) +  O(\|\xi\|_\F^3)\\
		& = f(\X + \xi ) + \langle \nabla f(\X + \xi),  \U_p \bSigma^{-1} \U_p^\top \rangle +  O(\|\xi\|_\F^3)\\
		& = f(\X) + \langle \nabla f(\X),\xi \rangle + \frac{1}{2} \nabla^2 f(\X)[\xi,\xi] + \langle \nabla f(\X), \U_p \bSigma^{-1} \U_p^\top  \rangle +  O(\|\xi\|_\F^3).
	\end{split}
\end{equation}
Since $\xi^p \X^\dagger \xi^p =  \U_p \bSigma^{-1} \U_p^\top$, where $\xi^p = P_{\U_\perp} \xi P_{\U} + P_{\U} \xi P_{\U_\perp}$, the second order term in \eqref{eq: Taylor-expan-PSD} is $ \frac{1}{2} \nabla^2 f(\X)[\xi,\xi] + \langle \nabla f(\X), \U_p \bSigma^{-1} \U_p^\top  \rangle$ and it equals to $ \frac{1}{2} \Hess f(\X)[\xi,\xi]$.

{\bf General case}:  Given small enough $ \xi = [\U \quad \U_\perp] \begin{bmatrix}
			\S & \D_2^\top\\
			\D_1 & \0
		\end{bmatrix} [\V \quad \V_\perp]^\top$, define $\U_p = \U_\perp \D_1$ and $\V_p = \V_\perp \D_2$. By Proposition \ref{prop: second-order-retraction} and Taylor expansion, we have
\begin{equation} \label{eq: Taylor-expan-general}
	\begin{split}
		f \circ R_\X^{(2)} (\xi) &= f(\X + \xi + \U_p \bSigma^{-1} \V_p^\top + O(\|\xi\|_\F^3) )\\
		& = f(\X + \xi + \U_p \bSigma^{-1} \V_p^\top) +  O(\|\xi\|_\F^3)\\
		& = f(\X + \xi ) + \langle \nabla f(\X + \xi),  \U_p \bSigma^{-1} \V_p^\top \rangle +  O(\|\xi\|_\F^3)\\
		& = f(\X) + \langle \nabla f(\X),\xi \rangle + \frac{1}{2} \nabla^2 f(\X)[\xi,\xi] + \langle \nabla f(\X), \U_p \bSigma^{-1} \V_p^\top  \rangle +  O(\|\xi\|_\F^3).
	\end{split}
\end{equation}	Since $\xi^p \X^\dagger \xi^p =  \U_p \bSigma^{-1} \V_p^\top$, where $P_{\U_\perp} \xi P_{\V} + P_{\U} \xi P_{\V_\perp}$, the second order term in \eqref{eq: Taylor-expan-general} is $ \frac{1}{2} \nabla^2 f(\X)[\xi,\xi] + \langle \nabla f(\X), \U_p \bSigma^{-1} \V_p^\top  \rangle$ and it equals to $\frac{1}{2} \Hess f(\X)[\xi,\xi]$. This finishes the proof of this proposition. \quad $\blacksquare$

We note the proof technique for deriving the Riemannian Hessian is analogous to the proof of \cite[Proposition 2.3]{vandereycken2013low}. Here we extend it to the setting for a general twice differentiable function $f$.

\subsection{Additional Proofs in Section \ref{sec: connection-PSD}}

{\noindent \bf Proof of Lemma \ref{lm: tangent-vector-equiva-PSD}.} Suppose $\X$ has the eigendecomposition $\U\bSigma\U^\top$ and $\P = \U^\top \Y$. Given $\xi = [\U \quad \U_\perp] \begin{bmatrix}
			\S & \D^\top\\
			\D & \0
		\end{bmatrix} [\U \quad \U_\perp]^\top$. For any $\A \in \scA_\Y^\xi$, it is easy to check $\Y\A^\top +\A\Y^\top = \xi$, so $\scA_\Y^\xi \subseteq \{\A: \Y\A^\top + \A\Y^\top = \xi\}$. For any $\A$ such that $\Y\A^\top +\A\Y^\top = \xi$, we have
\begin{equation*} 
    \begin{bmatrix}
			\S & \D^\top\\
			\D & \0
		\end{bmatrix} = \begin{bmatrix}
		            \U^\top \\
		            \U_\perp^\top
		\end{bmatrix} \xi [\U \quad \U_\perp]  =  \begin{bmatrix}
		            \U^\top \\
		            \U_\perp^\top
		\end{bmatrix} (\Y\A^\top + \A\Y^\top)  [\U \quad \U_\perp] =\begin{bmatrix}
			\P\A^\top\U + \U^\top\A\P^\top & \P \A^\top \U_\perp\\
			\U_\perp^\top \A \P^\top & \0
		\end{bmatrix} 
\end{equation*} by observing $\Y = \U \P$. 
This implies $\U_\perp\U_\perp^\top\A = \U_\perp\D\P^{-\top}$ and $\P\A^\top\U + \U^\top\A\P^\top = \S$. By denoting $\S_1 = \U^\top\A\P^\top$, we have $\S_1 + \S_1^\top=\S$ and $\U\U^\top \A = \U\S_1\P^{-\top}$. Finally, 
$\A = \U\U^\top\A + \U_\perp\U_\perp^\top\A = (\U\S_1 + \U_\perp\D)\P^{-\top}\in \scA_\Y^\xi.$
This proves $\scA_\Y^\xi \supseteq \{\A: \Y\A^\top + \A\Y^\top = \xi\}$ and finishes the proof. \quad $\blacksquare$
\vskip.3cm
{\noindent \bf Proof of Lemma \ref{lm: decomposition-Rp-r-PSD}.} First, it is easy to check the dimensions of $\scA_{\rnull}^{\Y}$ and $\scA_{\overline{\rnull}}^{\Y}$ are $ (r^2-r)/2 $ and $pr-(r^2-r)/2$, respectively. Since $ (r^2-r)/2  + pr-(r^2-r)/2 = pr$, to prove $\bbR^{p \times r} = \scA_{\rnull}^{\Y} \oplus \scA_{\overline{\rnull}}^{\Y}$, we only need to show $\scA_{\rnull}^{\Y}$ is orthogonal to $\scA_{\overline{\rnull}}^{\Y}$. Suppose $\A = \U\S\P^{-\top} \in \scA_{\rnull}^{\Y}$ and $\A' =(\U\S' + \U_\perp \D' )\P^{-\top} \in \scA_{\overline{\rnull}}^{\Y}$. Then
	\begin{equation*}
		\begin{split}
			\langle \A, \A' \rangle = \langle \S \P^{-\top}, \S' \P^{-\top} \rangle = \langle \S , \S' \P^{-\top} \P^{-1} \rangle \overset{(a)}= \langle \S , \S' \bSigma^{-1} \rangle \overset{(b)}= - \langle \S^\top , (\S' \bSigma^{-1})^\top \rangle = -\langle \A, \A' \rangle,
		\end{split}
	\end{equation*} where (a) is because $\P \P^\top = \bSigma$, (b) is because $\S +\S^\top = \0$, and $\S' \bSigma^{-1}$ is symmetric by the construction of  $\scA_{\rnull}^{\Y}$ and $\scA_{\overline{\rnull}}^{\Y}$, respectively. So we have $\langle \A, \A' \rangle = 0$ and this finishes the proof of this lemma.
\quad $\blacksquare$

\vskip.3cm
{\noindent \bf Proof of Corollary \ref{coro: landscape connection PSD}.} First, by the connection of Riemannian and Euclidean gradients in \eqref{eq: gradient-connect-PSD}, the connection of FOSPs under two formulations clearly holds. 

Suppose $\Y$ is a rank $r$ Euclidean SOSP of \eqref{eq: PSD factorization} and let $\X = \Y \Y^\top$. Given any  $\xi \in T_{\X}\cM_{r+}$, we have
\begin{equation*}
	\Hess f(\X)[\xi, \xi] \overset{ \eqref{eq: R-E-Hessian-PSD}}= \nabla^2 g(\Y)[\cL_\Y^{-1}(\xi), \cL_\Y^{-1}(\xi)] \geq 0,
\end{equation*} where the inequality is by the SOSP assumption on $\Y$. Combining the fact $\X$ is a Riemannian FOSP of \eqref{eq: PSD-manifold-formulation}, this shows $\X = \Y \Y^\top$ is a Riemannian SOSP of \eqref{eq: PSD-manifold-formulation}.

Next, let us show the other direction: suppose $\X$ is a Riemannian SOSP of \eqref{eq: PSD-manifold-formulation}, then for any $\Y$ such that $\Y \Y^\top = \X$, it is a Euclidean SOSP of \eqref{eq: PSD factorization}. To see this, first $\Y$ is of rank $r$ and we have shown $\Y$ is a Euclidean FOSP of \eqref{eq: PSD factorization}. Then by \eqref{eq: R-E-Hessian-PSD}, we have for any $\A \in \bbR^{p \times r}$:
\begin{equation*}
\begin{split}
	\nabla^2 g(\Y)[\A,\A] = \Hess f(\X)[\xi_\Y^\A,\xi_\Y^\A] \geq 0. 
\end{split}
\end{equation*} 

Suppose $\Y$ is a rank $r$ Euclidean strict saddle of \eqref{eq: PSD factorization} and let $\X = \Y \Y^\top $. It implies that there exists $\A \in \scA_{\rnull}^\Y$ such that $\nabla^2 g(\Y)[\A, \A] < 0$. Then by \eqref{eq: R-E-Hessian-PSD} $\nabla^2 g(\Y)[\A,\A] = \Hess f(\X)[\cL(\A),\cL(\A)] < 0$, and this implies that $\Hess f(\X)$ also has at least one eigenvalue. Thus, $\X$ is a Riemannian strict saddle. The proof for the other direction is similar and for simplicity, we omit it here. \quad $\blacksquare$

\subsection{Additional Proofs in Section \ref{sec: connection-general}}
{\noindent \bf Proof of Lemma \ref{lm: tangent-vector-equiva-general}}. Given any tangent vector $\xi = [\U \quad \U_\perp] \begin{bmatrix}
			\S & \D_2^\top\\
			\D_1 & \0
		\end{bmatrix} [\V \quad \V_\perp]^\top$ in $T_{\X}\cM_r$ , denote $\scA_1 =  \{\A = [\A_L^\top \quad\A_R^\top]^\top: \L \A_R^\top + \A_L \R^\top = \xi \}$ and $\scA_2 = \{\A = [\A_L^\top \quad\A_R^\top]^\top:  \L \A_R^\top + \A_L \R^\top = \xi\, \text{ and } \, \L^\top \A_L + \A_L^\top \L- \R^\top \A_R -\A_R^\top \R = \0 \}$. The rest of the proof is divided into two steps: in Step 1 we show the results on  $\scA^\xi_{\L,\R}$; in Step 2 we show the results on $\widetilde{\scA}^{\,\,\xi}_{\L,\R}$.
		
{\bf Step 1}. It is clear $\dim(\scA^\xi_{\L,\R}) = r^2$. For any $\A = [\A_L^\top \quad\A_R^\top]^\top \in \scA^\xi_{\L,\R}$, it is straightforward to check $\L \A_R^\top + \A_L \R^\top = \xi$, so $\scA^\xi_{\L,\R} \subseteq \scA_1$. For any $\A$ such that $\L \A_R^\top + \A_L \R^\top = \xi$, we have
\begin{equation} \label{eq: identity-tangent-vector-equiva-general}
\begin{split}
       \begin{bmatrix}
			\S & \D_2^\top\\
			\D_1 & \0
		\end{bmatrix} =  \begin{bmatrix}
            \U^\top \\
            \U_\perp^\top
    \end{bmatrix} \xi [\V \quad \V_\perp]  
    &=  \begin{bmatrix}
            \U^\top \\
            \U_\perp^\top
    \end{bmatrix} (\L \A_R^\top + \A_L \R^\top) [\V \quad \V_\perp] \\
    & = \begin{bmatrix}
			\P_1 \A_R^\top \V + \U^\top \A_L \P_2^\top  & \P_1 \A_R^\top \V_\perp \\
			\U_\perp^\top \A_L \P_2^\top & \0
		\end{bmatrix}
\end{split}
\end{equation}	by observing $\L = \U \P_1, \R = \V \P_2$. This implies $P_{\U_\perp} \A_L = \U_\perp \D_1 \P_2^{-\top}$, $P_{\V_\perp}\A_R = \V_\perp \D_2 \P_1^{-\top}$ and $\P_1 \A_R^\top \V + \U^\top \A_L \P_2^\top = \S$. By denoting $\S_1 = \U^\top \A_L \P_2^\top$ and $\S_2^\top = \V^\top \A_R \P_1^\top$, we have $\S_1 + \S_2 = \S$, $P_\U \A_L = \U \S_1 \P_2^{-\top}$ and $P_\V \A_R = \V \S_2^\top \P_1^{-\top}$. Finally, $\A_L = P_\U \A_L +  P_{\U_\perp} \A_L = (\U \S_1 + \U_\perp \D_1 )\P_2^{-\top}$, $ \A_R =P_\V \A_R + P_{\V_\perp} \A_R = ( \V \S_2^\top +\V_\perp \D_2 )\P_1^{-\top}$. So $\A = [\A_L^\top \quad\A_R^\top]^\top \in \scA^\xi_{\L,\R}$ and $\scA^\xi_{\L,\R} \supseteq \scA_1$. This proves the first result. 	
		
{\bf Step 2}. Let us begin by proving $\dim(\widetilde{\scA}^{\,\,\xi}_{\L,\R}) = (r^2-r)/2$. First, by simple computation, we have $\dim(\widetilde{\scA}^{\,\,\xi}_{\L,\R}) = \dim(\scS)$ where 
\begin{equation*} 
	 	\scS := \left\{\S_1 \in \bbR^{r \times r}: \P_1^\top \S_1 \P_2^{-\top} + (\P_1^\top \S_1 \P_2^{-\top})^\top + \P_1^{-1} \S_1 \P_2 + (\P_1^{-1} \S_1 \P_2)^\top = \P_1^{-1} \S \P_2 + (\P_1^{-1} \S \P_2)^\top \right\}. 
\end{equation*} Next, we show $\scS$ is of dimension $(r^2 - r)/2$. Construct the following linear map $\varphi_{\L,\R}: \S' \to \P_1^\top \S' \P_2^{-\top} + \P_1^{-1} \S' \P_2.$ We claim $\varphi_{\L,\R}$ is a bijective linear map over $\bbR^{r \times r}$:
\begin{itemize}[leftmargin=*]
    \item {\bf injective part}: suppose there exists $\S_1', \S_2' \in \bbR^{r \times r}$ such that $\S_1' \neq \S_2' $ and $\varphi_{\L,\R}(\S_1') = \varphi_{\L,\R}(\S_2')$. Then by definition of $\varphi_{\L,\R}$, we have $\P_1^\top (\S_1' - \S_2') \P_2^{-\top} + \P_1^{-1} (\S_1' - \S_2') \P_2 = \0$. It further implies $\P_1\P_1^\top (\S_1' - \S_2') + (\S_1' - \S_2') \P_2 \P_2^\top = \0$. This is a Sylvester equation with respect to $(\S_1'- \S_2')$ and we know from \cite[Theorem VII.2.1]{bhatia2013matrix} that it has a unique solution $\0$ due to the fact $\P_1\P_1^\top$ and $-\P_2 \P_2^\top$ have disjoint spectra. So we get $\S_1' = \S_2'$, a contradiction. 
    
    \item {\bf surjective part}: for any $\widetilde{\S} \in \bbR^{r \times r}$, we can find a unique $\widetilde{\S}'$ such that $\varphi_{\L,\R}(\widetilde{\S}') = \widetilde{\S}$. This follows from the facts: (1) $\{\S': \P_1^\top \S' \P_2^{-\top} + \P_1^{-1} \S' \P_2 = \widetilde{\S} \} = \{\S': \P_1 \P_1^\top \S' + \S' \P_2 \P_2^\top = \P_1 \widetilde{\S} \P_2^\top \}$; (2) $\P_1 \P_1^\top \S' + \S' \P_2 \P_2^\top = \P_1 \widetilde{\S} \P_2^\top$ is a Sylvester equation with respect to $\S'$ which has a unique solution again by \cite[Theorem VII.2.1]{bhatia2013matrix}.
\end{itemize}
Then we have $\scS = \{\varphi_{\L,\R}^{-1}(\S'): \S' + \S^{'\top} = \P_1^{-1} \S \P_2 + (\P_1^{-1} \S \P_2)^\top  \}$ and
\begin{equation*}
    \begin{split}
        \dim(\scS) &= \dim(\{\varphi_{\L,\R}^{-1}(\S'): \S' + \S^{'\top} = \P_1^{-1} \S \P_2 + (\P_1^{-1} \S \P_2)^\top  \}) \\
        &= \dim(\{\S': \S' + \S^{'\top} = \P_1^{-1} \S \P_2 + (\P_1^{-1} \S \P_2)^\top  \}) = (r^2-r)/2.
    \end{split}
\end{equation*}

Finally, we show the second result. For any $\A = [\A_L^\top \quad\A_R^\top]^\top \in \widetilde{\scA}^{\,\,\xi}_{\L,\R}$, it is straightforward to check $\L \A_R^\top + \A_L \R^\top = \xi$ and $\L^\top \A_L + \A_L^\top \L- \R^\top \A_R -\A_R^\top \R = \0$. So $\widetilde{\scA}^{\,\,\xi}_{\L,\R} \subseteq \scA_2$. For any $\A \in \scA_2$, following the same proof of \eqref{eq: identity-tangent-vector-equiva-general} we have $\A_L = (\U \S_1 + \U_\perp \D_1 )\P_2^{-\top}, \A_R = ( \V \S_2^\top +\V_\perp \D_2 )\P_1^{-\top}$ where $\S_1 = \U^\top \A_L \P_2^\top$, $\S_2^\top = \V^\top \A_R \P_1^\top$ and they satisfy $\S_1 + \S_2 = \S$. $\L^\top \A_L + \A_L^\top \L- \R^\top \A_R -\A_R^\top \R = \0$ further requires $\S_1,\S_2$ to satisfy $ \P_1^\top \S_1 \P_2^{-\top} + \P_2^{-1} \S_1^\top \P_1 - \P_2^\top \S_2^\top \P_1^{-\top} - \P_1^{-1} \S_2 \P_2 = \0$. So $\A = [\A_L^\top \quad\A_R^\top]^\top \in \widetilde{\scA}^{\,\,\xi}_{\L,\R}$ and $\widetilde{\scA}^{\,\,\xi}_{\L,\R} \supseteq \scA_2$. This finishes the proof of this lemma. \quad $\blacksquare$

\vskip.3cm
{\noindent \bf Proof of Lemma \ref{lm: decomposition-Rp1p2-general}}. We first consider the result of $\scA_{\rnull}^{\L,\R} $ and $ \scA_{\overline{\rnull}}^{\L,\R}$. It is easy to check $\scA_{\rnull}^{\L,\R}$ and $ \scA_{\overline{\rnull}}^{\L,\R}$ are of dimensions $r^2$ and $(p_1+p_2-r)r$, respectively. Since $r^2 + (p_1+p_2-r)r = (p_1 +p_2)r$, to prove $\bbR^{(p_1 + p_2) \times r} = \scA_{\rnull}^{\L,\R} \oplus \scA_{\overline{\rnull}}^{\L,\R}$, we only need to show $\scA_{\rnull}^{\L,\R}$ is orthogonal to $ \scA_{\overline{\rnull}}^{\L,\R}$. Indeed, for any $\A =\begin{bmatrix}
	\U \S \P_2^{-\top} \\
	-\V \S^\top \P_1^{-\top}
\end{bmatrix} \in \scA_{\rnull}^{\L,\R}$, and $\A' = \begin{bmatrix}
	(\U \S' \P_2 \P_2^\top + \U_\perp \D'_1 ) \P_2^{-\top} \\
	(\V \S^{'\top} \P_1 \P_1^\top  + \V_\perp \D'_2 ) \P_1^{-\top}
\end{bmatrix} \in \scA_{\overline{\rnull}}^{\L,\R}$, by simple calculations, we have $\langle \A, \A' \rangle = \langle \S, \S' \rangle - \langle \S, \S' \rangle = 0.$

Next, we prove the result of $\widetilde{\scA}_{\rnull}^{\,\,\L,\R}$ and $\widetilde{\scA}_{\overline{\rnull}}^{\,\,\L,\R}$. From the dimension of $\scS$ in Step 2 of the proof of Lemma \ref{lm: tangent-vector-equiva-general}, we have $\dim(\scS_{\L,\R}) = (r^2 - r)/2$. As a result of this, we  have $\widetilde{\scA}_{\rnull}^{\,\,\L,\R}$ is of dimension $(r^2-r)/2$. Thus, $(\S_1 - \S_2) \perp \scS_{\L,\R}$ in the definition of $\widetilde{\scA}_{\overline{\rnull}}^{\,\,\L,\R}$ adds $(r^2 - r)/2$ constraints and $\dim(\widetilde{\scA}_{\overline{\rnull}}^{\,\,\L,\R}) = (p_1+p_2)r - (r^2 - r)/2$. Now, to prove $\bbR^{(p_1 + p_2) \times r} = \widetilde{\scA}_{\rnull}^{\,\,\L,\R} \oplus \widetilde{\scA}_{\overline{\rnull}}^{\,\,\L,\R}$, we only need to show $\widetilde{\scA}_{\rnull}^{\,\,\L,\R}$ is orthogonal to $\widetilde{\scA}_{\overline{\rnull}}^{\,\,\L,\R}$. In fact, for any $\A = \begin{bmatrix}
	\U \S \P_2^{-\top} \\
	-\V \S^\top \P_1^{-\top}
\end{bmatrix} \in \widetilde{\scA}_{\rnull}^{\,\,\L,\R}$ and $\A' = \begin{bmatrix}
	(\U \S'_1 \P_2 \P_2^\top + \U_\perp \D'_1 ) \P_2^{-\top} \\
	(\V \S_2^{'\top} \P_1 \P_1^\top + \V_\perp \D'_2 ) \P_1^{-\top}
\end{bmatrix} \in \widetilde{\scA}_{\overline{\rnull}}^{\,\,\L,\R}$, we have $	\langle \A, \A' \rangle = \langle \S, \S_1'\rangle - \langle \S, \S_2' \rangle = 0$, where the second equality is because $\S \in \scS_{\L,\R}$ and $(\S_1' - \S_2') \perp \scS_{\L,\R}$ by the construction of $\widetilde{\scA}_{\rnull}^{\,\,\L,\R}$ and $\widetilde{\scA}_{\overline{\rnull}}^{\,\,\L,\R}$, respectively. This finishes the proof of this lemma. \quad $\blacksquare$

\vskip.3cm
{\noindent \bf Proof of Corollary \ref{coro: landscape connection general case}.} First, for any Euclidean FOSP $(\L,\R)$ of \eqref{eq: general factor with reg} or $(\L, \R)$ such that $\L^\top \L = \R^\top \R$, we have $\nabla g_\reg (\L,\R) = \nabla g(\L,\R)$ by \eqref{eq: reg-FOSP-property} and Proposition \ref{prop: gradient-exp}, respectively. The connection on FOSPs of different formulations can be easily obtained by the connection of Riemannian and Euclidean gradients given in \eqref{eq: gradient-connection-general}. Next, we show the equivalence on SOSPs of different formulations.

Suppose $\X$ is a Riemannian SOSP of \eqref{eq: general prob}, we claim any $(\L,\R)$ such that $\L\R^\top = \X$ is a Euclidean SOSP of \eqref{eq: general factor formu} and any $(\L,\R)$ such that $\L\R^\top = \X$ and $\L^\top \L= \R^\top \R$ is a Euclidean SOSP of \eqref{eq: general factor with reg}. To see it, first $(\L,\R)$ in both cases are Euclidean FOSP of \eqref{eq: general factor formu} and \eqref{eq: general factor with reg} as we mentioned before. For any $\A = [\A_L^\top \quad \A_R^\top]^\top \in \bbR^{(p_1 + p_2) \times r}$, by Theorems \ref{th: general case spec-connect 1} and \ref{th: general case spec-connect 2} we have
\begin{equation*}
	\begin{split}
	\nabla^2 g(\L, \R)[\A, \A] & \overset{ \eqref{eq: R-E-Hessian-general-1}} = \Hess f(\X)[\xi_{\L,\R}^\A, \xi_{\L,\R}^\A] \geq 0; \\
		\nabla^2 g_\reg(\L, \R)[\A, \A] & \overset{ \eqref{eq: R-E-Hessian-general-reg} }\geq \Hess f(\X)[ \xi_{\L,\R}^\A, \xi_{\L,\R}^\A ] \geq 0.
	\end{split}
\end{equation*}

Next we show the reverse direction: suppose $(\L,\R)$ is a rank $r$ Euclidean SOSP of \eqref{eq: general factor formu} or \eqref{eq: general factor with reg}, then $\X = \L \R^\top$ is a Riemannian SOSP of \eqref{eq: general prob}. To see this, for any $\xi \in T_\X \cM_r$, 
\begin{equation*}
	\begin{split}
		\Hess f(\L\R^\top)[\xi, \xi] & \overset{ \eqref{eq: R-E-Hessian-general-1} } = \nabla^2 g(\L,\R)[\cL_{\L,\R}^{-1}(\xi),\cL_{\L,\R}^{-1}(\xi)] \geq 0, \\
		\Hess f(\L\R^\top)[\xi, \xi] & \overset{ \eqref{eq: R-E-Hessian-general-reg} } = \nabla^2 g_\reg(\L,\R)[ \cL^{-1}_{\L,\R}(\xi) , \cL^{-1}_{\L,\R}(\xi)] \geq 0.
	\end{split}
\end{equation*} This shows $\X$ is a Riemannian SOSP of \eqref{eq: general prob}.

Suppose $(\L,\R)$ is a rank $r$ Euclidean strict saddle of \eqref{eq: general factor formu} or \eqref{eq: general factor with reg}, and let $\X = \L \R^\top$. Then by definition there exists $\A_1, \A_2 \in \bbR^{(p_1 + p_2) \times r}$ such that $\nabla^2 g(\L, \R)[\A_1, \A_1] < 0$ and $\nabla^2 g_\reg(\L, \R)[\A_2, \A_2] < 0$. Then 
\begin{equation*}
	\begin{split}
	 \Hess f(\X)[\xi_{\L,\R}^{\A_1}, \xi_{\L,\R}^{\A_1}] & \overset{ \eqref{eq: R-E-Hessian-general-1}}=  \nabla^2 g(\L, \R)[\A_1, \A_1] < 0; \\
		\Hess f(\X)[ \xi_{\L,\R}^{\A_2}, \xi_{\L,\R}^{\A_2} ] & \overset{ \eqref{eq: R-E-Hessian-general-reg} }\leq \nabla^2 g_\reg(\L, \R)[\A_2, \A_2]< 0.
	\end{split}
\end{equation*} This implies that $\Hess f(\X)$ has negative eigenvalues in both cases, i.e., $\X$ is a Riemannian strict saddle. The proof for the reverse direction is similar and for simplicity, we omit it here. \quad $\blacksquare$

\vskip.3cm 
{\noindent \bf Proof of Theorem \ref{th: spectrum-connection-two-factorizations}}. 
This proof is divided into two steps. In Step 1, we show \eqref{eq: two-Euclidean-Hessian-connection}; in Step 2, we give the spectrum bounds for the bijective map $\cJ$ and the spectrum connection between $\nabla^2 g_\reg (\L_\reg, \R_\reg)$ and $\nabla^2 g(\L,\R) $. 

{\bf Step 1.} First, since $\L_\reg \R_\reg^\top = \L\R^\top$, $\L_\reg$ and $\L$ share the same left singular subspace. Thus $\L\bDelta = \L\L^\dagger \L_\reg = \L_\reg$ and $\bDelta$ is of rank $r$. Meanwhile, by $\L\R^\top = \L_\reg \R_\reg^\top$, we have $ \bDelta \R_\reg^\top= \L^\dagger \L_\reg \R_\reg^\top = \L^\dagger \L\R^\top = \R^\top $. Moreover, as $(\L_\reg, \R_\reg)$ is a Euclidean FOSP of \eqref{eq: general factor with reg}, by \eqref{eq: reg-Hessian-on-FOSP} we have for any $\A = [\A_L^\top \quad \A_R^\top]^\top  \in \bbR^{(p_1 + p_2) \times r}$:
\begin{equation*}
	\nabla^2 g_\reg (\L_\reg, \R_\reg) [\A,\A] -  \mu \|\L_\reg^\top \A_L + \A_L ^\top \L_\reg - \R_\reg^\top \A_R - \A_R^\top \R_\reg \|_\F^2 = \nabla g^2 (\L_\reg,\R_\reg)[ \A, \A ].
\end{equation*} Next, we show $\nabla g^2 (\L_\reg,\R_\reg)[ \A, \A ]=\nabla g^2 (\L,\R)[ \cJ(\A), \cJ(\A) ]$. By Proposition \ref{prop: Hessian-exp} we have
\begin{equation*}
	\begin{split}
		& \quad \nabla^2 g(\L_\reg,\R_\reg)[\A,\A] \\
		&=  \nabla^2 f(\L_\reg \R_\reg^\top)[\L_\reg \A_R^\top + \A_L \R_\reg^\top,\L_\reg \A_R^\top + \A_L \R_\reg^\top] + 2 \langle \nabla f(\L_\reg \R_\reg  ^\top), \A_L\A_R^\top \rangle\\
		& = \nabla^2 f(\L\R^\top)[\L\bDelta \A_R^\top + \A_L \bDelta^{-1} \R^\top,\L\bDelta \A_R^\top + \A_L \bDelta^{-1} \R^\top] \\
		& \quad + 2 \langle \nabla f(\L\R  ^\top), \A_L \bDelta^{-1} \bDelta \A_R^\top \rangle\\
		& = \nabla^2 g(\L,\R)[\cJ(\A),\cJ(\A)].
	\end{split}
\end{equation*} This finishes the proof for the first part.

{\bf Step 2.} Next, we provide the spectrum bounds for the bijection operator. Suppose $\A =  [\A_L^\top \quad \A_R^\top]^\top$ and $ \cJ(\A) = [\A_L^{'\top} \quad \A_R^{'\top}]^\top $. Then
\begin{equation*}
	\begin{split}
		\|\cJ(\A)\|_\F^2 = \|\A_L'\|_\F^2 + \|\A_R'\|_\F^2 = \|\A_L \bDelta^{-1}\|_\F^2 + \|\A_R \bDelta^\top \|_\F^2 \leq   \left(\sigma_1(\bDelta) \vee (1/\sigma_{r}(\bDelta ))\right)^2 \|\A\|_\F^2,\\
		\|\A\|_\F^2 = \|\A_L\|_\F^2 + \|\A_R\|_\F^2 = \|\A_L' \bDelta\|_\F^2 + \|\A_R' \bDelta^{-\top} \|_\F^2 \leq   \left(\sigma_1(\bDelta) \vee (1/\sigma_{r}(\bDelta ))\right)^2 \|\cJ(\A)\|_\F^2.
	\end{split}
\end{equation*}

Finally, we provide a spectrum connection of two Euclidean Hessians at FOSPs. By \eqref{eq: two-Euclidean-Hessian-connection}, we have $ \nabla^2 g_\reg (\L_\reg, \R_\reg) \succcurlyeq \cJ^* \nabla g^2 (\L,\R) \cJ $. So the first inequality of \eqref{ineq: spectrum-bound-two-Euclidean-Hessian} follows from Lemma \ref{lm: two-diff-symmetric-matrix-spectrum-connection1}(ii) in the Appendix and \eqref{ineq: bijection-spectrum-bound-two-Euclidean-Hessian}. Also by \eqref{eq: reg-FOSP-property}, \eqref{eq: two-Euclidean-Hessian-connection} and Lemma \ref{lm: reg-term-Hessian-bound}, we have $ \nabla^2 g_\reg (\L_\reg, \R_\reg) - 8 \mu \sigma_1(\L_\reg \R_\reg^\top) \mathcal{I} \preccurlyeq \cJ^* \nabla g^2 (\L,\R) \cJ $ and the second inequality in \eqref{ineq: spectrum-bound-two-Euclidean-Hessian} follows from Lemma \ref{lm: two-diff-symmetric-matrix-spectrum-connection1}(i) and \eqref{ineq: bijection-spectrum-bound-two-Euclidean-Hessian}. This finishes the proof.
\quad $\blacksquare$

\subsection{Additional Proofs in Section \ref{sec: application}}
{\noindent \bf Proof of Theorem \ref{th: manifold-landscape-phase-retrieval}.} By Theorem I.1 and Theorem II.2 of \cite{li2019toward}, we have with probability at least $1 - \exp(-C'n)$, the factorization formulation $g(\x)$ in \eqref{eq: modi-formulation-phase-retrieval} has the following geometric landscape properties: (1) $\x^*$ is the unique Euclidean SOSP of $g(\x)$; (2) for any other non-zero Euclidean FOSP $\x$ of $g(\x)$, it satisfies $\lambda_{\min}( \nabla^2 g(\x) ) \leq -3 \|\x^*\|_2^2 = -3 \sigma_1(\X^*)$ under the assumptions of Theorem \ref{th: manifold-landscape-phase-retrieval}.

By Corollary \ref{coro: landscape connection PSD}, we have $\X^* = \x^* \x^{*\top}$ is the unique Riemannian SOSP of \eqref{eq: phase-retrieval-manifold}. In addition, by Theorem \ref{th: RHessian-EHessian PSD}, for any other Riemannian FOSP $\X$ of \eqref{eq: phase-retrieval-manifold}, we have
\begin{equation*}
	\lambda_{\min}(\Hess f(\X)) \leq \frac{1}{4\sigma_1(\X)} \lambda_{\min}( \nabla^2 g(\x) ) \leq - \frac{3 \sigma_1(\X^*) }{4\sigma_1(\X)},
\end{equation*} where $\x$ is any Euclidean FOSP satisfying $\x \x^\top = \X$. \quad $\blacksquare$

\vskip.3cm 
{\noindent \bf Proof of Theorem \ref{th: global-optimality-generic-f-manifold-formulation}}.  
First, \cite{zhu2018global} considered the geometric landscape of \eqref{eq: general factor with reg} when $f$ satisfies the $(2r,4r)$-restricted strong convexity and smoothness property. Under the assumptions of Theorem \ref{th: global-optimality-generic-f-manifold-formulation}, Theorem 3 of \cite{zhu2018global} shows any Euclidean SOSP $(\L,\R)$ of the regularized factorization formulation satisfies $\L\R^\top = \X^*$. By Corollary \ref{coro: landscape connection general case} of this paper, we further conclude if the input rank $r = r^*$ in \eqref{eq: general prob}, then $\X^*$ is the unique Riemannian SOSP of \eqref{eq: general prob} and if $r > r^*$, there is no Riemannian SOSP of \eqref{eq: general prob}.

At the same time, by Theorem 3 of \cite{zhu2018global}, any Euclidean FOSP $(\L,\R)$ of \eqref{eq: general factor with reg} that is not a SOSP must be a strict saddle and satisfy
\begin{equation*}
	\lambda_{\min}( \nabla^2 g_\reg(\L,\R) ) \leq \left\{ \begin{array}{lc}
		-0.08 \alpha_1 \sigma_r(\X^*), & \text{ if } r = r^*;\\
		-0.05 \alpha_1 \cdot ( \sigma_{r^c}^2 (\W) \wedge 2 \sigma_{r^*}(\X^*) ), & \text{ if } r > r^*,
	\end{array} \right.
\end{equation*} where $\W = [\L^\top \quad \R^\top]^\top$ and $r^c$ is the rank of $\W$. Under the manifold formulation \eqref{eq: general prob}, by Theorem \ref{th: general case spec-connect 2}, any Riemannian FOSP $\X$ that is not a Riemannian SOSP must satisfy
\begin{equation*}
\begin{split}
	\lambda_{\min}( \Hess f(\X) )&\leq \lambda_{\min}( \nabla^2 g_\reg(\L',\R') )/2\sigma_1(\X) \\
	&\leq \left\{ \begin{array}{lc}
		-0.08 \alpha_1 \sigma_r(\X^*)/(2\sigma_1(\X)), & \text{ if } r = r^*;\\
		-0.05 \alpha_1 \cdot ( \sigma_{r}^2 (\W') \wedge 2 \sigma_{r^*}(\X^*) )/(2\sigma_1(\X)), & \text{ if } r > r^*,
	\end{array} \right.
\end{split}
\end{equation*} where $\W' = [\L^{'\top} \quad \R^{'\top}]^\top$ and $(\L', \R')$ is a rank $r$ Euclidean FOSP of \eqref{eq: general factor with reg} satisfying $\L'\R^{'\top} = \X$. Finally, we only need to compute $\sigma_{r}^2 (\W')$. By Lemma \ref{lm: balanced-factor-property} we have $\L'= \U \P$ and $\R' = \V \P$ for some invertible $\P \in \bbR^{r \times r}$, where $\U, \V$ are the left and right singular subspaces of $\X$. So $\sigma_{r} (\W') = \sigma_r([\L^{'\top} \quad \R^{'\top}]^\top) = \sqrt{2} \sigma_r(\P) = \sqrt{2 \sigma_r(\X)} $. This finishes the proof of this theorem.
\quad $\blacksquare$

\vskip.3cm
{\noindent \bf Proof of Theorem \ref{th: global-optimality-factorization-noregularization}}. Under the assumptions of Theorem \ref{th: global-optimality-factorization-noregularization}, by Theorem 3 of \cite{zhu2018global} we have for a rank $r$ Euclidean FOSP $(\L_\reg,\R_\reg)$ of the regularized formulation \eqref{eq: general factor with reg}, it is either a Euclidean SOSP satisfying $\L_\reg \R_\reg^\top = \X^*$ or a strict saddle with $\lambda_{\min}( \nabla^2 g_\reg(\L_\reg,\R_\reg) ) \leq -0.08 \alpha_1\sigma_r(\X^*)$.

By Corollary \ref{coro: landscape connection general case} and Theorem \ref{th: spectrum-connection-two-factorizations}, we have for any rank $r$ Euclidean FOSP $(\L,\R)$ of \eqref{eq: general factor formu}, it is either a Euclidean SOSP satisfying $\L\R^\top = \X^*$ or a strict saddle with $$\lambda_{\min}( \nabla^2 g(\L,\R) ) \leq \theta_{\bDelta}^2 \lambda_{\min}( \nabla^2 g_\reg(\L'_\reg,\R'_\reg) ) \leq -0.08 \theta_{\bDelta}^2 \alpha_1\sigma_r(\X^*),$$ where $\theta_{\bDelta}:= (1/\sigma_1(\bDelta)) \wedge \sigma_{r}(\bDelta)$, $\bDelta = \L^\dagger \L'_\reg$ and $(\L'_\reg,\R'_\reg)$ is a rank $r$ Euclidean FOSP of \eqref{eq: general factor with reg} satisfying $\L'_\reg \R_\reg^{'\top} = \L\R^\top =: \X$.

Finally, we give a lower bound for $\theta_{\bDelta}$. Notice $\L\bDelta = \L'_\reg$, and
\begin{equation*}
\begin{split}
	&\sigma_1(\bDelta) = \sigma_1(\L^\dagger \L'_\reg) \leq \sigma_1(\L^\dagger) \sigma_1(\L'_\reg) \overset{ \eqref{eq: reg-FOSP-property}, \text{Lemma } \ref{lm: balanced-factor-property} }= \sigma_1^{1/2}(\X)/\sigma_r(\L), \\
	&\sigma_r^{1/2}(\X)\overset{ \eqref{eq: reg-FOSP-property}, \text{Lemma } \ref{lm: balanced-factor-property} }=\sigma_r(\L'_\reg) = \sigma_r(\L\bDelta) = \inf_{\x:\|\x\|_2 = 1} \|\L\bDelta \x \|_2 \leq \sigma_1(\L) \inf_{\x:\|\x\|_2 = 1} \|\bDelta \x \|_2 = \sigma_1(\L) \sigma_r(\bDelta ) .
\end{split}
\end{equation*} We have $\theta_{\bDelta}:= (1/\sigma_1(\bDelta)) \wedge \sigma_{r}(\bDelta) \geq (\sigma_r(\L)/\sigma^{1/2}_1(\X)) \wedge ( \sigma^{1/2}_r(\X)/\sigma_1(\L))$. This finishes the proof of this theorem.
\quad $\blacksquare$

\section{Additional Lemmas} \label{sec: additional-lemmas}
Recall $\lambda_k(\cdot)$ and $\sigma_k(\cdot)$ are the $k$th largest eigenvalue and $k$th largest singular value of matrix $(\cdot)$. Also $\lambda_{\max}(\cdot)$, $\lambda_{\min}(\cdot)$ denote the largest and least eigenvalue of matrix $(\cdot)$.
\begin{Lemma}\label{lm: two-symmetric-matrix-spectrum-connection}
	Suppose $\A \in \bbS^{p \times p}$ is symmetric and $\P \in \bbR^{p \times p}$ is invertible. Then $\lambda_k(\P^\top \A \P)$ is sandwiched between $\sigma_p^2(\P) \lambda_k(\A )$ and $\sigma_1^2(\P) \lambda_k(\A )$ for $k = 1,\ldots, p$.
\end{Lemma}
{\noindent \bf Proof.} Suppose $\u_1,\ldots,\u_p$ are eigenvectors corresponding to $\lambda_1(\A),\ldots,\lambda_p(\A)$ and $\v_1,\ldots,\v_p$ are eigenvectors corresponding to $\lambda_1(\P^\top\A\P),\ldots,\lambda_p(\P^\top\A\P)$. For $k = 1,\ldots, p$, define
\begin{equation*}
	\begin{split}
		&\mathcal{U}_k = \rspan\{\u_1,\ldots,\u_k\}, \quad \mathcal{U}'_k = \rspan\{\P^{-1}\u_1,\ldots,\P^{-1}\u_k\},\\
		&\mathcal{V}_k = \rspan\{\v_1,\ldots,\v_k\}, \quad \mathcal{V}'_k = \rspan\{\P\v_1,\ldots,\P\v_k\}.
	\end{split}
\end{equation*}

Let us first consider the case that $\lambda_k(\A) \geq 0$. By Lemma \ref{lm: max-min-theorem}, we have
\begin{equation} \label{ineq: spectrum-ineq1}
	\begin{split}
		\lambda_k(\P^\top \A \P) \geq \min_{\u \in \cU_k',\u \neq \0 } \frac{\u^\top \P^\top \A \P \u}{\|\u\|_2^2} = \min_{\u \in \cU_k,\u \neq \0 } \frac{\u^\top \A  \u}{\|\P^{-1}\u\|_2^2} \geq  \min_{\u \in \cU_k,\u \neq \0 }\frac{\lambda_k(\A)\|\u\|_2^2}{\|\P^{-1}\u\|_2^2} \geq \lambda_k(\A) \sigma^2_p(\P)\geq 0.
	\end{split}
\end{equation} On the other hand, we have
\begin{equation}\label{ineq: spectrum-ineq2}
	\begin{split}
		\lambda_k(\A) \overset{\text{Lemma } \ref{lm: max-min-theorem}}\geq \min_{\u \in \cV_k',\u \neq \0 } \frac{\u^\top \P^{-\top} \P^\top \A \P \P^{-1} \u}{\|\u\|_2^2} = \min_{\v \in \cV_k,\v \neq \0 } \frac{\v^\top \P^\top \A \P \v}{\|\P\v\|_2^2} &\geq \min_{\v \in \cV_k,\v \neq \0 } \frac{\lambda_k(\P^\top \A \P) \|\v\|_2^2 }{\|\P\v\|_2^2}\\
		 &\overset{\eqref{ineq: spectrum-ineq1} }\geq \frac{\lambda_k(\P^\top \A \P) }{ \sigma^2_1(\P)}.
	\end{split}
\end{equation} So we have proved the result for the case that $\lambda_k(\A) \geq 0$. When $\lambda_k(\A) < 0$, we have $\lambda_{p+1-k}(-\A) = -\lambda_k(\A) > 0$. Following the same proof of \eqref{ineq: spectrum-ineq1} and \eqref{ineq: spectrum-ineq2}, we have
\begin{equation*}
	\begin{split}
		-\lambda_k(\P^\top \A \P) = \lambda_{p+1-k}(-\P^\top \A \P) \geq \sigma_p^2(\P) \lambda_{p+1-k}(-\A) = -\sigma_p^2(\P)\lambda_k(\A) > 0,\\
		-\lambda_k(\A) =  \lambda_{p+1-k}(-\A) \geq  \lambda_{p+1-k}(-\P^\top \A \P)/\sigma^2_1(\P) = -\lambda_k(\P^\top \A \P)/\sigma^2_1(\P).
	\end{split}
\end{equation*} This finishes the proof of this lemma. \quad $\blacksquare$

\begin{Lemma}\label{lm: two-diff-symmetric-matrix-spectrum-connection1}
	Suppose $\A \in \bbS^{p \times p}, \B \in \bbS^{q \times q}$ are symmetric matrices with $q \geq p$ and $\P \in \bbR^{q \times p}, \Q \in \bbR^{p \times q}$.
	\begin{itemize}
		\item[(i)] If $\P^\top \B \P \succcurlyeq \A $, then $\lambda_k(\B)\sigma^2_1(\P) \vee \lambda_k(\B)\sigma^2_p(\P) \geq \lambda_k(\A)$ holds for $k = 1,\ldots,p$.
		\item[(ii)] If $\P^\top \B \P \preccurlyeq \A$, then $\lambda_{k+q-p}(\B)\sigma_1^2(\P)\wedge \lambda_{k+q-p}(\B)\sigma_p^2(\P) \leq \lambda_{k}(\A)$ holds for $k = 1,\ldots,p$.
		\item[(iii)] If $\Q^\top \A \Q \preccurlyeq \B $, then $\lambda_{\min}(\B) \geq  \sigma^2_1(\Q) \lambda_{\min}(\A) \wedge 0 $.
		\item[(iv)] If $\Q^\top \A \Q \succcurlyeq \B $, then $\lambda_{1}(\B) \leq  \sigma^2_1(\Q) \lambda_{\max}(\A) \vee 0$.
	\end{itemize}
\end{Lemma}
{\noindent \bf Proof.} We first prove the first and the second claims under the assumption that $\sigma_p(\P)>0$, i.e., all $p$ columns of $\P$ are linearly independent. 

Suppose $\u_1,\ldots,\u_p$ are eigenvectors corresponding to $\lambda_1(\A),\ldots,\lambda_p(\A)$, respectively and let $\mathcal{U}_k = \rspan\{\u_1,\ldots,\u_k\}$. Then
\begin{equation*}
	\lambda_k(\B) \overset{(a) }\geq \inf_{\u \in \mathcal{U}_k } \frac{\u^\top \P^\top \B \P \u }{\|\P \u\|_2^2} \geq  \inf_{\u \in \mathcal{U}_k }  \frac{\u^\top \A \u }{\|\P \u\|_2^2} \geq \inf_{\u \in \mathcal{U}_k } \frac{\lambda_k(\A) \|\u\|_2^2}{\|\P \u\|_2^2} \geq \left\{ \begin{array}{c c}
		\lambda_k(\A)/\sigma^2_1(\P), & \text{ if } \lambda_k (\A) \geq 0;\\
		\lambda_k(\A)/\sigma^2_p(\P), & \text{ if } \lambda_k(\A) < 0.
	\end{array} \right.
\end{equation*} Here (a) is because $\{ \P \u_1,\ldots, \P \u_k \}$ forms a $k$ dimensional subspace in $\bbR^{q}$ and Lemma \ref{lm: max-min-theorem}. 

To see the second claim under $\sigma_p(\P)>0$, suppose $\v_1,\ldots,\v_q$ are eigenvectors corresponding to $\lambda_1(\B),\ldots,\lambda_q(\B)$ and let $\mathcal{V}_{k+q-p} = \rspan\{\v_1,\ldots,\v_{k + q-p}\}$. 
\begin{equation} \label{ineq: two-diff-size-spectrum-connection-ineq1}
	\begin{split}
		\lambda_k(\A) \overset{(a)}\geq &  \inf_{\v:\P \v \in \cV_{k+q-p}} \frac{\v^\top \A \v}{\|\v\|_2^2} \geq \inf_{\v:\P \v \in \cV_{k+q-p}} \frac{\v^\top \P^\top \B \P \v}{\|\v\|_2^2} \geq \inf_{\v:\P \v \in \cV_{k+q-p}} \frac{\lambda_{k+q-p}(\B) \|\P\v\|_2^2}{\| \v\|_2^2} \\
		\geq & \left\{ \begin{array}{c c}
		\sigma^2_p(\P)\lambda_{k+q-p}(\B), & \text{ if } \lambda_{k+q-p} (\B) \geq 0\\
		\sigma^2_1(\P)\lambda_{k+q-p}(\B), & \text{ if } \lambda_{k+q-p}(\B) < 0
	\end{array} \right.
	\end{split}
\end{equation} Here (a) is because of Lemma \ref{lm: max-min-theorem} and the fact $\{\v: \P \v \in \cV_{k+q-p} \}$ has dimension at least $k$. 

When $\sigma_p(\P)=0$, we construct a series of matrices $\P_l$ such that $\lim_{l\to\infty}\P_l = \P$ and $\sigma_p(\P_l)>0$. According to the previous proofs, 
$$\lambda_k(\B)\sigma^2_1(\P_l) \vee \lambda_k(\B)\sigma^2_p(\P_l) \geq \lambda_k(\P_l^\top \B \P_l),$$
$$\lambda_{k+q-p}(\B)\sigma_1^2(\P_l)\wedge \lambda_{k+q-p}(\B)\sigma_p^2(\P_l) \leq \lambda_{k}(\P_l^\top \B \P_l).$$
Since $\sigma_k(\cdot)$ and $\lambda_k(\cdot)$ are continuous functions of the input matrix, by taking $l\to \infty$, we have
$$\lambda_k(\B)\sigma^2_1(\P) \vee \lambda_k(\B)\sigma^2_p(\P) \geq \lambda_k(\P^\top \B \P) \overset{(a)}\geq \lambda_k(\A), \text{ under the assumption of Claim 1;}$$
$$\lambda_{k+q-p}(\B)\sigma_1^2(\P)\wedge \lambda_{k+q-p}(\B)\sigma_p^2(\P) \leq \lambda_{k}(\P^\top \B \P) \overset{(a)}\leq \lambda_k(\A), \text{ under the assumption of Claim 2.}$$ Here in (a) we use the fact for any two $p_1$-by-$p_1$ symmetric matrices $\W_1,\W_2$, $\W_1 \succcurlyeq \W_2 $ implies $\lambda_k(\W_1) \geq \lambda_k(\W_2)$ for any $k \in [p_1] $. This finishes the proof for the first two claims.

To prove the third claim, suppose $\v_{\min}$ is the eigenvector corresponding to the smallest eigenvalue of $\B$, then
\begin{equation*}
	\lambda_{\min}(\B)  = \v_{\min}^\top \B \v_{\min} \geq  \v_{\min}^\top \Q^\top \A \Q \v_{\min} \geq \lambda_{\min}(\A) \| \Q \v_{\min}\|_2^2 \geq \left\{ \begin{array}{c c}
		0, & \text{ if } \lambda_{\min} (\A) \geq 0;\\
		\sigma^2_1(\Q)\lambda_{\min}(\A),& \text{ if } \lambda_{\min}(\A) < 0.
	\end{array} \right.
\end{equation*}

To prove the last claim, suppose $\v_{\max}$ is the eigenvector corresponding to the largest eigenvalue of $\B$, then
\begin{equation*}
	\lambda_{1}(\B)  = \v_{\max}^\top \B \v_{\max} \leq  \v_{\max}^\top \Q^\top \A \Q \v_{\max} \leq \lambda_{\max}(\A) \| \Q \v_{\max}\|_2^2 \leq \left\{ \begin{array}{c c}
		0, & \text{ if } \lambda_{\max} (\A) < 0;\\
		\sigma^2_1(\Q)\lambda_{\max}(\A),& \text{ if } \lambda_{\max}(\A) \geq 0.
	\end{array} \right.
\end{equation*}

This finishes the proof of this lemma.
\quad $\blacksquare$

\begin{Lemma}{\rm(Max-min Theorem for Eigenvalues \cite[Corollary III.1.2]{bhatia2013matrix} )}\label{lm: max-min-theorem}
For any $p$-by-$p$ real symmetric matrix $\A$ with eigenvalues $\lambda_1 \geq \lambda_2 \geq \cdots \geq \lambda_p$. If $\scC_k$ denotes the set of subspaces of $\bbR^p$ of dimension $k$, then $\lambda_k = \max_{C \in \scC_k} \min_{\u \in C, \u \neq \0} \u^\top \A \u /\|\u\|_2^2.$
\end{Lemma}

\begin{Lemma}\label{lm: reg-term-Hessian-bound}
	Suppose $\L\in \bbR^{p_1 \times r}$ and $\R \in \bbR^{p_2\times r}$. Then for any $[\A_L^\top \quad \A_R^\top]^\top \in \bbR^{(p_1+p_2) \times r}$, 
	\begin{equation*}
		\|\L^\top \A_L + \A_L^\top \L- \R^\top \A_R - \A_R^\top \L\|_\F^2 \leq 8(\sigma_1(\L) \vee \sigma_1(\R) )^2 (\|\A_R\|_\F^2 + \|\A_L\|_\F^2).
	\end{equation*}
\end{Lemma}
{\bf \noindent Proof.}
\begin{equation*}
	\begin{split}
		\|\L^\top \A_L + \A_L^\top \L- \R^\top \A_R - \A_R^\top \L\|_\F^2  &\leq 2(\|\L^\top \A_L + \A_L^\top \L\|_\F^2 + \|\R^\top \A_R + \A_R^\top \L\|_\F^2)\\
		& \leq 2( 4\|\L^\top \A_L\|_\F^2 + 4\|\R^\top \A_R\|_\F^2  )\\
		& \leq 8(\sigma_1(\L) \vee \sigma_1(\R) )^2 (\|\A_R\|_\F^2 + \|\A_L\|_\F^2).
	\end{split}
\end{equation*}
This finishes the proof. \quad $\blacksquare$

\begin{Lemma} \label{lm: balanced-factor-property}
	Suppose $\L\in \bbR^{p_1 \times r}, \R \in \bbR^{p_2 \times r}$ are two rank $r$ matrices and $\L^\top \L= \R^\top \R$. Let $\U \bSigma \V^\top$ be a SVD of $\L\R^\top$. Then we have $\L= \U \P, \R = \V \P$ for some $r$-by-$r$ full rank matrix $\P$ satisfying $\P \P^\top = \bSigma$.
\end{Lemma}
{\noindent \bf Proof}. First since $\L\R^\top$ has SVD $\U \bSigma \V^\top$, we have $\L= \U \P_1$ and $\R = \V \P_2$. Next we show $\P_1 = \P_2$. Since $\P_1 \P_2^\top = \bSigma$, we have
\begin{equation*} 
	\bSigma^2 = \P_1 \P_2^\top \P_2 \P_1^\top \overset{(a)}= \P_1 \P_1^\top \P_1 \P_1^\top \overset{(b)}\Longrightarrow \bSigma = \P_1 \P_1^\top.
\end{equation*} Here (a) is because $\L^\top \L= \R^\top \R$ implies $\P_1^\top \P_1 = \P_2^\top \P_2$; and (b) is because a PSD matrix has a unique principal square root \citep{johnson2001uniqueness}. This finishes the proof of this lemma.
\quad $\blacksquare$

\end{document}